\documentclass[12pt,oneside,english]{amsart}
\usepackage[T1]{fontenc}
\usepackage[latin9]{inputenc}
\usepackage{geometry}
\usepackage{array}
\usepackage{longtable}
\usepackage{float}
\usepackage{amstext}
\usepackage{amsthm}
\usepackage{amssymb}
\usepackage{stmaryrd}
\usepackage{graphicx}
\usepackage{setspace}
\usepackage{wasysym}

\makeatletter

\providecommand{\tabularnewline}{\\}

\numberwithin{equation}{section}
\numberwithin{figure}{section}
\numberwithin{table}{section}
\theoremstyle{plain}
\newtheorem{thm}{\protect\theoremname}
\theoremstyle{definition}
\newtheorem{defn}[thm]{\protect\definitionname}
\theoremstyle{plain}
\newtheorem{lem}[thm]{\protect\lemmaname}

\makeatother

\usepackage{babel}
\providecommand{\definitionname}{Definition}
\providecommand{\lemmaname}{Lemma}
\providecommand{\theoremname}{Theorem}

\begin{document}
\title[Centres of centralizers of nilpotent elements in Lie superalgebras]{Centres of centralizers of nilpotent elements in exceptional Lie superalgebras}
\author{Leyu Han}
\begin{abstract}
Let $\mathfrak{g}=\mathfrak{g}_{\bar{0}}\oplus\mathfrak{g}_{\bar{1}}$
be a finite-dimensional simple Lie superalgebra of type $D(2,1;\alpha)$,
$G(3)$ or $F(4)$ over $\mathbb{C}$. Let $G$ be the simply connected
semisimple algebraic group over $\mathbb{C}$ such that $\mathrm{Lie}(G)=\mathfrak{g}_{\bar{0}}$.
Suppose $e\in\mathfrak{g}_{\bar{0}}$ is nilpotent. We describe the
centralizer $\mathfrak{g}^{e}$ of $e$ in{\normalsize{} }$\mathfrak{g}$
and its centre $\mathfrak{z}(\mathfrak{g}^{e})$ especially. We also
determine the labelled Dynkin diagram for $e$. We prove theorems
relating the dimension of $\left(\mathfrak{z}(\mathfrak{g}^{e})\right)^{G^{e}}$
and the labelled Dynkin diagram.
\end{abstract}

\maketitle
\begin{singlespace}

\section{Introduction}
\end{singlespace}

\begin{singlespace}
\noindent Let $\mathfrak{g}=\mathfrak{g}_{\bar{0}}\oplus\mathfrak{g}_{\bar{1}}$
be a simple basic classical Lie superalgebra over $\mathbb{C}$ and
let $e\in\mathfrak{g}_{\bar{0}}$ be nilpotent. This paper forms a
project to investigate the centralizer $\mathfrak{g}^{e}=\{x\in\mathfrak{g}:[x,e]=0\}$
of $e$ in $\mathfrak{g}$ and the centre of centralizer $\mathfrak{z}(\mathfrak{g}^{e})=\{x\in\mathfrak{g}^{e}:[x,y]=0\text{ for all }y\in\mathfrak{g}^{e}\}$
of $e$ in $\mathfrak{g}$. There has been a lot of research on the
centralizer and the centre of centralizer of nilpotent element in
the case of simple Lie algebras, but a lot less is known in the case
of Lie superalgebras. The aim of this paper is to lift the level of
understanding in the area of Lie superalgebra to the similar level
to that in the area of Lie algebra. More precisely, the present paper
is planned to be the first of two papers in which we calculate bases
of $\mathfrak{g}^{e}$ and $\mathfrak{z}(\mathfrak{g}^{e})$. In this
paper, we consider the case when $\mathfrak{g}$ is a exceptional
Lie superalgebra, while the other dealing with the other basic classical
Lie superalgebras $\mathfrak{sl}(m|n)$ and $\mathfrak{osp}(m|2n)$.
\end{singlespace}

\begin{singlespace}
For nonsuper cases, work on $\mathfrak{g}^{e}$ dates back to 1966,
when Springer \cite{Springer1966} worked with a simple algebraic
group $G$ and considered the centralizer $G^{u}$ of a unipotent
element $u\in G$. Many mathematicians studied $G^{u}$ for different
types of $G$ after that, the reader is referred to the introduction
of \cite{Lawther2008} for an overview of the other research of $G^{e}$
and $\mathfrak{g}^{e}$. Jantzen gave an explicit account on the structure
of $\mathfrak{g}^{e}$ for classical Lie algebras $\mathfrak{g}$
in \cite{Jantzen2004a}. Seitz \cite{Seitz2004} pointed out the dimension
of $Z(G^{u})$ is of considerable interest. In \cite{Lawther2008},
Lawther\textendash Testerman studied the centralizer $G^{u}$, especially
its centre $Z(G^{u})$ and determined the dimension of the Lie algebra
of $Z(G^{u})$ over a field of characteristic $0$ or a good prime.
Using a $G$-equivariant homeomorphism, Lawther\textendash Testerman
worked with a nilpotent element $e\in\mathrm{Lie}(G)$ rather than
$u$. The study of the centre $\mathfrak{z}(\mathfrak{g}^{e})$ for
classical Lie algebras $\mathfrak{g}$ over a field of characteristic
$0$ was undertaken by Yakimova in \cite{Yakimova2009} and Lawther\textendash Testerman
\cite{Lawther2008} made use of work of Yakimova in \cite{Yakimova2009}
to deal with classical cases. 
\end{singlespace}

Centralizers of nilpotent elements $e$ in Lie superalgebras $\mathfrak{g}$
for the case where $\mathfrak{g}=\mathfrak{gl}(m|n)$ was done in
\cite{Wang2009} over a field of prime characteristic. In \cite{Hoyt2012},
Hoyt claimed that the construction is identical in characteristic
zero and further describe the construction of $\mathfrak{g}^{e}$
for $\mathfrak{g}=\mathfrak{osp}(m|2n)$. The dimension of $\mathfrak{z}(\mathfrak{g}^{e})$
for exceptional Lie superalgebras remains a mystery and we attempt
to shed some light upon this mystery here.

In the rest of this introduction, we give a more detailed survey of
our results.

\begin{singlespace}
In the present paper, let $\mathfrak{g}=\mathfrak{g}_{\bar{0}}\oplus\mathfrak{g}_{\bar{1}}$
be a finite-dimensional simple Lie superalgebra of type $D(2,1;\alpha)$,
$G(3)$ or $F(4)$ over $\mathbb{C}$. Let $G$ be the simply connected
semisimple algebraic group over $\mathbb{C}$ such that $\mathrm{Lie}(G)=\mathfrak{g}_{\bar{0}}$.
Then there is a representation $\rho:G\rightarrow\mathrm{GL}(\mathfrak{g}_{\bar{1}})$
such that $d_{\rho}:\mathrm{Lie}(G)\rightarrow\mathfrak{gl}(\mathfrak{g}_{\bar{1}})$
determines the adjoint action of $\mathfrak{g}_{\bar{0}}$ on $\mathfrak{g}_{\bar{1}}$.
Note that $G$ is given explicitly in Table \ref{tab:G}. 
\end{singlespace}

\begin{table}[H]
\begin{singlespace}
\noindent \begin{centering}
\begin{tabular}{|c||c|}
\hline 
$\mathfrak{g}$ & $G$\tabularnewline
\hline 
\hline 
$D(2,1;\alpha)$ & $\mathrm{SL}_{2}(\mathbb{C})\times\mathrm{SL}_{2}(\mathbb{C})\times\mathrm{SL}_{2}(\mathbb{C})$\tabularnewline
\hline 
\hline 
$G(3)$ & $\mathrm{SL}_{2}(\mathbb{C})\times G_{2}$\tabularnewline
\hline 
\hline 
$F(4)$ & $\mathrm{SL}_{2}(\mathbb{C})\times\mathrm{Spin}_{7}(\mathbb{C})$\tabularnewline
\hline 
\end{tabular}
\par\end{centering}
\end{singlespace}
\caption{\label{tab:G}Algebraic group $G$}
\end{table}

Let $e\in\mathfrak{g}_{\bar{0}}$ be nilpotent, we investigate the
centralizer $\mathfrak{g}^{e}$ of $e$ in $\mathfrak{g}$, especially
its centre $\mathfrak{z}(\mathfrak{g}^{e})$. In particular, we give
bases for $\mathfrak{g}^{e}$ and $\mathfrak{z}(\mathfrak{g}^{e})$
in Tables \ref{tab:results in D(2,1)}, \ref{tab:G(3)} and \ref{tab:F(4)}.
Write $G^{e}=\{g\in G:geg^{-1}=e\}$ for the centralizer of $e$ in
$G$. We also find a basis for $\left(\mathfrak{z}(\mathfrak{g}^{e})\right)^{G^{e}}=\{x\in\mathfrak{z}(\mathfrak{g}^{e}):gxg^{-1}=x\text{ for all }g\in G^{e}\}$
in \S5.5 and \S6.7. 

Note that $e$ lies in an $\mathfrak{sl}(2)$-triple $\{e,h,f\}\subseteq\mathfrak{g}_{\bar{0}}$
by Jacobson\textendash Morozov Theorem. We use $h$ to determine labelled
Dynkin diagrams with respect to $e$. A full definition of labelled
Dynkin diagrams with respect to $e$ is given in \S2.3. In contrary
to Lie algebra case, in general $e$ determines more than one labelled
Dynkin diagram, dependent on a choice of positive roots. 

Write $\mathfrak{g}=\bigoplus_{j\in\mathbb{Z}}\mathfrak{g}(j)$ as
its ad$h$-eigenspace decomposition, we can decompose $\mathfrak{g}^{e}$
into the direct sum of ad$h$-eigenspaces, i.e. $\mathfrak{g}^{e}=\bigoplus_{j\geq0}\mathfrak{g}^{e}(j)$.
We also describe the $\mathfrak{g}^{e}(0)$-module structure on each
$\mathfrak{g}^{e}(j)$ for $j>0$ in Tables \ref{tab:g^e(0)-D(2,1;)},
\ref{tab:g^e(0)-G(3)} and \ref{tab:g^e(0)-F(4)}.

\begin{singlespace}
Our results can therefore be viewed as extensions of those obtained
by Lawther and Testerman in \cite{Lawther2008} to case of Lie superalgebras
over a field of characteristic zero. They obtain four theorems as
a consequence of their work. In this paper, we obtain analogues of
Theorems 2\textendash 4 in \cite{Lawther2008} for exceptional Lie
superalgebras. We view $\left(\mathfrak{z}(\mathfrak{g}^{e})\right)^{G^{e}}$
as the correct replacement for $Z(G^{u})$ since $\mathrm{Lie}(Z(G^{e}))=\left(\mathfrak{z}(\mathfrak{g}^{e})\right)^{G^{e}}$
for a field of characteristic zero.

Fix $\varDelta$ to be a labelled Dynkin diagram with respect to $e$.
Let $n_{i}(\varDelta)$ be the number of nodes which have labels equal
to $i$ in $\varDelta$. An interesting fact is that the choice of
$\varDelta$ does not affect the following theorems and labels in
$\varDelta$ can only be $0,1$ or $2$. 
\end{singlespace}

We first consider the case where $\varDelta$ only have even labels.
\begin{thm}
\begin{singlespace}
\noindent Let $\mathfrak{g}=\mathfrak{g}_{\bar{0}}\oplus\mathfrak{g}_{\bar{1}}$
be a Lie superalgebra of type $D(2,1;\alpha)$, $G(3)$ or $F(4)$
and $e\in\mathfrak{g}_{\bar{0}}$ be nilpotent. Let $G$ be the algebraic
group defined as in Table \ref{tab:G}. Assume $n_{1}(\varDelta)=0$.
Then $\dim\left(\mathfrak{z}(\mathfrak{g}^{e})\right)^{G^{e}}=n_{2}(\varDelta)=\dim\mathfrak{z}(\mathfrak{g}^{h})$.
\end{singlespace}
\end{thm}

\begin{singlespace}
\noindent In order to state Theorem 2, we define the sub-labelled
Dynkin diagram $\varDelta_{0}$ to be\textit{ the $2$-free core of
$\varDelta$} where $\varDelta_{0}$ is obtained by removing all nodes
with labels equal to $2$ from $\varDelta$. 
\end{singlespace}
\begin{thm}
\begin{singlespace}
Let $\mathfrak{g}=\mathfrak{g}_{\bar{0}}\oplus\mathfrak{g}_{\bar{1}}$
be a Lie superalgebra of type $D(2,1;\alpha)$, $G(3)$ or $F(4)$
and $e\in\mathfrak{g}_{\bar{0}}$ be nilpotent. Let $G$ be the algebraic
group defined as in Table \ref{tab:G}. Let $\varDelta_{0}$ be the
$2$-free core of $\varDelta$. Let $\mathfrak{g}_{0}$ be the subalgebra
of $\mathfrak{g}$ generated by the root vectors corresponding to
the simple roots in $\varDelta_{0}$. Then $\mathfrak{g}_{0}$ is
a direct sum of simple Lie superalgebras and there exists a nilpotent
orbit in $(\mathfrak{g}_{0})_{\bar{0}}$ having labelled Dynkin diagram
$\varDelta_{0}$. Let $G_{0}$ be the simply connected semisimple
algebraic group with respect to $(\mathfrak{g}_{0})_{\bar{0}}$. Suppose
$e_{0}\in(\mathfrak{g}_{0})_{\bar{0}}$ is a representative of this
orbit, then:

1. $\dim\mathfrak{g}^{e}-\dim\mathfrak{g}_{0}^{e_{0}}=n_{2}(\varDelta)$;

2. $\dim\left(\mathfrak{z}(\mathfrak{g}^{e})\right)^{G^{e}}-\dim\left(\mathfrak{z}(\mathfrak{g}_{0}^{e_{0}})\right)^{G_{0}^{e_{0}}}=n_{2}(\varDelta)$.
\end{singlespace}
\end{thm}

\begin{singlespace}
Our last result gives a more general result relating $\dim\left(\mathfrak{z}(\mathfrak{g}^{e})\right)^{G^{e}}$
and $\varDelta$. In this statement we use the notation for nilpotent
elements as introduced later in \S4.1, \S5.1 and \S6.1.
\end{singlespace}
\begin{thm}
\begin{singlespace}
\noindent Let $\mathfrak{g}=\mathfrak{g}_{\bar{0}}\oplus\mathfrak{g}_{\bar{1}}$
be a Lie superalgebra of type $D(2,1;\alpha)$, $G(3)$ or $F(4)$
and $e\in\mathfrak{g}_{\bar{0}}$ be nilpotent. Let $a_{1},...,a_{l}$
be the labels in $\varDelta$. Then
\begin{equation}
\dim\left(\mathfrak{z}(\mathfrak{g}^{e})\right)^{G^{e}}=\left\lceil \frac{1}{2}\sum_{i=1}^{l}a_{i}\right\rceil +\varepsilon\label{eq:15 theorem 3}
\end{equation}

\noindent where the value of $\varepsilon$ is equal to $0$ with
the following exceptions: $\varepsilon=-1$ when $\mathfrak{g}=D(2,1;\alpha)$,
$e=E_{1}+E_{2}+E_{3}$ or $\mathfrak{g}=F(4)$, $e=E+e_{(7)}$.
\end{singlespace}
\end{thm}

This paper is organized as follows: In Section \ref{sec:preliminaries},
we recall some basic vocabulary of Lie superalgebras such as basic
classical Lie superalgebras, root systems and Borel subalgebras. We
also give a full definition of the labelled Dynkin diagram for a system
of simple roots. In Section \ref{sec:general-method}, we recall some
generalities on $\mathfrak{g}^{e}$ and $\mathfrak{z}(\mathfrak{g}^{e})$
which help calculations later. In Sections \ref{sec:D(2,1;)}\textendash \ref{sec:F(4)},
we recall a construction of exceptional Lie superalgebras $\mathfrak{g}=D(2,1;\alpha)$,
$G(3)$, $F(4)$ and use this to explicitly determine the centralizers
$\mathfrak{g}^{e}$ and centres $\mathfrak{z}(\mathfrak{g}^{e})$
of centralizers of nilpotent even elements $e$ in $\mathfrak{g}$.
We also give all possible Dynkin diagrams for $\mathfrak{g}$ and
further determine the labelled Dynkin diagrams with respect to $e$.
\begin{singlespace}

\section{Preliminaries\label{sec:preliminaries}}
\end{singlespace}

\subsection{Basic classical Lie superalgebras}

\begin{singlespace}
\noindent Recall that a finite-dimensional simple Lie superalgebra
$\mathfrak{g}=\mathfrak{g}_{\bar{0}}\oplus\mathfrak{g}_{\bar{1}}$
over $\mathbb{C}$ is called a basic classical Lie superalgebra if
$\mathfrak{g}_{\bar{0}}$ is a reductive Lie algebra and there exists
a non-degenerate even invariant supersymmetric bilinear form $(\cdotp,\cdotp)$
on $\mathfrak{g}$. Kac gives the classification of finite-dimensional
complex simple Lie superalgebras in \cite[Theorem 5]{Kac1977}. He
argues that simple basic classical Lie superalgebras that are not
Lie algebras consist of classical types which are infinite families
and three exceptional types. In this paper, we consider the Lie superalgebras
$D(2,1;\alpha)$, $G(3)$ and $F(4)$, which are called the exceptional
Lie superalgebras. We will describe a construction of $D(2,1;\alpha)$,
$G(3)$ and $F(4)$ in \S4.1, \S5.1 and \S6.3 respectively. 
\end{singlespace}

\begin{singlespace}
Let $\mathfrak{g}=\mathfrak{g}_{\bar{0}}\oplus\mathfrak{g}_{\bar{1}}$
be a finite-dimensional basic classical Lie superalgebra and $\text{\ensuremath{\mathfrak{h}}}$
be a Cartan subalgebra of $\mathfrak{g}$. The following concepts
can be found in \cite[Chapter 1]{Cheng2012}. There is a root space
decomposition $\mathfrak{g}=\mathfrak{h}\oplus\bigoplus_{\alpha\in\mathfrak{h}^{*}}\mathfrak{g_{\alpha}}$
where $\mathfrak{h}=\mathfrak{g}_{0}=\{x\in\mathfrak{g}:[h,x]=0\ \text{for\ all}\ h\in\mathfrak{h}\}=\mathfrak{g}^{\mathfrak{h}}$
and $\mathfrak{g_{\alpha}}:=\{x\in\mathfrak{g}:[h,x]=\alpha(h)x\ \text{for\ all}\ h\in\mathfrak{\mathfrak{h}}\}$
are the corresponding weight spaces. Note that $\mathfrak{g_{\alpha}}$
is 1-dimensional. The set $\Phi=\{\alpha\in\mathfrak{h^{*}:\alpha\neq}0,\mathfrak{g}_{\alpha}\neq0\}$
is called the root system of $\mathfrak{g}$ and we say that $\mathfrak{g}_{\alpha}$
is the root space corresponding to the root $\alpha\in\Phi$. The
even and odd roots are defined to be $\Phi_{\bar{0}}=:\{\alpha\text{\ensuremath{\in}}\Phi:\mathfrak{g}_{\alpha}\subseteq\mathfrak{g}_{\bar{0}}\}$
and $\Phi_{\bar{1}}=:\{\alpha\in\Phi:\mathfrak{g}_{\alpha}\subseteq\mathfrak{g}_{\bar{1}}\}$.
A subset $\Phi^{+}\subseteq\Phi$ is a system of positive roots if
for each root $\alpha\in\Phi$ there is exactly one of $\alpha$,$-\alpha$
contained in $\Phi^{+}$; and for any two distinct roots $\alpha,\beta\in\Phi^{+}$,
$\alpha+\beta\in\Phi$ implies that $\alpha+\beta\in\Phi^{+}$. Given
a system of positive roots $\Phi^{+}$, elements of $-\Phi^{+}$ form
a system of negative roots. Note that $\Phi^{+}=\Phi_{\bar{0}}^{+}\cup\Phi_{\bar{1}}^{+}$.
A system of simple roots $\varPi=\{\alpha_{1},...,\alpha_{l}\}\subseteq\Phi^{+}$
consists of all elements that cannot be written as the sum of two
elements of $\Phi^{+}$. Note that $l$ does not depend on choice
of $\varPi$ and we call it the rank of $\mathfrak{g}$.
\end{singlespace}

\subsection{Dynkin diagrams\label{subsec:Dynkin-diagrams}}

\begin{singlespace}
\noindent Let $\mathfrak{g}=\mathfrak{g}_{\bar{0}}\oplus\mathfrak{g}_{\bar{1}}$
be a basic classical Lie superalgebra with a Cartan subalgebra $\mathfrak{h}\subseteq\mathfrak{g}_{\bar{0}}$.
There exists a\textit{ }triangular decomposition\textit{ }$\mathfrak{g}=\mathfrak{n}^{-}\oplus\mathfrak{h}\oplus\mathfrak{n}^{+}$
where $\mathfrak{n}^{+}$ (resp. $\mathfrak{n}^{-}$) is a subalgebra
such that $[\mathfrak{h},\mathfrak{n}^{+}]\subseteq\mathfrak{n}^{+}$
(resp. $[\mathfrak{h},\mathfrak{n}^{-}]\subseteq\mathfrak{n}^{-}$)
and $\dim\mathfrak{n}^{+}=\dim\mathfrak{n}^{-}$, see \cite[Section 2]{Penkov}.
The solvable subalgebra $\mathfrak{b}=\mathfrak{h}\oplus\mathfrak{n}^{+}$
is called a Borel subalgebra of $\mathfrak{g}$. We work with Borel
subalgebras up to conjugacy by $G$. Note that there are in general
many inequivalent conjugacy classes of Borel subalgebras and every
Borel subalgebra containing $\mathfrak{h}$ determines a corresponding
system of positive roots $\Phi^{+}$. Consequently $\mathfrak{b}$
determines a system of simple roots $\varPi$. Then for each conjugacy
class of Borel subalgebras of $\mathfrak{g}$, a simple root system
can be transformed into an equivalent one with the same Dynkin diagram
under the transformation of the Weyl group $W$ of $\mathfrak{g}$,
see \cite[Subsection 2.3]{Frappat1989}. 
\end{singlespace}

\begin{singlespace}
We next recall the concept of the Dynkin diagram as defined for example
in \cite[Section 2.2]{Frappat1989}. We know that there exists a non-degenerate
even invariant supersymmetric bilinear form $(\cdotp,\cdotp)$ on
$\mathfrak{g}$. One can check that $(\cdotp,\cdotp)$ restricts to
a non-degenerate symmetric bilinear form on $\mathfrak{h}$. Therefore,
there exists an isomorphism from $\mathfrak{h}$ to $\mathfrak{h}^{*}$
which provides a symmetric bilinear form on $\mathfrak{h}^{*}$. Then
the \textit{Dynkin diagram} of a Lie superalgebra $\mathfrak{g}$
with a simple root system $\varPi$ is a graph where the vertices
are labelled by $\varPi$ and there are $\mu_{\alpha\beta}$ lines
between the vertices labelled by simple roots $\alpha_{i}$ and $\alpha_{j}$
such that:
\begin{equation}
\mu_{\alpha\beta}=\begin{cases}
\ensuremath{\vert(\alpha_{i},\alpha_{j})\ensuremath{\vert}} & \text{if }(\alpha_{i},\alpha_{i})=(\alpha_{j},\alpha_{j})=0,\\
\frac{2\ensuremath{\vert}(\alpha_{i},\alpha_{j})\ensuremath{\vert}}{min\{\vert(\alpha_{i},\alpha_{i})\vert,\ensuremath{\vert}(\alpha_{j},\alpha_{j})\ensuremath{\vert}\}} & \text{if }(\alpha_{i},\alpha_{i})(\alpha_{j},\alpha_{j})\neq0,\\
\frac{2\ensuremath{\vert}(\alpha_{i},\alpha_{j})\ensuremath{\vert}}{min_{(\alpha_{k},\alpha_{k})\neq0}\ensuremath{\vert}(\alpha_{k},\alpha_{k})\ensuremath{\vert}} & \text{if }(\alpha_{i},\alpha_{i})\neq0,\ (\alpha_{j},\alpha_{j})=0\ \text{\text{and}\ }\alpha_{k}\in\Phi.
\end{cases}\label{eq:lines-=0003BC}
\end{equation}
 We say a root $\alpha\in\Phi$ is \textit{isotropic} if $(\alpha,\alpha)=0$
and is \textit{non-isotropic} if $(\alpha,\alpha)\neq0$. We associate
a white node $\ocircle$ to each even root, a grey node $\varotimes$
to each odd isotropic root and a black node $\newmoon$ to each odd
non-isotropic root. Moreover, when $\mu_{\alpha\beta}>1$, we put
an arrow pointing from the vertex labelled by $\alpha_{i}$ to the
vertex labelled by $\alpha_{j}$ if $(\alpha_{i},\alpha_{i})(\alpha_{j},\alpha_{j})\neq0$
and $(\alpha_{i},\alpha_{i})>(\alpha_{j},\alpha_{j})$ or if $(\alpha_{i},\alpha_{i})=0,(\alpha_{j},\alpha_{j})\neq0$
and $\vert(\alpha_{j},\alpha_{j})\vert<2$, or pointing from the vertex
labelled by $\alpha_{j}$ to the vertex labelled by $\alpha_{i}$
if $(\alpha_{i},\alpha_{i})=0,(\alpha_{j},\alpha_{j})\neq0$ and $\vert(\alpha_{j},\alpha_{j})\vert>2$.
If the value of $\mu_{\alpha\beta}$ is not a natural number, then
we label the edge between vertices corresponding to roots $\alpha$
and $\beta$ with $\mu_{\alpha\beta}$ instead of drawing multiple
lines between them. Note that the Dynkin diagram depends on $\varPi$
up to conjugacy, thus Dynkin diagrams of $\mathfrak{g}$ for different
choices of simple roots can be different.
\end{singlespace}

\subsection{Labelled Dynkin diagrams\label{subsec:Lablled-Dynkin-diagrams}}

\noindent Let $e\in\mathfrak{g}_{\bar{0}}$ be nilpotent. There exists
an $\mathfrak{sl}(2)$-triple $\{e,h,f\}\subseteq\mathfrak{g}_{\bar{0}}$
by the Jacobson\textendash Morozov Theorem, see for example \cite[Theorem 3.3.1]{Collingwood1993}.
An $\mathfrak{sl}(2)$-triple determines a grading on $\mathfrak{g}$
according to the eigenvalues of ad$h,$ thus we can decompose $\mathfrak{g}$
into its ad$h$-eigenspaces $\mathfrak{g}=\bigoplus_{j\in\mathbb{Z}}\mathfrak{g}(j)$
where $\mathfrak{g}(j)=\{x\in\mathfrak{g}:[h,x]=jx\}$. We can choose
a Borel subalgebra $\mathfrak{b}\subseteq\bigoplus_{j\geq0}\mathfrak{g}(j)$
and a Cartan subalgebra $\mathfrak{h}\subseteq\mathfrak{b}$. Then
we obtain the corresponding system of positive roots $\Phi^{+}$ and
a system of simple roots $\varPi=\{\alpha_{1},...,\alpha_{l}\}$ which
will give a Dynkin diagram of $\mathfrak{g}$. Furthermore, for each
$i=1,...,l$, note that $\mathfrak{g}_{\alpha_{i}}$ is the root space
corresponding to $\alpha_{i}$ and $\mathfrak{g}_{\alpha_{i}}\subseteq\mathfrak{g}(j_{i})$
for some $j_{i}\geq0$. Hence, we have $\alpha_{i}(h)\geq0$ for $i=1,...,l$. 
\begin{defn}
\begin{singlespace}
\noindent The \textit{labelled Dynkin diagram} $\varDelta$ of $e$
determined by $\varPi$ is given by taking the Dynkin diagram of $\mathfrak{g}$
and labelling each node $\alpha$ with $\alpha(h)$.
\end{singlespace}
\end{defn}

\begin{singlespace}

\section{Generalities on $\mathfrak{g}^{e}$ and $\mathfrak{z}(\mathfrak{g}^{e})$\label{sec:general-method}}
\end{singlespace}

\begin{singlespace}
\noindent Let $\mathfrak{g}=\mathfrak{g}_{\bar{0}}\oplus\mathfrak{g}_{\bar{1}}$
be a basic classical Lie superalgebra and $e\in\mathfrak{g}_{\bar{0}}$
be nilpotent. In this section, we give an overview of some general
methods for calculating $\mathfrak{g}^{e}$ and $\mathfrak{z}(\mathfrak{g}^{e})$.
\end{singlespace}

\begin{singlespace}
Note that any element $x\in\mathfrak{g}$ can be written as $x=x_{\bar{0}}+x_{\bar{1}}$
such that $x_{\bar{i}}\in\mathfrak{g}_{\bar{i}}$. For a nilpotent
element $e\in\mathfrak{g}_{\bar{0}}$, if $[x,e]=0$ then $[x,e]=[x_{\bar{0}},e]+[x_{\bar{1}},e]=0$.
This implies that $[x_{\bar{0}},e]=[x_{\bar{1}},e]=0$ since $[x_{\bar{0}},e]\in\mathfrak{g}_{\bar{0}}$
and $[x_{\bar{1}},e]\in\mathfrak{g}_{\bar{1}}$. Hence $\mathfrak{g}^{e}=\mathfrak{g}_{\bar{0}}^{e}\oplus\mathfrak{g}_{\bar{1}}^{e}$.

For a nilpotent element $e\in\mathfrak{g}_{\bar{0}}$, recall that
there exists an $\mathfrak{sl}(2)$-triple $\{e,h,f\}\subseteq\mathfrak{g}_{\bar{0}}$
as noted in \S2.3 and any two $\mathfrak{sl}(2)$-triples containing
$e$ are conjugate under the action of the group $G^{e}$. We have
that $\mathfrak{g}$ is a module for $\mathfrak{s}=\left\langle e,h,f\right\rangle $
via the adjoint action. Define $V^{\mathfrak{sl}}(d)$ to be the $(d+1)$-dimensional
simple $\mathfrak{sl}(2)$-module with highest weight $d$. By the
representation theory of $\mathfrak{sl}(2)$, we can decompose $\mathfrak{g}$
into a direct sum of finite-dimensional $\mathfrak{s}$-submodules
$\mathfrak{g}^{i}$ and each of them is isomorphic to $V^{\mathfrak{sl}}(d_{i})$
for some $d_{i}\in\mathbb{Z}$ and $d_{i}\geq0$. The element $h$
of the $\mathfrak{sl}(2)$-triple is semisimple and the eigenvalues
of $h$ on $\mathfrak{g}^{i}$ are $d_{i},d_{i}-2,...,-(d_{i}-2),-d_{i}$.
The only vectors in $\mathfrak{g}^{i}$ annihilated by $e$ are the
multiples of the highest weight vector, i.e. if $\mathfrak{g}^{i}$
has basis $\{x_{d_{i}}^{i},x_{d_{i}-2}^{i},...,x_{-d_{i}+2}^{i},x_{-d_{i}}^{i}\}$
for $i=1,2,...,r$, then the vectors annihilated by $e$ are $\left\langle x_{d_{i}}^{i}\right\rangle $.
Thus the vectors in $\mathfrak{g}$ centralized by $e$ are $\left\langle x_{d_{1}}^{1},x_{d_{2}}^{2},...,x_{d_{r}}^{r}\right\rangle $
and they have ad$h$-eigenvalues $d_{1},d_{2},...,d_{r}$. Hence,
from the $\mathrm{ad}h$-eigenspace decomposition of $\mathfrak{g}$
we determine the $\mathrm{ad}h$ eigenvalues of elements of $\mathfrak{g}^{e}$.

We consider $\mathfrak{sl}(2)$ frequently in the following sections,
so we fix the notation $\mathfrak{sl}(2)=\left\langle E,H,F\right\rangle $
where
\[
E=\begin{pmatrix}0 & 1\\
0 & 0
\end{pmatrix},H=\begin{pmatrix}1 & 0\\
0 & -1
\end{pmatrix},F=\begin{pmatrix}0 & 0\\
1 & 0
\end{pmatrix}.
\]
The commutator relations between basis elements for $\mathfrak{sl}(2)$
are $[H,E]=2E$, $[H,F]=-2F$ and $[E,F]=H$. Let $V$ be a two-dimensional
vector space with basis $v_{1}=(1,0)^{t}$ and $v_{-1}=(0,1)^{t}$.

When describing the $\mathfrak{g}^{e}(0)$-module structure on each
$\mathfrak{g}^{e}(j)$ for $j>0$, we need the following lemma:
\end{singlespace}
\begin{lem}
\begin{singlespace}
\noindent \textup{\label{lem:osp(1,2)}Let $A=\mathfrak{g}=\mathfrak{g}_{\bar{0}}\oplus\mathfrak{g}_{\bar{1}}$
be a Lie superalgebra where $\{u_{-2},u_{0},u_{2}\}$ is a basis of
$\mathfrak{g}_{\bar{0}}$ and $\{u_{-1},u_{1}\}$ is a basis of $\mathfrak{g}_{\bar{1}}$
such that: (1) $[u_{0},u_{i}]=a_{i}u_{i}$; (2) $[u_{1},u_{1}]=au_{2}$
and $[u_{-1},u_{-1}]=bu_{-2}$; (3) $[u_{2},u_{-2}]=cu_{0}$ for $a_{i},a,b,c\neq0$.
Then $A$ is simple and $A\cong\mathfrak{osp}(1|2)$ .}
\end{singlespace}
\end{lem}

\begin{singlespace}
\noindent \begin{proof}Let $I$ be an non-zero ideal of $A$. Then
$I$ is a direct sum of $\mathrm{ad}u_{0}$ eigenspaces, thus $u_{i}\in I$
for some $i$. If $i=0$, then condition (1) implies that $I=A$.
If $i=\pm2$, then condition $(3)$ implies that $u_{0}\in I$ and
thus $I=A$. If $i=\pm1$, then condition (2) implies that $u_{-2}$
or $u_{2}$ lies in $I$. Thus $u_{0}\in I$ and $I=A$. Therefore,
we have that $A$ is simple. According to the classification Theorem
of simple Lie superalgebras in \cite[Theorem 5]{Kac1977}, we deduce
that $A\cong\mathfrak{osp}(1|2)$.\end{proof}
\end{singlespace}

\begin{singlespace}
We consider the representations of $\mathfrak{osp}(1|2)$ frequently
in the following sections. As shown in \cite[Section 2]{Pais1975},
all finite-dimensional representations of $\mathfrak{osp}(1|2)$ are
completely reducible. Also in \cite[Section 2]{Pais1975} the irreducible
representations of $\mathfrak{osp}(1|2)$ are constructed. We recall
that the irreducible representations of $\mathfrak{osp}(1|2)$ are
parametrized by $l\in\{\frac{a}{2}:a\in\mathbb{Z}_{\geq0}\}$ and
we write $V^{\mathfrak{osp}}(l)$ for the representation corresponding
to $l$. Then $\dim V^{\mathfrak{osp}}(l)=4l+1$. We know that $\mathfrak{osp}(1|2)$
is $5$-dimensional with basis $\{u_{-2},u_{-1},u_{0},u_{1},u_{2}\}$.
The eigenvalues of $u_{0}$ on $V^{\mathfrak{osp}}(l)$ are $l,l-\frac{1}{2},\dots,-l$.

From now on let us denote $\mathfrak{z}=\mathfrak{z}(\mathfrak{g}^{e})$.
Given $x=x_{\bar{0}}+x_{\bar{1}}\in\mathfrak{z}$, for any $y=y_{\bar{0}}+y_{\bar{1}}\in\mathfrak{g}^{e}$,
we have $[x,y_{\bar{0}}]=[x_{\bar{0}},y_{\bar{0}}]+[x_{\bar{1}},y_{\bar{0}}]=0$.
Since $[x_{\bar{0}},y_{\bar{0}}]\in\mathfrak{g}_{\bar{0}}^{e}$ and
$[x_{\bar{1}},y_{\bar{0}}]\in\mathfrak{g}_{\bar{1}}^{e}$, we have
$[x_{\bar{0}},y_{\bar{0}}]=[x_{\bar{1}},y_{\bar{0}}]=0$. Similarly
we have $[x_{\bar{0}},y_{\bar{1}}]=[x_{\bar{1}},y_{\bar{1}}]=0$.
Therefore, we know that $x_{\bar{0}},x_{\bar{1}}\in\mathfrak{z}$
and thus $\mathfrak{z}=\mathfrak{z}_{\bar{0}}\oplus\mathfrak{z}_{\bar{1}}$.
Moreover, we can decompose $\mathfrak{z}$ into the direct sum of
ad$h$-eigenspaces in each case, i.e. $\mathfrak{z}=\text{\ensuremath{\bigoplus}}_{j\geq0}\mathfrak{z}(j)$
for all ad$h$-eigenvalue $j\geq0$. 
\end{singlespace}

In \S5.5 and \S6.7, we consider the adjoint action of group $G^{e}$
on $\mathfrak{z}$ for $\mathfrak{g}=G(3)$ and $F(4)$. According
to \cite[Section 5.10]{Jantzen2004a}, we have that $G^{e}$ is the
semidirect product as an algebraic group of the reductive group $C^{e}=G^{e}\cap G^{h}$
and $R^{e}$, the unipotent radical of $G^{e}$. Denote the connected
component of $G^{e}$ (resp. $C^{e}$) containing the identity by
$(G^{e})^{\circ}$ (resp. $(C^{e})^{\circ}$). The group $R^{e}$
is connected by \cite[Section 5.10]{Jantzen2004a}, thus we get $G^{e}/(G^{e})^{\circ}\cong C^{e}/(C^{e})^{\circ}$.
Since the adjoint action of $\mathfrak{g}^{e}$ on itself is the differential
of the adjoint of $G^{e}$ on $\mathfrak{g}^{e}$, we have $\mathfrak{z}=\{x\in\mathfrak{g}^{e}:g\cdot x=x\text{ for all }g\in(G^{e})^{\circ}\}$.
Therefore, there is an action of $G^{e}/(G^{e})^{\circ}$ on $\mathfrak{z}$
and $\mathfrak{z}^{G^{e}}\subseteq\mathfrak{z}^{G^{e}/(G^{e})^{\circ}}\subseteq\mathfrak{z}$.
In order to determine $\mathfrak{z}^{G^{e}}$, it suffices to consider
the action of elements of $G^{e}/(G^{e})^{\circ}$ on $\mathfrak{z}$.
\begin{singlespace}

\section{Exceptional Lie superalgebras $D(2,1;\alpha)$\label{sec:D(2,1;)}}
\end{singlespace}

\begin{singlespace}
\noindent In this section, we describe an explicit construction for
the Lie Superalgebras $\mathfrak{g}=D(2,1;\alpha)$ following \cite{M.Scheunert1976}
and \cite[Section 4.2]{Musson2012}. We also give representatives
of nilpotent orbits $e\in D(2,1;\alpha)_{\bar{0}}$. We use the explicit
construction in \S4.1 to determine $\mathfrak{g}^{e}$ and $\mathfrak{z}(\mathfrak{g}^{e})$.
The labelled Dynkin diagram $\varDelta$ with respect to each nilpotent
elements are drawn afterwards.
\end{singlespace}

\subsection{Structure of Lie superalgebras $D(2,1;\alpha)$\label{subsec:Structure-of-D(2,1;)}}

\begin{singlespace}
\noindent The Lie superalgebras $D(2,1;\alpha)$ with $\alpha\in\mathbb{C}$\textbackslash$\left\{ 0,1\right\} $
form a one-parameter family of superalgebras of dimension $17$. In
\cite{M.Scheunert1976}, Scheunert denotes these algebra by $\Gamma(\sigma_{1},\sigma_{2},\sigma_{3})$
where $\sigma_{1},\sigma_{2},\sigma_{3}$ are complex numbers such
that $\sigma_{1}+\sigma_{2}+\sigma_{3}=0$. According to \cite[Section 4.2]{Musson2012},
the Lie superalgebra $\Gamma(\sigma_{1},\sigma_{2},\sigma_{3})$ is
simple if and only if $\sigma_{i}\neq0$ for $i=1,2,3$ . If there
exists another triple $(\sigma_{1}^{'},\sigma_{2}^{'},\sigma_{3}^{'})$
such that $\Gamma(\sigma_{1},\sigma_{2},\sigma_{3})\cong\Gamma(\sigma_{1}^{'},\sigma_{2}^{'},\sigma_{3}^{'})$,
then there must exist a permutation $\rho$ of $\{1,2,3\}$ and a
nonzero complex number $c$ such that $\sigma_{i}^{'}=c\sigma_{\rho(i)}$
for $i=1,2,3$. This implies that the $\Gamma(\sigma_{1},\sigma_{2},\sigma_{3})$
form a one-parameter family and we have that $\Gamma(\sigma_{1},\sigma_{2},\sigma_{3})=D(2,1;\alpha)$
for a specific choice of $\sigma_{1},\sigma_{2},\sigma_{3}$. For
any $\alpha\in\mathbb{C}\setminus\{0,-1\}$, we have $D(2,1;\alpha)=\Gamma(1+\alpha,-1,-\alpha)\cong\Gamma(\frac{1+\alpha}{\alpha},-1,-\frac{1}{\alpha})\cong\Gamma(-\alpha,-1,1+\alpha)$. 
\end{singlespace}

\begin{singlespace}
For $i=1,2,3$, take $V_{i}$ to be a copy of $V$ where $V$ is defined
in \S3. Let $\psi_{i}$ be the non-degenerate skew-symmetric bilinear
form on $V_{i}$ defined by $\psi_{i}(v_{1},v_{-1})=1$. We also define
a bilinear map $p_{i}:V_{i}\times V_{i}\rightarrow\mathfrak{sl}(2)$
by $p_{i}(x,y)(z)=\psi_{i}(y,z)x-\psi_{i}(z,x)y$ for $x,y,z\in V_{i}$.
We can calculate that $p_{i}(v_{1},v_{1})=2E$, $p_{i}(v_{1},v_{-1})=-H$
and $p_{i}(v_{-1},v_{-1})=-2F$.

By definition, $\mathfrak{g}=D(2,1;\alpha)=\mathfrak{g}_{\bar{0}}\oplus\mathfrak{g}_{\bar{1}}$,
where
\[
\mathfrak{g}_{\bar{0}}=\mathfrak{sl}(2)\oplus\mathfrak{sl}(2)\oplus\mathfrak{sl}(2)\text{ and }\mathfrak{g}_{\bar{1}}=V_{1}\otimes V_{2}\otimes V_{3}.
\]

\end{singlespace}

\begin{singlespace}
\noindent Note that $\mathfrak{g}_{\bar{0}}$ is a Lie algebra thus
has Lie bracket $\left[\cdot,\cdot\right]:\mathfrak{g}_{\bar{0}}\times\mathfrak{g}_{\bar{0}}\rightarrow\mathfrak{g}_{\bar{0}}$.
Let $x=(x_{1},x_{2},x_{3})\in\mathfrak{g}_{\bar{0}}$ and $v=v_{1}\otimes v_{2}\otimes v_{3}\in\mathfrak{g}_{\bar{1}}$,
then the bracket $\left[\cdot,\cdot\right]:\mathfrak{g}_{\bar{0}}\times\mathfrak{g}_{\bar{1}}\rightarrow\mathfrak{g}_{\bar{1}}$
is defined by $[x,v]:=x\cdotp v=x_{1}v_{1}\otimes v_{2}\otimes v_{3}+v_{1}\otimes x_{2}v_{2}\otimes v_{3}+v_{1}\otimes v_{2}\otimes x_{3}v_{3}.$
According to equation (4.2.1) in \cite{Musson2012}, the bracket $\left[\cdot,\cdot\right]:\mathfrak{g}_{\bar{1}}\times\mathfrak{g}_{\bar{1}}\rightarrow\mathfrak{g}_{\bar{0}}$
is given by 
\begin{align*}
[v_{1}\otimes v_{2}\otimes v_{3},u_{1}\otimes u_{2}\otimes u_{3}] & =\sigma_{1}\psi_{2}(v_{2},u_{2})\psi_{3}(v_{3},u_{3})p_{1}(v_{1},u_{1})\\
+\sigma_{2}\psi_{1}(v_{1},u_{1})\psi_{3}(v_{3},u_{3})p_{2}(v_{2},u_{2}) & +\sigma_{3}\psi_{1}(v_{1},u_{1})\psi_{2}(v_{2},u_{2})p_{3}(v_{3},u_{3}),
\end{align*}
where $v_{i},u_{i}\in V_{i}$.
\end{singlespace}

\begin{singlespace}
Next we give a basis for $\mathfrak{g}$. We first fix the following
notation: let $E_{1}=(E,0,0)$, $E_{2}=(0,E,0)$ and $E_{3}=(0,0,E)$.
Similarly, we denote $F_{1}=(F,0,0)$, $F_{2}=(0,F,0)$, $F_{3}=(0,0,F)$,
$H_{1}=(H,0,0)$, $H_{2}=(0,H,0)$ and $H_{3}=(0,0,H)$. Clearly,
$\mathfrak{g}_{\bar{0}}$ has a basis $\{E_{i},H_{i},F_{i}:i=1,2,3\}$
and $\mathfrak{g}_{\bar{1}}$ has a basis $\{v_{i}\otimes v_{j}\otimes v_{k}:i,j,k=\pm1\}$.
In the remaining subsections, we write $v_{i,j,k}$ in place of $v_{i}\otimes v_{j}\otimes v_{k}$
for $i,j,k\in\{\pm1\}$.
\end{singlespace}

\subsection{Root system and Dynkin diagrams for $D(2,1;\alpha)$\label{subsec:Root-system-D(2,1)}}

\begin{singlespace}
\noindent We follow the construction of the root system of $D(2,1;\alpha)$
given in \cite[Appendix A]{Iohara2001}. A Lie superalgebra of type
$D(2,1;\alpha)$ has root system 
\[
\Phi_{\bar{0}}=\{\pm2\beta_{1},\pm2\beta_{2},\pm2\beta_{3}\}\text{ and\ }\Phi_{\bar{1}}=\{i\beta_{1}+j\beta_{2}+k\beta_{3}:i,j,k=\pm1\},
\]

\noindent where $\{\beta_{1},\beta_{2},\beta_{3}\}$ is an orthogonal
basis such that $(\beta_{1},\beta_{1})=\frac{1}{2}$, $(\beta_{2},\beta_{2})=-\frac{1}{2}\alpha-\frac{1}{2}$
and $(\beta_{3},\beta_{3})=\frac{1}{2}\alpha$. The corresponding
root vectors are listed below:
\end{singlespace}
\begin{singlespace}
\noindent \begin{center}
\begin{longtable}[c]{|c||c||c||c|}
\hline 
\multicolumn{1}{|c||}{Roots} & $2\beta_{i},i=1,2,3$ & $-2\beta_{i},i=1,2,3$ & $i\beta_{1}+j\beta_{2}+k\beta_{3}:i,j,k\in\{\pm1\}$\tabularnewline
\hline 
\endfirsthead
\hline 
Root vectors & $E_{i}$ & $F_{i}$ & $v_{i,j,k}$ \tabularnewline
\hline 
\end{longtable}
\par\end{center}
\end{singlespace}

\noindent We can check that all the odd roots in $D(2,1;\alpha)$
are isotropic. In the following table, we give all possible Dynkin
diagrams with respect to different systems of simple roots based on
\cite[Section 2.20]{Frappat1996}.

\begin{longtable}[c]{|>{\centering}m{7cm}||>{\centering}m{5cm}|}
\caption{Dynkin diagrams for $D(2,1;\alpha)$}
\tabularnewline
\endfirsthead
\hline 
Simple systems $\varPi=\{\alpha_{1},\alpha_{2},\alpha_{3}\}$ & Dynkin diagrams\tabularnewline
\hline 
\hline 
$\{2\beta_{1},-\beta_{1}+\beta_{2}-\beta_{3},2\beta_{3}\}$ & Figure 4.1

\includegraphics[scale=0.8]{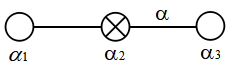}\tabularnewline
\hline 
\hline 
 $\{2\beta_{1},-\beta_{1}-\beta_{2}+\beta_{3},2\beta_{2}\}$ & Figure 4.2

\includegraphics[scale=0.8]{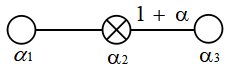}\tabularnewline
\hline 
\hline 
$\{2\beta_{3},\beta_{1}-\beta_{2}-\beta_{3},2\beta_{2}\}$ & Figure 4.3

\includegraphics[scale=0.8]{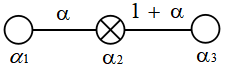}\tabularnewline
\hline 
\hline 
 $\{-\beta_{1}+\beta_{2}+\beta_{3},\beta_{1}-\beta_{2}+\beta_{3},\beta_{1}+\beta_{2}-\beta_{3}\}$ & Figure 4.4

\includegraphics{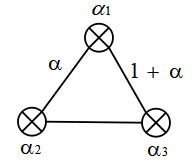}\tabularnewline
\hline 
\end{longtable}

\subsection{Centres of centralizers of nilpotent elements $e\in D(2,1;\alpha)$
and labelled Dynkin diagrams with respect to $e$\label{subsec:Centres-of-centralizers-D(2,1)}}

\begin{singlespace}
\noindent Let $\mathfrak{g}=D(2,1;\alpha)=\mathfrak{g}_{\bar{0}}\oplus\mathfrak{g}_{\bar{1}}$.
A nilpotent element $e\in\mathfrak{g}_{\bar{0}}$ is of the form $(e_{1},e_{2},e_{3})$
where $e_{i}\in\mathfrak{sl}(2)$ for $i\in\{1,2,3\}$. We know that
representatives of nilpotent elements in $\mathfrak{sl}(2)$ up to
conjugation by $\mathrm{SL}(2)$ are $0\ \text{and}\ E$. We give
basis elements for $\mathfrak{g}^{e}$ and $\mathfrak{z}(\mathfrak{g}^{e})$
and labelled Dynkin diagrams $\varDelta$ with respect to $e$ when
$e=0,E_{1},E_{1}+E_{2},E_{1}+E_{2}+E_{3}$ in Table \ref{tab:results in D(2,1)}.
Note that the cases $e=E_{2}$, $e=E_{3}$ are similar to $e=E_{1}$
and cases $e=E_{2}+E_{3}$, $e=E_{1}+E_{3}$ are similar to $e=E_{1}+E_{2}$.
Hence, any other case is similar to one of above. The numbers in the
column labelled by ``$\varDelta$'' represent labels $a_{i}$ corresponding
to $\alpha_{i}$ for $i=1,2,3$ in labelled Dynkin diagram with respect
to $e$.

\noindent 
\begin{table}[H]
\begin{singlespace}
\noindent \begin{centering}
\begin{tabular}{|>{\centering}m{1.5cm}||>{\centering}m{8cm}||>{\centering}m{1cm}||>{\centering}m{3.3cm}|}
\hline 
$e$ & $\mathfrak{g}^{e}$  & $\mathfrak{z}(\mathfrak{g}^{e})$ &  $\varDelta$ \tabularnewline
\hline 
\hline 
$0$ & $\mathfrak{g}$ & $\{0\}$ & Figures 4.1, 4.2, 4.3, 4.4: All labels are zeros.\tabularnewline
\hline 
\hline 
\textbf{$E_{1}$} & $\langle E_{1},E_{2},H_{2},F_{2},E_{3},H_{3},F_{3},v_{i,j,k}:j,k=\pm1\rangle$ & $\langle e\rangle$ & Figure 4.3: $0,1,0$\tabularnewline
\hline 
\hline 
$E_{1}+E_{2}$ & \begin{singlespace}
\noindent $\langle E_{1},E_{2},E_{3},H_{3},F_{3},v_{1,1,1},v_{1,1,-1},v_{1,-1,1}-v_{-1,1,1},v_{1,-1,-1}-v_{-1,1,-1}\rangle$
\end{singlespace}
 & $\langle e\rangle$ & Figure 4.1: $2,0,0$ 

Figure 4.3: $0,0,2$ 

Figure 4.4: $0,0,2$\tabularnewline
\hline 
\hline 
$E_{1}+E_{2}+E_{3}$ & $\langle E_{1},E_{2},E_{3},v_{1,1,1},v_{1,1,-1}-v_{-1,1,1},v_{1,-1,1}-v_{-1,1,1}\rangle$ & $\langle e\rangle$ & Figure 4.4: $1,1,1$\tabularnewline
\hline 
\end{tabular}
\par\end{centering}
\end{singlespace}
\caption{\label{tab:results in D(2,1)}$\mathfrak{g}^{e}$, $\mathfrak{z}(\mathfrak{g}^{e})$
and $\varDelta$ for $\mathfrak{g}=D(2,1;\alpha)$}

\end{table}
 
\end{singlespace}

\begin{singlespace}
Let $V^{\mathfrak{sl}}(j)$ be an $\mathfrak{sl}(2)$-module with
highest weight $j$ and $V^{\mathfrak{osp}}(j)$ be an $\mathfrak{osp}(1|2)$-module
with highest weight $j$. We also describe the $\mathfrak{g}^{e}(0)$-module
structure on each $\mathfrak{g}^{e}(j)$ for $j>0$ in Table \ref{tab:g^e(0)-D(2,1;)}.
\end{singlespace}
\begin{singlespace}
\noindent \begin{center}
\begin{longtable}[c]{|>{\centering}m{2.5cm}|>{\centering}m{2cm}|>{\centering}m{9cm}|}
\caption{\label{tab:g^e(0)-D(2,1;)}The $\mathfrak{g}^{e}(0)$-module structure
on $\mathfrak{g}^{e}(j)$ for $j>0$ }
\tabularnewline
\endfirsthead
\hline 
$e$ & $\mathfrak{g}^{e}(0)$  & $\mathfrak{g}^{e}(j)$ for $j>0$\tabularnewline
\hline 
\hline 
$0$ & $\mathfrak{g}^{e}$ & $0$\tabularnewline
\hline 
\hline 
\textbf{$E_{1}$} & $\mathfrak{sl}(2)\oplus\mathfrak{sl}(2)$ & $\mathfrak{g}^{e}(1)=V^{\mathfrak{sl}}(1)\otimes V^{\mathfrak{sl}}(1)$,
$\mathfrak{g}^{e}(2)=V^{\mathfrak{sl}}(0)\otimes V^{\mathfrak{sl}}(0)$\tabularnewline
\hline 
\hline 
$E_{1}+E_{2}$ & $\mathfrak{osp}(1|2)$ & $\mathfrak{g}^{e}(2)=V^{\mathfrak{osp}}(0)\oplus V^{\mathfrak{osp}}(1)$\tabularnewline
\hline 
\hline 
$E_{1}+E_{2}+E_{3}$ & $\{0\}$ & $\dim\mathfrak{g}^{e}(1)=2$, $\dim\mathfrak{g}^{e}(2)=3$, $\dim\mathfrak{g}^{e}(3)=1$\tabularnewline
\hline 
\end{longtable}
\par\end{center}
\end{singlespace}

\begin{singlespace}
Let $e=E_{1}+E_{2}$. In the remaining part of this subsection, we
calculate $\mathfrak{g}^{e}$, $\mathfrak{z}(\mathfrak{g}^{e})$ and
draw the labelled Dynkin diagrams with respect to $e$.

We easily calculate that $\mathfrak{g}_{\bar{0}}^{e}=\left\langle E_{1},E_{2},E_{3},H_{3},F_{3}\right\rangle .$
To determine $\mathfrak{g}_{\bar{1}}^{e}$, assume $x=\sum a_{i,j,k}v_{i,j,k}\in\mathfrak{g}_{\bar{1}}^{e}$
where $a_{i,j,k}\in\mathbb{C},i,j,k\in\{\pm1\}$. By calculating $[E_{1}+E_{2},x]=0$,
we obtain that $a_{1,1,k}$ are arbitrary, $a_{1,-1,k}=-a_{-1,1,k}$
for $k=\pm1$ and $a_{-1,-1,k}=0$. Hence a basis of $\mathfrak{g}_{\bar{1}}^{e}$
is $\{v_{1,1,1},v_{1,1,-1},v_{1,-1,1}-v_{-1,1,1},v_{1,-1,-1}-v_{-1,1,-1}\}$.
Therefore, we have that $\dim\mathfrak{g}^{e}=5+4=9$.

By computing commutator relations between basis elements for $\mathfrak{g}^{e}(0)$,
we deduce that $\mathfrak{g}^{e}(0)\cong\mathfrak{osp}(1|2)$ according
to Lemma \ref{lem:osp(1,2)} where $F_{3},v_{1,-1,-1}-v_{-1,1,-1},H_{3},v_{1,-1,1}-v_{-1,1,1},E_{3}$
correspond to $u_{-2},u_{-1},u_{0},u_{1},u_{2}$ in Lemma \ref{lem:osp(1,2)}.
Moreover, we obtain that $\mathfrak{g}^{e}(2)=V^{\mathfrak{osp}}(0)\oplus V^{\mathfrak{osp}}(1)$.
Hence, we have that $\mathfrak{z}=\mathfrak{z}(0)\oplus\mathfrak{z}(2)\subseteq\left(\mathfrak{g}^{e}(0)\right)^{\mathfrak{g}^{e}(0)}\oplus\left(\mathfrak{g}^{e}(2)\right)^{\mathfrak{g}^{e}(0)}=\left\langle E_{1}+E_{2}\right\rangle $.
We know that $e=E_{1}+E_{2}\in\mathfrak{z}$. Therefore, $\mathfrak{z}=\left\langle E_{1}+E_{2}\right\rangle $
and $\dim\mathfrak{z}=1$.

Next we look for the labelled Dynkin diagrams with respect to $e$.
We find an element $h=(H,H,0)$ such that $h$ belongs to an $\mathfrak{sl}(2)$-triple
$\{e,h,f\}$ in $\mathfrak{g}_{\bar{0}}$. By calculating the ad$h$-eigenvalue
on each root vector, we have that roots in $\mathfrak{g}(>0)$ are
$\{2\beta_{1},2\beta_{2},\beta_{1}+\beta_{2}+\beta_{3},\beta_{1}+\beta_{2}-\beta_{3}\}$
and roots in $\mathfrak{g}(0)$ are $\Phi(0)=\{\pm2\beta_{3},i\beta_{1}-i\beta_{2}+k\beta_{3}:i,k=\pm1\}$.
Hence, we have that $\mathfrak{g}(0)\cong\mathfrak{gl}(1\mid2)$ and
there are three systems of simple roots of $\mathfrak{g}(0)$: $\varPi_{1}(0)=\{-\beta_{1}+\beta_{2}-\beta_{3},2\beta_{3}\}$,
$\varPi_{2}(0)=\{2\beta_{3},\beta_{1}-\beta_{2}-\beta_{3}\}$ and
$\varPi_{3}(0)=\{-\beta_{1}+\beta_{2}+\beta_{3},\beta_{1}-\beta_{2}+\beta_{3}\}$
up to conjugacy. By extending $\varPi_{i}(0)$ to $\varPi$ for $i=1,2,3$,
we get three systems of positive roots $\Phi_{i}^{+}$ and simple
roots $\varPi_{i}$. Therefore, there are three conjugacy classes
of Borel subalgebras such that $\mathfrak{b}=\mathfrak{h}\oplus\bigoplus_{\alpha\in\Phi^{+}}\mathfrak{g}_{\alpha}\subseteq\bigoplus_{j\geq0}\mathfrak{g}(j)$.
Hence, the systems of simple roots are: 

$\varPi_{1}=\{\alpha_{1}=2\beta_{1},\alpha_{2}=-\beta_{1}+\beta_{2}-\beta_{3},\alpha_{3}=2\beta_{3}\}$.
We compute $\mu_{12}=1$ and $\mu_{23}=\alpha$ using Formula (\ref{eq:lines-=0003BC}).
Therefore, the labelled Dynkin diagram with repect to $\varPi_{1}$
is the Dynkin diagram in Figure 4.1 with labels $2,0,0$.

$\varPi_{2}=\{\alpha_{1}=2\beta_{3},\alpha_{2}=\beta_{1}-\beta_{2}-\beta_{3},\alpha_{3}=2\beta_{2}\}$.
We compute that $\mu_{12}=\alpha$ and $\mu_{23}=1+\alpha$ using
Formula (\ref{eq:lines-=0003BC}). Therefore, the labelled Dynkin
diagram with repect to $\varPi_{2}$ is the Dynkin diagram in Figure
4.3 with labels $0,0,2$.

$\varPi_{3}=\{\alpha_{1}=-\beta_{1}+\beta_{2}+\beta_{3},\alpha_{2}=\beta_{1}-\beta_{2}+\beta_{3},\alpha_{3}=\beta_{1}+\beta_{2}-\beta_{3}\}$.
We compute that $\mu_{12}=\alpha,\mu_{13}=1+\alpha$ and $\mu_{23}=2$
using Formula (\ref{eq:lines-=0003BC}). Therefore, the labelled Dynkin
diagram with repect to $\varPi_{3}$ is the Dynkin diagram in Figure
4.4 with labels $0,0,2$.
\end{singlespace}

\subsection{Analysis of results\label{subsec:Analysis-D(2,1;)}}

\begin{singlespace}
\noindent Note that $e=E_{1}+E_{2}$ is the only case in which the
corresponding labelled Dynkin diagram $\varDelta$ has no label equal
to $1$. For this case, we have $n_{2}(\varDelta)=1$, $\dim\mathfrak{z}(\mathfrak{g}^{e})=1$
and $\mathfrak{g}^{h}=\mathfrak{g}(0)\cong\mathfrak{gl}(1|2)$ by
\S4.3. Hence, $\mathfrak{g}^{h}$ has centre of dimension $1$. Therefore,
$\dim\mathfrak{z}(\mathfrak{g}^{h})=n_{2}(\varDelta)=\dim\mathfrak{z}(\mathfrak{g}^{e})$. 

\noindent In order to prove Theorem $2$, we only need to look at
the case $e=E_{1}+E_{2}$ as the remaining cases do not have labels
equal to $2$ so that the $2$-free core $\varDelta_{0}$ of $\varDelta$
is the same as $\varDelta$. For this case, we have $\dim\mathfrak{g}^{e}=9$,
$\dim\mathfrak{z}(\mathfrak{g}^{e})=1$ and $\varDelta_{0}=\varDelta_{\mathfrak{g}(0)}$.
Hence, we deduce that the corresponding Lie superalgebra $\mathfrak{g}_{0}=\mathfrak{sl}(2|1)$
and the nilpotent orbit $e_{0}$ with respect to $\varDelta_{0}$
is equal to $0$. Therefore, we have that $\dim\mathfrak{g}_{0}^{e_{0}}=8$
and $\dim\mathfrak{g}^{e}-\dim\mathfrak{g}_{0}^{e_{0}}=n_{2}(\varDelta)=1$.
Similarly, we have $\dim\mathfrak{z}(\mathfrak{g}^{e})-\dim\mathfrak{z}(\mathfrak{g}_{0}^{e_{0}})=n_{2}(\varDelta)=1$
because $\dim\mathfrak{z}(\mathfrak{g}_{0}^{e_{0}})=0$.
\end{singlespace}

\begin{singlespace}
Theorem 3 for $D(2,1;\alpha)$ can be obtained immediately from Table
\ref{tab:results in D(2,1)}, i.e. we have that $\dim\mathfrak{z}(\mathfrak{g}^{e})=\left\lceil \frac{1}{2}\sum_{i=1}^{3}a_{i}\right\rceil +\varepsilon$
where $\varepsilon=0$ for $e=0,E_{1},E_{1}+E_{2}$ and $\varepsilon=-1$
for $e=E_{1}+E_{2}+E_{3}$.
\end{singlespace}

\section{The Exceptional Lie superalgebra $G(3)$\label{sec:G(3)}}

\subsection{Structure of the Lie superalgebra $G(3)$\label{subsec:Structure-of-G(3)}}

Let $V_{2}=V$ where $V$ is defined in \S3. Let $G_{2}$ be the
exceptional Lie algebra and $V_{7}=\left\langle e_{3},e_{2},e_{1},e_{0},e_{-1},e_{-2},e_{-3}\right\rangle $
be its $7$-dimensional module. A construction of $G(3)$ can be found
in \cite[Chapter 4]{Musson2012} and we recall this construction below.
Recall that $G(3)=\mathfrak{g}=\mathfrak{g}_{\bar{0}}\oplus\mathfrak{g}_{\bar{1}}$
where 
\[
\mathfrak{g}_{\bar{0}}=\mathfrak{sl}(2)\oplus G_{2}\text{ and }\mathfrak{g}_{\bar{1}}=V_{2}\otimes V_{7}.
\]
 We view $G_{2}\subseteq\mathfrak{gl}(V_{7})$ and then $\mathfrak{g}_{\bar{0}}$
has a basis $\{E,H,F,h_{1},h_{2},x_{i},y_{i}:i=1,\dots,6\}$ where
\[
x_{1}=\begin{pmatrix}0 & -1 & 0 & 0 & 0 & 0 & 0\\
0 & 0 & 0 & 0 & 0 & 0 & 0\\
0 & 0 & 0 & 1 & 0 & 0 & 0\\
0 & 0 & 0 & 0 & -2 & 0 & 0\\
0 & 0 & 0 & 0 & 0 & 0 & 0\\
0 & 0 & 0 & 0 & 0 & 0 & 1\\
0 & 0 & 0 & 0 & 0 & 0 & 0
\end{pmatrix},\ x_{2}=\begin{pmatrix}0 & 0 & 0 & 0 & 0 & 0 & 0\\
0 & 0 & 1 & 0 & 0 & 0 & 0\\
0 & 0 & 0 & 0 & 0 & 0 & 0\\
0 & 0 & 0 & 0 & 0 & 0 & 0\\
0 & 0 & 0 & 0 & 0 & -1 & 0\\
0 & 0 & 0 & 0 & 0 & 0 & 0\\
0 & 0 & 0 & 0 & 0 & 0 & 0
\end{pmatrix},
\]
\[
y_{1}=\begin{pmatrix}0 & 0 & 0 & 0 & 0 & 0 & 0\\
-1 & 0 & 0 & 0 & 0 & 0 & 0\\
0 & 0 & 0 & 0 & 0 & 0 & 0\\
0 & 0 & 2 & 0 & 0 & 0 & 0\\
0 & 0 & 0 & -1 & 0 & 0 & 0\\
0 & 0 & 0 & 0 & 0 & 0 & 0\\
0 & 0 & 0 & 0 & 0 & 1 & 0
\end{pmatrix},\ y_{2}=\begin{pmatrix}0 & 0 & 0 & 0 & 0 & 0 & 0\\
0 & 0 & 0 & 0 & 0 & 0 & 0\\
0 & 1 & 0 & 0 & 0 & 0 & 0\\
0 & 0 & 0 & 0 & 0 & 0 & 0\\
0 & 0 & 0 & 0 & 0 & 0 & 0\\
0 & 0 & 0 & 0 & -1 & 0 & 0\\
0 & 0 & 0 & 0 & 0 & 0 & 0
\end{pmatrix},
\]
\[
h_{1}=\text{diag}(1,-1,2,0,-2,1,-1)\ \text{and}\ h_{2}=\text{diag}(0,1,-1,0,1,-1,0),
\]
and $x_{3}=[x_{1},x_{2}]$, $x_{4}=[x_{1},x_{3}]$, $x_{5}=[x_{1},x_{4}]$,
$x_{6}=[x_{5},x_{2}]$. The remaining negative root vectors in the
basis of $G_{2}$ can be generated by $y_{1}$ and $y_{2}$ in a similar
way. A basis of $\mathfrak{g}_{\bar{1}}$ is $\{v_{i}\otimes e_{j}:i=\pm1,j=0,\pm1,\pm2,\pm3\}$.

\begin{singlespace}
We know that $\mathfrak{g}_{\bar{0}}$ is a Lie algebra and the bracket
$\left[\cdotp,\cdotp\right]:\mathfrak{g}_{\bar{0}}\times\mathfrak{g}_{\bar{1}}\rightarrow\mathfrak{g}_{\bar{1}}$
is given by $[x+y,u\otimes w]=xu\otimes w+u\otimes yw$ for $x\in\mathfrak{sl}(2),y\in G_{2},u\in V_{2}$
and $w\in V_{7}$. The bracket $\left[\cdotp,\cdotp\right]:\mathfrak{g}_{\bar{1}}\times\mathfrak{g}_{\bar{1}}\rightarrow\mathfrak{g}_{\bar{0}}$
is given in \cite[Theorem 4.4.5]{Musson2012}. For $v_{i},v_{k}\in V_{2},e_{j},e_{l}\in V_{7}$
, we have 
\begin{equation}
[v_{i}\otimes e_{j},v_{k}\otimes e_{l}]=\psi_{2}(v_{i},v_{k})p_{7}(e_{j},e_{l})+\psi_{7}(e_{j},e_{l})p_{2}(v_{i},v_{k}),\label{eq:G(3)}
\end{equation}
where $\psi_{2}$ is a non-degenerate skew-symmetric bilinear form
on $V_{2}$ such that $\psi_{2}(v_{1},v_{-1})=1$ and $p_{2}:V_{2}\times V_{2}\rightarrow\mathfrak{sl}(2)$
is given by $p_{2}(x,y)(z)=4\left(\psi_{1}(y,z)x-\psi_{1}(z,x)y\right)$.
We calculate that $p_{2}(v_{1},v_{-1})=-4H$, $p_{2}(v_{1},v_{1})=8E$
and $p_{2}(v_{-1},v_{-1})=-8F$. The mappings $\psi_{7}$ and $p_{7}$
are defined in \cite[Theorem 4.4.5]{Musson2012}. We can explicitly
calculate that $\psi_{7}(e_{j},e_{-j})=2$, $\psi_{7}(e_{0},e_{0})=-1$
and $\psi_{7}(e_{i},e_{j})=0$ if $i\neq j$. As in \cite[Subsection 4.7.9]{Musson2012},
we can calculate that $p_{7}(e_{-3},e_{3})=-8h_{1}-12h_{2}$. According
to the graded Jacobi identity, we have that 
\begin{align*}
[x_{1},[v_{1}\otimes e_{-3},v_{-1}\otimes e_{3}]] & =[v_{1}\otimes e_{-3},[x_{1},v_{-1}\otimes e_{3}]]+[[x_{1},v_{1}\otimes e_{-3}],v_{-1}\otimes e_{3}]\\
 & =[v_{1}\otimes e_{-3},0]+[v_{1}\otimes e_{-2},v_{-1}\otimes e_{3}]=p_{7}(e_{-2},e_{3}).
\end{align*}
Moreover, $[x_{1},[v_{1}\otimes e_{-3},v_{-1}\otimes e_{3}]]=[x_{1},-8H-8h_{1}-12h_{2}]=4x_{1}$.
Therefore, we deduce that $p_{7}(e_{-2},e_{3})=4x_{1}$. Using similar
methods we obtain the explicit mapping $p_{7}:V_{7}\times V_{7}\rightarrow G_{2}$
which is given in the following table:
\end{singlespace}

\begin{singlespace}
\noindent 
\begin{table}[H]
\begin{singlespace}
\begin{longtable}[c]{|>{\centering}m{1cm}||>{\centering}m{1.7cm}||>{\centering}m{1.5cm}||>{\centering}m{1.6cm}||>{\centering}m{1.1cm}||>{\centering}m{1.8cm}||>{\centering}m{1.7cm}||>{\centering}m{1.5cm}|}
\hline 
 & $e_{3}$ & $e_{2}$ & $e_{1}$ & $e_{0}$ & $e_{-1}$ & $e_{-2}$ & $e_{-3}$\tabularnewline
\hline 
\hline 
$e_{3}$ & $0$ & $-2x_{6}$ & $2x_{5}$ & $2x_{4}$ & $-4x_{3}$ & $-4x_{1}$ & $8h_{1}+12h_{2}$\tabularnewline
\hline 
\hline 
$e_{2}$ & $2x_{6}$ & $0$ & $-2x_{4}$ & $-4x_{3}$ & $12x_{2}$ & $4h_{1}+12h_{2}$ & $-4y_{1}$\tabularnewline
\hline 
\hline 
$e_{1}$ & $-2x_{5}$ & $2x_{4}$ & $0$ & $4x_{1}$ & $4h_{1}$ & $12y_{2}$ & $4y_{3}$\tabularnewline
\hline 
\hline 
$e_{0}$ & $-2x_{4}$ & $4x_{3}$ & $-4x_{1}$ & $0$ & $4y_{1}$ & $4y_{3}$ & $2y_{4}$\tabularnewline
\hline 
\hline 
$e_{-1}$ & $4x_{3}$ & $-12x_{2}$ & $-4h_{1}$ & $-4y_{1}$ & $0$ & $-2y_{4}$ & $-2y_{5}$\tabularnewline
\hline 
\hline 
$e_{-2}$ & $4x_{1}$ & $-4h_{1}-12h_{2}$ & $-12y_{2}$ & $-4y_{3}$ & $2y_{4}$ & $0$ & $-2y_{6}$\tabularnewline
\hline 
\hline 
$e_{-3}$ & $-8h_{1}-12h_{2}$ & $4y_{1}$ & $-4y_{3}$ & $-2y_{4}$ & $2y_{5}$ & $2y_{6}$ & $0$\tabularnewline
\hline 
\end{longtable}
\end{singlespace}

\caption{\label{tab:p7 G2}$p_{7}:V_{7}\times V_{7}\rightarrow G_{2}$}

\end{table}

\end{singlespace}

\subsection{Root system and Dynkin diagrams for $G(3)$\label{subsec:root-system-G(3)}}

\begin{singlespace}
\noindent We follow the description of the root system of $\mathfrak{g}=G(3)$
given in \cite[Chapter 4]{Musson2012}. Let $\mathfrak{h}=\left\langle H,h_{1},h_{2}\right\rangle $
be the Cartan subalgebra of $\mathfrak{g}$. Note that the roots of
$G(3)$ can be expressed in terms of $\delta,\varepsilon_{1,}\varepsilon_{2},\varepsilon_{3}\in\mathfrak{h}^{*}$
where $\varepsilon_{1}+\varepsilon_{2}+\varepsilon_{3}=0$. The root
system $\Phi=\Phi_{\bar{0}}\cup\Phi_{\bar{1}}$ is given by
\[
\Phi_{\bar{0}}=\{\pm2\delta,\varepsilon_{i}-\varepsilon_{j},\pm\varepsilon_{i}:1\leq i,j\leq3\}\text{ and }\Phi_{\bar{1}}=\{\pm\delta\pm\varepsilon_{i},\pm\delta:1\leq i\leq3\}
\]

\noindent where the bilinear form $(\cdotp,\cdotp)$ on $\mathfrak{h}^{*}$
is defined by $(\delta,\delta)=2$, $(\varepsilon_{i},\varepsilon_{j})=1-3\delta_{ij}$,
and $(\delta,\varepsilon_{i})=0$.
\end{singlespace}

\begin{singlespace}
The table below lists all roots together with corresponding root vectors:
\end{singlespace}
\noindent \begin{center}
\begin{tabular}{|c||>{\centering}m{0.6cm}|>{\centering}m{0.6cm}||>{\centering}m{1.3cm}||>{\centering}m{0.6cm}||>{\centering}m{0.6cm}||>{\centering}m{1.3cm}||c||c||c|}
\hline 
Roots & $2\delta$ & $\varepsilon_{1}$ & $\varepsilon_{2}-\varepsilon_{1}$ & $\varepsilon_{2}$ & $-\varepsilon_{3}$ & $\varepsilon_{1}-\varepsilon_{3}$ & $\varepsilon_{2}-\varepsilon_{3}$ & $i\delta-\varepsilon_{3}$ & $i\delta+\varepsilon_{j}$\tabularnewline
\hline 
Root vectors & $E$ & $x_{1}$ & $x_{2}$ & $x_{3}$ & $x_{4}$ & $x_{5}$ & $x_{6}$ & $v_{i}\otimes e_{3}$ & $v_{i}\otimes e_{j}$\tabularnewline
\hline 
\hline 
Roots & $-2\delta$ & $-\varepsilon_{1}$ & $\varepsilon_{1}-\varepsilon_{2}$ & $-\varepsilon_{2}$ & $\varepsilon_{3}$ & $\varepsilon_{3}-\varepsilon_{1}$ & $\varepsilon_{3}-\varepsilon_{2}$ & $i\delta+\varepsilon_{3}$ & $i\delta$\tabularnewline
\hline 
Root vectors & $F$ & $y_{1}$ & $y_{2}$ & $y_{3}$ & $y_{4}$ & $y_{5}$ & $y_{6}$ & $v_{i}\otimes e_{-3}$ & $v_{i}\otimes e_{0}$\tabularnewline
\hline 
\end{tabular}
\par\end{center}

\begin{singlespace}
\noindent where $i\in\{1,-1\}$, $j\in\{2,1,-1,-2\}$ and we define
$\varepsilon_{j}=-\varepsilon_{-j}$ for $j<0$. We further deduce
that the odd roots $\pm\delta\pm\varepsilon_{i}$ are isotropic and
$\pm\delta$ are non-isotropic.
\end{singlespace}

\begin{singlespace}
The following table covers all possible Dynkin diagrams with respect
to different systems of simple roots based on \cite[Section 2.19]{Frappat1996}.
\end{singlespace}

\begin{longtable}[c]{|>{\centering}m{6cm}||>{\centering}p{5cm}|}
\caption{Dynkin diagrams for $G(3)$}
\tabularnewline
\endfirsthead
\hline 
Simple systems $\varPi=\{\alpha_{1},\alpha_{2},\alpha_{3}\}$ & Dynkin diagrams\tabularnewline
\hline 
\hline 
\begin{singlespace}
\noindent $\{\delta+\varepsilon_{3},\varepsilon_{1},\varepsilon_{2}-\varepsilon_{1}\}$
\end{singlespace}
 & Figure 5.1
\noindent \centering{}\includegraphics{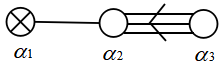}\tabularnewline
\hline 
\hline 
$\{-\delta-\varepsilon_{3},\delta-\varepsilon_{2},\varepsilon_{2}-\varepsilon_{1}\}$ & Figure 5.2
\noindent \centering{}\includegraphics{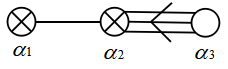}\tabularnewline
\hline 
\hline 
$\{\delta,-\delta+\varepsilon_{1},\varepsilon_{2}-\varepsilon_{1}\}$ & Figure 5.3
\noindent \centering{}\includegraphics{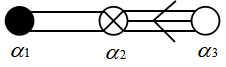}\tabularnewline
\hline 
\hline 
$\{\varepsilon_{1},-\delta+\varepsilon_{2},\delta-\varepsilon_{1}\}$ & Figure 5.4
\noindent \centering{}\includegraphics{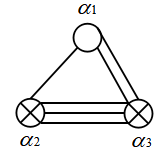}\tabularnewline
\hline 
\end{longtable}

\subsection{Centres of centralizers of nilpotent elements $e$ in $G(3)$ and
labelled Dynkin diagrams with respect to $e$\label{subsec:centres-of-centralizers--G(3)}}

\begin{singlespace}
\noindent Let $e=e_{\mathfrak{sl}(2)}+e_{G_{2}}\in\mathfrak{g}_{\bar{0}}$
be nilpotent where $e_{\mathfrak{sl}(2)}\in\mathfrak{sl}(2)$ and
$e_{G_{2}}\in G_{2}$. According to \cite[Section 11]{Lawther2008},
we know that representatives of nilpotent orbits in $\mathfrak{sl}(2)$
are $0,E$ and representatives of nilpotent orbits in $G_{2}$ are
$0,x_{2},x_{1},x_{2}+x_{5},x_{1}+x_{2}$ up to the adjoint action
of $G=\mathrm{SL}_{2}(\mathbb{C})\times K$ where $K$ is the Lie
group of type $G_{2}$. Hence, there are in total $10$ possibilities
for $e$. It is clear that $\mathfrak{sl}(2)^{E}=\langle E\rangle$
and $\mathfrak{sl}(2)^{0}=\mathfrak{sl}(2)$. We give basis elements
for $\mathfrak{g}^{e}$ and $\mathfrak{z}(\mathfrak{g}^{e})$ and
the labelled Dynkin diagrams $\varDelta$ with respect to $e$ in
Table \ref{tab:G(3)}. Note that the numbers in the column labelled
``$\varDelta$'' represent labels $a_{i}$ corresponding to $\alpha_{i}$
for $i=1,2,3$ in labelled Dynkin diagram with respect to $e$.
\end{singlespace}
\begin{singlespace}
\noindent \begin{center}
\begin{longtable}[c]{|>{\centering}m{1.5cm}||>{\centering}m{6cm}||>{\centering}m{3cm}||>{\centering}m{3.5cm}|}
\caption{\label{tab:G(3)}$\mathfrak{g}^{e}$, $\mathfrak{z}(\mathfrak{g}^{e})$
and $\varDelta$ for $\mathfrak{g}=G(3)$}
\tabularnewline
\endfirsthead
\hline 
$e$ & $\mathfrak{g}^{e}$ & $\mathfrak{z}(\mathfrak{g}^{e})$ & $\varDelta$\tabularnewline
\hline 
\hline 
$E+(x_{1}+x_{2})$ & $\langle E,x_{1}+x_{2},x_{6},v_{1}\otimes e_{3},v_{1}\otimes e_{2}+v_{-1}\otimes e_{3}\rangle$ & $\langle e,x_{6},v_{1}\otimes e_{3}\rangle$ & Figure 5.3: $1,1,2$\tabularnewline
\hline 
\hline 
$E+x_{2}$ & \begin{singlespace}
\noindent \centering{}$\langle E,2h_{1}+3h_{2},x_{2},y_{1},x_{3},x_{6},y_{5},x_{4},y_{4},v_{1}\otimes e_{2},v_{1}\otimes e_{-1},v_{1}\otimes e_{3},v_{1}\otimes e_{0},v_{1}\otimes e_{-3},v_{1}\otimes e_{1}-v_{-1}\otimes e_{2},v_{1}\otimes e_{-2}+v_{-1}\otimes e_{-1}\rangle$
\end{singlespace}
 & $\langle e\rangle$ & Figure 5.1: $0,0,1$

Figure 5.2: $0,0,1$

Figure 5.4: $0,0,1$\tabularnewline
\hline 
\hline 
$E+x_{1}$ & $\langle E,x_{1},x_{5},y_{2},x_{6},y_{6},h_{1}+2h_{2},v_{1}\otimes e_{1},v_{1}\otimes e_{3},v_{1}\otimes e_{-2},v_{1}\otimes e_{0}-v_{-1}\otimes e_{1},v_{1}\otimes e_{2}+v_{-1}\otimes e_{3},v_{1}\otimes e_{-3}-v_{-1}\otimes e_{-2}\rangle$ & $\langle e\rangle$ & Figure 5.2: $1,0,0$

Figure 5.4: $1,0,0$\tabularnewline
\hline 
\hline 
$E+(x_{2}+x_{5})$ & \begin{singlespace}
\noindent \centering{}$\langle E,x_{6},x_{3},x_{4},x_{2}+x_{5},v_{1}\otimes e_{3},v_{1}\otimes e_{2},v_{1}\otimes e_{0},6v_{-1}\otimes e_{3}-v_{1}\otimes e_{-1},v_{1}\otimes e_{1}-v_{-1}\otimes e_{2}\rangle$
\end{singlespace}
 & $\langle e,x_{6}\rangle$ & Figure 5.4: $0,1,1$\tabularnewline
\hline 
\hline 
$E$ & \begin{singlespace}
\noindent \centering{}$\left\langle E\right\rangle \oplus\mathrm{G}_{2}\oplus\langle v_{1}\otimes e_{3},v_{1}\otimes e_{2},v_{1}\otimes e_{1},v_{1}\otimes e_{0},v_{1}\otimes e_{-1},v_{1}\otimes e_{-2},v_{1}\otimes e_{-3}\rangle$
\end{singlespace}
 & $\langle e\rangle$ & Figure 5.1: $1,0,0$\tabularnewline
\hline 
\hline 
$x_{1}+x_{2}$ & \begin{singlespace}
\noindent \centering{}$\langle E,H,F,x_{1}+x_{2},x_{6},v_{1}\otimes e_{3},v_{-1}\otimes e_{3}\rangle$
\end{singlespace}
 & $\langle e,x_{6}\rangle$ & Figure 5.3: $0,2,2$\tabularnewline
\hline 
\hline 
$x_{2}$ & \begin{singlespace}
\noindent \centering{}$\langle E,H,F,2h_{1}+3h_{2},x_{2},y_{1},x_{3},x_{6},y_{5},x_{4},y_{4},v_{1}\otimes e_{2},v_{-1}\otimes e_{2},v_{1}\otimes e_{-1},v_{-1}\otimes e_{-1},v_{1}\otimes e_{-3},v_{-1}\otimes e_{-3},v_{1}\otimes e_{3},v_{-1}\otimes e_{3},v_{1}\otimes e_{0},v_{-1}\otimes e_{0}\rangle$
\end{singlespace}
 & $\langle e\rangle$ & Figure 5.3: $0,0,1$\tabularnewline
\hline 
\hline 
$x_{1}$ & $\langle E,H,F,x_{1},x_{5},y_{2},x_{6},y_{6},h_{1}+2h_{2},v_{1}\otimes e_{1},v_{-1}\otimes e_{1},v_{1}\otimes e_{3},v_{-1}\otimes e_{3},v_{1}\otimes e_{-2},v_{-1}\otimes e_{-2}\rangle$ & $\langle e\rangle$ & Figure 5.3: $0,1,0$\tabularnewline
\hline 
\hline 
$x_{2}+x_{5}$ & \begin{singlespace}
\noindent \centering{}$\langle E,H,F,x_{6},x_{3},x_{4},x_{2}+x_{5},v_{1}\otimes e_{3},v_{-1}\otimes e_{3},v_{1}\otimes e_{2},v_{-1}\otimes e_{2},v_{1}\otimes e_{0},v_{-1}\otimes e_{0}\rangle$
\end{singlespace}
 & $\langle e,x_{6}\rangle$ & Figure 5.3: $0,0,2$

Figure 5.4: $0,2,0$\tabularnewline
\hline 
\hline 
$0$ & $\mathfrak{g}$ & $\{0\}$ & Figures 5.1, 5.2, 5.3, 5.4. All labels are zeros.\tabularnewline
\hline 
\end{longtable}
\par\end{center}
\end{singlespace}

\begin{singlespace}
For each nilpotent element $e$, we find a semisimple element $h$
such that $h$ lies in an $\mathfrak{sl}(2)$-triple in $\mathfrak{g}_{\bar{0}}$
that contains $e$. We also calculate the $\mathfrak{g}^{e}(0)$-module
structure on each $\mathfrak{g}^{e}(j)$ for $j>0$ in the table below.
Let $V^{\mathfrak{sl}}(j)$ be an $\mathfrak{sl}(2)$-module with
highest weight $j$ and $V^{\mathfrak{osp}}(j)$ be an $\mathfrak{osp}(1|2)$-module
with highest weight $j$. Note that for $e=x_{2}$, the $\mathfrak{g}^{e}(0)$-module
structure on $\mathfrak{g}^{e}(j)$ is not included as it requires
the construction of $\mathfrak{osp}(3|2)$ representations. 
\end{singlespace}
\begin{singlespace}
\noindent \begin{center}
\begin{longtable}[c]{|>{\centering}m{2cm}||>{\centering}m{2.5cm}||>{\centering}m{2cm}||>{\centering}m{7cm}|}
\caption{\label{tab:g^e(0)-G(3)}The $\mathfrak{g}^{e}(0)$-module structure
on $\mathfrak{g}^{e}(j)$ for $j>0$ }
\tabularnewline
\endfirsthead
\hline 
$e$ & $h$ & $\mathfrak{g}^{e}(0)$ & $\mathfrak{g}^{e}(j),j>0$\tabularnewline
\hline 
\hline 
$E+(x_{1}+x_{2})$ & $H+(6h_{1}+10h_{2})$ & $0$ & $\dim\mathfrak{g}^{e}(10)=\dim\mathfrak{g}^{e}(7)=\dim\mathfrak{g}^{e}(5)=1$,$\dim\mathfrak{g}^{e}(2)=2$\tabularnewline
\hline 
\hline 
$E+x_{2}$ & $H+h_{2}$ & $\mathfrak{osp}(1|2)$ & $\mathfrak{g}^{e}(1)=V^{\mathfrak{osp}}(3)$, $\mathfrak{g}^{e}(2)=V^{\mathfrak{osp}}(0)\oplus V^{\mathfrak{osp}}(1)$\tabularnewline
\hline 
\hline 
$E+x_{1}$ & $H+h_{1}$ & $\mathfrak{osp}(1|2)$ & $\mathfrak{g}^{e}(1)=\mathfrak{g}^{e}(3)=V^{\mathfrak{osp}}(0)$,
$\mathfrak{g}^{e}(2)=V^{\mathfrak{osp}}(0)\oplus V^{\mathfrak{osp}}(1)$\tabularnewline
\hline 
\hline 
$E+(x_{2}+x_{5})$ & $H+(2h_{1}+4h_{2})$ & $0$ & $\dim\mathfrak{g}^{e}(4)=1,\dim\mathfrak{g}^{e}(3)=2,$ $\dim\mathfrak{g}^{e}(2)=4,\dim\mathfrak{g}^{e}(1)=3$\tabularnewline
\hline 
\hline 
$E$ & $H$ & $G_{2}$ & $\mathfrak{g}^{e}(1)=V_{7},\mathfrak{g}^{e}(2)=\left\langle E\right\rangle $\tabularnewline
\hline 
\hline 
$x_{1}+x_{2}$ & $6h_{1}+10h_{2}$ & $\mathfrak{sl}(2)$ & $\mathfrak{g}^{e}(2)=\mathfrak{g}^{e}(10)=V^{\mathfrak{sl}}(0)$,
$\mathfrak{g}^{e}(6)=V^{\mathfrak{sl}}(1)$\tabularnewline
\hline 
\hline 
$x_{2}$ & $h_{2}$ & $\mathfrak{osp}(3|2)$ & Omitted.\tabularnewline
\hline 
\hline 
$x_{1}$ & $h_{1}$ & $\mathfrak{sl}(2)\oplus\mathfrak{sl}(2)$ & $\mathfrak{g}^{e}(1)=V^{\mathfrak{sl}}(1)\otimes V^{\mathfrak{sl}}(1),\mathfrak{g}^{e}(3)=V^{\mathfrak{sl}}(0)\otimes V^{\mathfrak{sl}}(1),$$\mathfrak{g}^{e}(2)=\left(V^{\mathfrak{sl}}(0)\otimes V^{\mathfrak{sl}}(0)\right)\oplus\left(V^{\mathfrak{sl}}(1)\otimes V^{\mathfrak{sl}}(0)\right).$\tabularnewline
\hline 
\hline 
$x_{2}+x_{5}$ & $2h_{1}+4h_{2}$ & $\mathfrak{osp}(1|2)$ & $\mathfrak{g}^{e}(2)=V^{\mathfrak{osp}}(0)\oplus V^{\mathfrak{osp}}(1)\oplus V^{\mathfrak{osp}}(1)$,
$\mathfrak{g}^{e}(4)=V^{\mathfrak{osp}}(0)$\tabularnewline
\hline 
\hline 
$0$ & $0$ & $\mathfrak{g}$ & $0$\tabularnewline
\hline 
\end{longtable}
\par\end{center}
\end{singlespace}

\begin{singlespace}
In the remaining part of this subsection, we explain explicit calculations
for finding $\mathfrak{g}^{e}$ and $\mathfrak{z}(\mathfrak{g}^{e})$
and obtain the corresponding labelled Dynkin diagrams for nilpotent
element $E+x_{2}$. The results of remaining cases are obtained using
the same approach.

When $e=E+x_{2}$, we already know that $\mathfrak{sl}(2)^{E}=\langle E\rangle$,
now we are going to work out $\mathrm{G}_{2}^{x_{2}}$. Observe that
$h_{G_{2}}=\text{diag}(0,1,-1,0,1,-1,0)=h_{2}$ belongs to an $\mathfrak{sl}(2)$-triple
$\{e_{G_{2}},h_{G_{2}},f_{G_{2}}\}$ in $G_{2}$ containing $e_{G_{2}}=x_{2}$.
Then we can work out the ad$h_{G_{2}}$-eigenspaces with non-negative
eigenvalues of $G_{2}$:
\end{singlespace}
\noindent \begin{center}
\begin{tabular}{|c|c|c|c|}
\hline 
Eigenvalues of $h_{G_{2}}$ & $0$ & $1$ & $2$\tabularnewline
\hline 
\hline 
Eigenvectors & $h_{1},h_{2},x_{4},y_{4}$  & $y_{1},x_{3},x_{6},y_{5}$ & $x_{2}$\tabularnewline
\hline 
\end{tabular}
\par\end{center}

\begin{singlespace}
\noindent This demonstrates that $G_{2}^{x_{2}}\cong G_{2}^{x_{2}}(0)\oplus G_{2}^{x_{2}}(1)\oplus G_{2}^{x_{2}}(1)$.
It is clear that $G_{2}^{x_{2}}(2)=\langle x_{2}\rangle$ and $G_{2}^{x_{2}}(1)=\langle y_{1},x_{3},x_{6},y_{5}\rangle$.
Note that $G_{2}^{x_{2}}(0)$ has dimension $3$. Since $\left[2h_{1}+3h_{2},x_{2}\right]=0$
and $\left[x_{4},x_{2}\right]=0=\left[y_{4},x_{2}\right]$, we have
that $G_{2}^{x_{2}}(0)=\langle x_{4},y_{4},2h_{1}+3h_{2}\rangle$.
Therefore, we have that $\mathfrak{g}_{\bar{0}}^{e}$ has a basis
$\{E,2h_{1}+3h_{2},x_{4},y_{4},x_{2},y_{1},x_{3},x_{6},y_{5}\}$.
\end{singlespace}

\begin{singlespace}
Now we calculate $\mathfrak{g}_{\bar{1}}^{e}$. We look at the $\mathfrak{sl}(2)$-triple
$\{e,h,f\}\subseteq\mathfrak{g}_{\bar{0}}$ and work out all non-negative
ad$h$-eigenspaces in $\mathfrak{g}_{\bar{1}}$.
\end{singlespace}
\begin{singlespace}
\noindent \begin{center}
\begin{tabular}{|>{\centering}m{3cm}||c|>{\centering}m{3cm}|>{\centering}m{4cm}|}
\hline 
ad$h$-eigenvalues & $2$ & $1$ & $0$\tabularnewline
\hline 
\hline 
ad$h$-eigenvectors & $v_{1}\otimes e_{2}$, $v_{1}\otimes e_{-1}$ & $v_{1}\otimes e_{3}$, $v_{1}\otimes e_{0}$, $v_{1}\otimes e_{-3}$ & $v_{1}\otimes e_{1}$, $v_{1}\otimes e_{-2}$, $v_{-1}\otimes e_{2}$,
$v_{-1}\otimes e_{-1}$\tabularnewline
\hline 
\end{tabular}
\par\end{center}
\end{singlespace}

\begin{singlespace}
\noindent The above table implies that $\mathfrak{g}_{\bar{1}}^{e}\cong\mathfrak{g}_{\bar{1}}^{e}(0)\oplus\mathfrak{g}_{\bar{1}}^{e}(1)\oplus\mathfrak{g}_{\bar{1}}^{e}(2)$
where $\mathfrak{g}_{\bar{1}}^{e}(2)$ has a basis $\{v_{1}\otimes e_{2},v_{1}\otimes e_{-1}\}$
and $\mathfrak{g}_{\bar{1}}^{e}(1)$ has a basis $\{v_{1}\otimes e_{3},v_{1}\otimes e_{0},v_{1}\otimes e_{-3}\}$.
To determine $\mathfrak{g}_{\bar{1}}^{e}(0)$, we need to find elements
of the form $x=a_{1,1}v_{1}\otimes e_{1}+a_{1,-2}v_{1}\otimes e_{-2}+a_{-1,2}v_{-1}\otimes e_{2}+a_{-1,-1}$$v_{-1}\otimes e_{-1}$
that are centralized by $e$. Then $\left[e,x\right]=0$ gives that
$a_{1,1}=-a_{-1,2}$ and $a_{1,-2}=a_{-1,-1}$. Hence $\mathfrak{g}_{\bar{1}}^{e}(0)$
has a basis $\{v_{1}\otimes e_{1}-v_{-1}\otimes e_{2},v_{1}\otimes e_{-2}+v_{-1}\otimes e_{-1}\}$.
Therefore, $\mathfrak{g}_{\bar{1}}^{e}$ has a basis $\{v_{1}\otimes e_{2},v_{1}\otimes e_{-1},v_{1}\otimes e_{3},v_{1}\otimes e_{0},$$v_{1}\otimes e_{-3},v_{1}\otimes e_{1}-v_{-1}\otimes e_{2},v_{1}\otimes e_{-2}+v_{-1}\otimes e_{-1}\}$.
In conclusion, we have that $\dim\mathfrak{g}^{e}=9+7=16$.
\end{singlespace}

\begin{singlespace}
By computing commutator relations between basis elements for $\mathfrak{g}^{e}(0)$,
we deduce that $\mathfrak{g}^{e}(0)\cong\mathfrak{osp}(1|2)$ according
to Lemma \ref{lem:osp(1,2)} where $y_{4},v_{1}\otimes e_{-2}+v_{-1}\otimes e_{-1},2h_{1}+3h_{2},v_{1}\otimes e_{1}-v_{-1}\otimes e_{2},x_{4}$
correspond to $u_{-2},u_{-1},u_{0},u_{1},u_{2}$ in Lemma \ref{lem:osp(1,2)}.
Moreover, we obtain that $\mathfrak{g}^{e}(1)=V^{\mathfrak{osp}}(3)$
and $\mathfrak{g}^{e}(2)=V^{\mathfrak{osp}}(0)\oplus V^{\mathfrak{osp}}(1)$.
Hence, we have that 
\[
\mathfrak{z}=\mathfrak{z}(0)\oplus\mathfrak{z}(1)\oplus\mathfrak{z}(2)\subseteq\left(\mathfrak{g}^{e}(0)\right)^{\mathfrak{g}^{e}(0)}\oplus\left(\mathfrak{g}^{e}(1)\right)^{\mathfrak{g}^{e}(0)}\oplus\left(\mathfrak{g}^{e}(2)\right)^{\mathfrak{g}^{e}(0)}=\langle E+x_{2}\rangle.
\]
 Note that $E+x_{2}\in\mathfrak{z}$, therefore $\mathfrak{z}=\langle E+x_{2}\rangle$
and it has dimension $1$.

Next we look at the labelled Dynkin diagrams with respect to $e$.
We obtain that roots in $\mathfrak{g}(>0)$ are $\{2\delta,\varepsilon_{2}-\varepsilon_{1},-\varepsilon_{1},\varepsilon_{2},\varepsilon_{2}-\varepsilon_{3,}\varepsilon_{3}-\varepsilon_{1},\delta+\varepsilon_{2},\delta-\varepsilon_{1},\delta-\varepsilon_{3},\delta,\delta+\varepsilon_{3}\}$
and roots in $\mathfrak{g}(0)$ are $\Phi(0)=\{\pm\varepsilon_{3},\pm(\delta+\varepsilon_{1}),\pm(\delta-\varepsilon_{2})\}$.
Hence, there are three systems of simple roots of $\mathfrak{g}(0)$:
$\varPi_{1}(0)=\{-\varepsilon_{3},-\delta-\varepsilon_{1}\}$, $\varPi_{2}(0)=\{-\delta+\varepsilon_{2},\delta+\varepsilon_{1}\}$
and $\varPi_{3}(0)=\{\delta-\varepsilon_{2},-\varepsilon_{3}\}$ up
to conjugacy. By extending $\varPi_{i}(0)$ to $\varPi$ for $i=1,2,3$,
we get three systems of positive roots $\Phi_{i}^{+}$ and simple
roots $\varPi_{i}$ and thus there are three conjugacy classes of
Borel subalgebras satisfying $\mathfrak{b}=\mathfrak{h}\oplus\bigoplus_{\alpha\in\Phi^{+}}\mathfrak{g}_{\alpha}\subseteq\bigoplus_{j\geq0}\mathfrak{g}(j)$.
Hence, the systems of simple roots are:

$\varPi_{1}=\{-\varepsilon_{3},-\delta-\varepsilon_{1},\delta+\varepsilon_{3}\}$.
We compute $\mu_{12}=1$, $\mu_{13}=2$ and $\mu_{23}=3$ using Formula
(\ref{eq:lines-=0003BC}). Therefore, the corresponding labelled Dynkin
diagram is the Dynkin diagram in Figure 5.4 with labels $0,0,1$. 
\end{singlespace}

$\varPi_{2}=\{-\delta+\varepsilon_{2},\delta+\varepsilon_{1},\varepsilon_{3}-\varepsilon_{1}\}$.
We compute $\mu_{12}=1$ and $\mu_{23}=3$ using Formula (\ref{eq:lines-=0003BC}).
Therefore, the corresponding labelled Dynkin diagram is the Dynkin
diagram in Figure 5.2 with labels $0,0,1$.

$\varPi_{3}=\{\delta-\varepsilon_{2},-\varepsilon_{3},\varepsilon_{3}-\varepsilon_{1}\}$.
We compute $\mu_{12}=1$ and $\mu_{23}=3$ using Formula (\ref{eq:lines-=0003BC}).
Therefore, the corresponding labelled Dynkin diagram is the Dynkin
diagram in Figure 5.1 with labels $0,0,1$.

\subsection{Analysis of results\label{subsec:Analysis-of-results-G3}}

\begin{singlespace}
\noindent Let $\mathfrak{h}=\left\langle H,h_{1},h_{2}\right\rangle \subseteq\mathfrak{g}$.
Denote a simple root system for $\mathfrak{g}^{h}$ by $\varPi_{h}$.
In order to prove Theorem 1 for $G(3),$ we consider two cases in
which the corresponding labelled Dynkin diagram has no label equal
to $1$. They are $e=x_{1}+x_{2}$ and $e=x_{2}+x_{5}$.
\end{singlespace}

\begin{singlespace}
When $e=x_{1}+x_{2}\in G(3)$, we have $\dim\mathfrak{z}(\mathfrak{g}^{e})=2$
and $n_{2}(\varDelta)=2$. Note that $\mathfrak{g}^{h}$ is generated
by root vectors $e_{\pm2\delta}$, $e_{\pm\delta}$ and $\mathfrak{h}$.
Hence, we can find a simple root system $\varPi_{h}=\{\delta\}$ and
thus $\mathfrak{g}^{h}=\mathfrak{z}(\mathfrak{g}^{h})\oplus\mathfrak{osp}(1|2)$
according to \cite[Subsection 3.4.1]{Musson2012}. Then $\mathfrak{z}(\mathfrak{g}^{h})=\{t\in\mathfrak{h}:\delta(t)=0\}$.
Hence, $\dim\mathfrak{z}(\mathfrak{g}^{h})=2=n_{2}(\varDelta)=\dim\mathfrak{z}(\mathfrak{g}^{e})$.

When $e=x_{2}+x_{5}\in G(3)$, we have $\dim\mathfrak{z}(\mathfrak{g}^{e})=2$
but $n_{2}(\varDelta)=1$. Note that $\mathfrak{g}^{h}$ is generated
by root vectors $e_{\pm2\delta},e_{\pm\delta},e_{\pm\varepsilon_{1}},e_{\pm(\delta+\varepsilon_{1})},e_{\pm(\delta-\varepsilon_{1})}$
and $\mathfrak{h}$, thus we can find a simple root system $\varPi_{h}=\{\varepsilon_{1},\delta-\varepsilon_{1}\}$
and thus $\mathfrak{g}^{h}=\mathfrak{z}(\mathfrak{g}^{h})\oplus\mathfrak{osp}(3|2)$
according to \cite[Subsection 3.4.1]{Musson2012}. Then $\mathfrak{z}(\mathfrak{g}^{h})=\{t\in\mathfrak{h}:\varepsilon_{1}(t)=(\delta-\varepsilon_{1})(t)=0\}$.
Hence, $\dim\mathfrak{z}(\mathfrak{g}^{h})=1=n_{2}(\varDelta)$ but
$\dim\mathfrak{z}(\mathfrak{g}^{h})\neq\dim\mathfrak{z}(\mathfrak{g}^{e})$.
We will further discuss this case in \S5.5 and complete verification
of Theorem 1 for $G(3)$.

In order to prove Theorem 2 for $G(3)$, we look at three cases below
and the remaining cases do not have labels equal to $2$ so that the
$2$-free core $\varDelta_{0}$ of $\varDelta$ is the same as $\varDelta$.

When $e=E+(x_{1}+x_{2})\in G(3)$, we have $\dim\mathfrak{g}^{e}=5$,
$\dim\mathfrak{z}(\mathfrak{g}^{e})=3$ and $\varDelta$ is given
in Table \ref{tab:G(3)}. From \cite[Subsection 3.4.1]{Musson2012}
we have that the Lie superalgebra corresponds to $\varDelta_{0}$
is $\mathfrak{g}_{0}=\mathfrak{osp}(3|2)$. An explicit explanation
of construction of $\mathfrak{osp}(3|2)$ can be found in \cite[Section 2.3]{Musson2012}.
We choose the representative of nilpotent orbit $e_{0}\in\mathfrak{osp}(3|2)_{\bar{0}}$
to be
\end{singlespace}

\begin{singlespace}
\noindent 
\[
\begin{pmatrix}0 & 1 & 0 & 0 & 0\\
0 & 0 & -1 & 0 & 0\\
0 & 0 & 0 & 0 & 0\\
0 & 0 & 0 & 0 & 1\\
0 & 0 & 0 & 0 & 0
\end{pmatrix},
\]
 thus $e$ has the Jordan type $(3|2)$. Hence, we obtain that $\dim\mathfrak{g}_{0}^{e_{0}}=4$
according to \cite[Subsection 3.2.2]{Hoyt2012}. Furthermore, we calculate
that $\dim\mathfrak{z}(\mathfrak{g}_{0}^{e_{0}})=2$. Therefore, $\dim\mathfrak{g}^{e}-\dim\mathfrak{g}_{0}^{e_{0}}=n_{2}(\varDelta)=1$
and $\dim\mathfrak{z}(\mathfrak{g}^{e})-\dim\mathfrak{z}(\mathfrak{g}_{0}^{e_{0}})=n_{2}(\varDelta)=1$
for this case.
\end{singlespace}

\begin{singlespace}
When $e=x_{1}+x_{2}\in G(3)$, we have $\dim\mathfrak{g}^{e}=7$,
$\dim\mathfrak{z}(\mathfrak{g}^{e})=2$ and $\varDelta$ is given
in Table \ref{tab:G(3)}. According to \cite[Subsection 3.4.1]{Musson2012},
we have that the Lie superalgebra corresponding to $\varDelta_{0}$
is $\mathfrak{g}_{0}=\mathfrak{osp}(1|2)$ and $e_{0}=0$. Thus we
know that $\mathfrak{g}_{0}^{e_{0}}=\mathfrak{g}_{0}$ and $\mathfrak{z}(\mathfrak{g}_{0}^{e_{0}})=\mathfrak{z}(\mathfrak{g}_{0})=\{0\}$.
Hence, we have $\dim\mathfrak{g}_{0}^{e_{0}}=\dim\mathfrak{g}_{0}=5$
and $\dim\mathfrak{z}(\mathfrak{g}_{0}^{e_{0}})=\dim\mathfrak{z}(\mathfrak{g}_{0})=0$.
Therefore, $\dim\mathfrak{g}^{e}-\dim\mathfrak{g}_{0}^{e_{0}}=n_{2}(\varDelta)=2$
and $\dim\mathfrak{z}(\mathfrak{g}^{e})-\dim\mathfrak{z}(\mathfrak{g}_{0}^{e_{0}})=n_{2}(\varDelta)=2$
for this case. 

When $e=x_{2}+x_{5}\in G(3)$, we have $\dim\mathfrak{g}^{e}=13$,
$\dim\mathfrak{z}(\mathfrak{g}^{e})=2$ and $\varDelta$ is given
in Table \ref{tab:G(3)}. From \cite[Subsection 3.4.1]{Musson2012}
we have that the Lie superalgebra corresponding to $\varDelta_{0}$
is $\mathfrak{g}_{0}=\mathfrak{osp}(3|2)$ and $e_{0}=0$. Similar
to the above case, we have $\dim\mathfrak{g}_{0}^{e_{0}}=\dim\mathfrak{g}_{0}=12$
and $\dim\mathfrak{z}(\mathfrak{g}_{0}^{e_{0}})=\dim\mathfrak{z}(\mathfrak{g}_{0})=0$.
Therefore, $\dim\mathfrak{g}^{e}-\dim\mathfrak{g}_{0}^{e_{0}}=n_{2}(\varDelta)=1$
and $\dim\mathfrak{z}(\mathfrak{g}^{e})-\dim\mathfrak{z}(\mathfrak{g}_{0}^{e_{0}})=2\neq n_{2}(\varDelta)$
for this case. We will further discuss this case in \S5.5.
\end{singlespace}

\subsection{Adjoint action on $G(3)$\label{subsec:Adjoint-action-of-G(3)}}

\begin{singlespace}
Let $K$ be a simple Lie group of type $G_{2}$, then $\mathrm{Lie}(K)$
is the Lie algebra of type $G_{2}$. Let $G=\mathrm{SL_{2}(\mathbb{C})}\times K$.
For a nilpotent element $e=e_{\mathfrak{sl}}+e_{G_{2}}\in\mathfrak{g}_{\bar{0}}$
where $e_{\mathfrak{sl}}\in\mathfrak{sl}(2)$ and $e_{G_{2}}\in\mathrm{Lie}(K)$,
we determine $\left(\mathfrak{z}(\mathfrak{g}^{e})\right)^{G^{e}}$
in this subsection. Write $(K^{e_{G_{2}}})^{\circ}$ for the connected
component of $K^{e_{G_{2}}}$ containing the identity and let $K^{e_{G_{2}}}(0)=K^{h_{G_{2}}}\cap K^{e_{G_{2}}}$.
An explicit structure of $K^{e_{G_{2}}}$ has been given in \cite[Section 11]{Lawther2008}.
\end{singlespace}

For $e_{\mathfrak{sl}}=0$, the centralizer in $\mathrm{\mathrm{SL_{2}(\mathbb{C})}}$
is connected. Thus it suffices to only look at the action of $K^{e}$
on $\mathfrak{z}\subseteq\mathfrak{z}_{\bar{0}}$. In this case, we
know that $G^{e}/(G^{e})^{\circ}\cong K^{e}/(K^{e})^{\circ}\cong K^{e}(0)/(K^{e}(0))^{\circ}$.
When $e=x_{1},x_{2},x_{1}+x_{2}$, we have that $K^{e}(0)/(K^{e}(0))^{\circ}=1$
by \cite[Section 11]{Lawther2008} and thus $\left(\mathfrak{z}(\mathfrak{g}^{e})\right)^{G^{e}}=\mathfrak{z}(\mathfrak{g}^{e})$. 

When $e=x_{2}+x_{5}$, the component group does not centralize $x_{6}$
according to \cite[page 73]{Lawther2008}. Hence, we deduce that $\left(\mathfrak{z}(\mathfrak{g}^{e})\right)^{G^{e}}=\left\langle e\right\rangle $.
Therefore, based on \S5.4, we have that $\dim\left(\mathfrak{z}(\mathfrak{g}^{e})\right)^{G^{e}}=\dim\mathfrak{z}(\mathfrak{g}^{h})=1$
and $\dim\left(\mathfrak{z}(\mathfrak{g}^{e})\right)^{G^{e}}-\dim\left(\mathfrak{z}(\mathfrak{g}_{0}^{e_{0}})\right)^{G_{0}^{e_{0}}}=n_{2}(\varDelta)=1$
for this case.

For $e_{\mathfrak{sl}}=E$, we know that $G^{e}=(\{\pm1\}\ltimes R^{E})\times K^{e_{G_{2}}}$
where $R^{E}$ is a connected normal subgroup of $G^{e}$. When $e=E+x_{2}$,
$E+x_{1}$, $E+(x_{2}+x_{5})$, we have $\mathfrak{z}\subseteq\mathfrak{z}_{\bar{0}}$.
We know that $\pm1\in\mathrm{\mathrm{SL_{2}(\mathbb{C})}}^{E}$ act
trivially on $\mathfrak{z}_{\bar{0}}$ and $\left(\mathfrak{z}(\mathfrak{g}^{e})\right)^{K^{e_{G_{2}}}}=\left\langle e\right\rangle $
by \cite[Section 11]{Lawther2008}. Hence, we have that $\left(\mathfrak{z}(\mathfrak{g}^{e})\right)^{G^{e}}=\mathfrak{z}(\mathfrak{g}^{e})$
for $e=E+x_{2}$, $E+x_{1}$ and $\left(\mathfrak{z}(\mathfrak{g}^{e})\right)^{G^{e}}=\left\langle e\right\rangle $
for $e=E+(x_{2}+x_{5})$. 

When $e=E+(x_{1}+x_{2})$, the component group of $\mathrm{SL(2)}^{^{E}}$
has order $2$ and we consider the element $g=-1\in G^{e}/(G^{e})^{\circ}$.
We know that $g$ acts trivially on $\mathfrak{z}_{\bar{0}}$, thus
$e,x_{6}\in\left(\mathfrak{z}(\mathfrak{g}^{e})\right)^{g}$. However,
$v_{1}\otimes e_{3}\notin\left(\mathfrak{z}(\mathfrak{g}^{e})\right)^{g}$
since the action of $g$ on $v_{1}\otimes e_{3}$ sends it to $-v_{1}\otimes e_{3}$.
Hence, we have that $\left(\mathfrak{z}(\mathfrak{g}^{e})\right)^{G^{e}}\subseteq\left(\mathfrak{z}(\mathfrak{g}^{e})\right)^{g}=\left\langle e,x_{6}\right\rangle $.
Therefore, we have that $\left(\mathfrak{z}(\mathfrak{g}^{e})\right)^{G^{e}}=\left\langle e,x_{6}\right\rangle $.
Based on \S5.4, we know that $G_{0}=\mathrm{O}_{3}(\mathbb{C})\times\mathrm{Sp}_{2}(\mathbb{C})$.
By considering the element $g'=-1\in G_{0}^{e_{0}}/(G_{0}^{e_{0}})^{\circ}$,
we obtain that $\left(\mathfrak{z}(\mathfrak{g}_{0}^{e_{0}})\right)^{G_{0}^{e_{0}}}\subseteq\left(\mathfrak{z}(\mathfrak{g}_{0}^{e_{0}})\right)^{g'}=\left\langle e_{0}\right\rangle $.
Therefore, we deduce that $\dim\left(\mathfrak{z}(\mathfrak{g}^{e})\right)^{G^{e}}-\dim\left(\mathfrak{z}(\mathfrak{g}_{0}^{e_{0}})\right)^{G_{0}^{e_{0}}}=n_{2}(\varDelta)=1$
for this case. The above argument completes the proof of Theorems
1 and 2 for $G(3)$.

\begin{singlespace}
By combining results in Table \ref{tab:G(3)}, we have that $\dim\left(\mathfrak{z}(\mathfrak{g}^{e})\right)^{G^{e}}=\left\lceil \frac{1}{2}\sum_{i=1}^{3}a_{i}\right\rceil $
which proves the statement of Theorem 3 for $G(3)$.
\end{singlespace}
\begin{singlespace}

\section{The Exceptional Lie superalgebra $F(4)$\label{sec:F(4)}}
\end{singlespace}

\subsection{$\mathfrak{so}(V,\beta)$ embedded into $\mathrm{C}(V,\beta)$\label{subsec:Clifford-algebra}}

\begin{singlespace}
\noindent Let $V$ be a finite-dimensional complex vector space with
a basis $\{v_{i}:i=1,\dots,m\}$ and a symmetric bilinear form $\beta$
and let $\mathrm{C}(V,\beta)$ be the Clifford algebra for $(V,\beta)$.
An explicit definition of the Clifford algebra can be found in \cite[Chapter 6]{Goodman2009}.
Recall that for any $x,y\in V$, we have that $\{x,y\}=\beta(x,y)1$
for $x,y\in V$ where $\{x,y\}=xy+yx$ is the anticommutator of $x,y$.
Then $\mathrm{C}(V,\beta)$ is spanned by $1$ and the products $v_{i_{1}}\dots v_{i_{l}}$
for $1\leq i_{1}<i_{2}<\dots<i_{l}\leq m$.
\end{singlespace}

\begin{singlespace}
For $u,w\in V$, we define $R_{u,w}\in\mathrm{End}(V)$ by $R_{u,w}(v)=\beta(w,v)u-\beta(u,v)w$.
For any $x,y\in V$, we can check that 
\[
\beta(R_{u,w}(x),y)=\beta(w,x)\beta(u,y)-\beta(u,x)\beta(w,y)=-\beta(x,R_{u,w}(y)).
\]
Hence, we know that $R_{u,w}\in\mathfrak{so}(V,\beta)$. The linear
transformations $R_{u,w}$ for $u,w\in V$ form a basis for $\mathfrak{so}(V,\beta)$,
see \cite[Lemma 6.2.1]{Goodman2009}. We only consider the case when
$\dim V=7$ in this section. Now let $\dim V=7$, there exists a decomposition
$V=W\oplus\left\langle e_{0}\right\rangle \oplus W^{*}$ where $W,W^{*}$
is a pair of dual maximal isotropic subspaces of $V$ corresponding
to $\beta$ and $\{e_{1},e_{2},e_{3}\}$ (resp. $\{e_{-1},e_{-2},e_{-3}\}$)
is a basis for $W$ (resp. $W^{*}$). The basis is chosen such that
$\beta(e_{i},e_{-j})=\delta_{ij}$ for $i,j\in\{1,2,3\}$, $\beta(e_{0},e_{0})=2$
and $\beta(e_{0},W)=\beta(e_{0},W^{*})=0$. Then 
\[
R_{e_{i},e_{-j}}=e_{i,j}-e_{-j,-i},R_{e_{i},e_{j}}=e_{i,-j}-e_{j,-i}\text{ for }i<j,R_{e_{-i},e_{-j}}=e_{-i,j}-e_{-j,i}\text{ for }i>j,
\]
\[
\text{and }R_{e_{i},e_{0}}=2e_{i,0}-e_{0,-i},\ R_{e_{-i},e_{0}}=2e_{-i,0}-e_{0,i}
\]
 form a basis for $\mathfrak{so}(V,\beta)$ where $e_{i,j}$ is the
elementary transformation which sends $e_{i}$ to $e_{j}$ and the
rest of basis vectors to $0$.

Next we determine representatives of nilpotent orbits in $\mathfrak{so}(7)$
using the above notation. According to \cite[Section 1.6]{Jantzen2004a},
we know that any nilpotent orbit in $\mathfrak{so}(7)$ has a Jordan
type $\lambda\in\{(7),(5,1^{2}),(3^{2},1),(3,2^{2}),(3,1^{4}),(2^{2},1^{3}),(1^{7})\}$.
By using the orthogonal Dynkin pyramid of $\lambda$ that is defined
in \cite[Section 6]{Elashvili2005}, we are able to give a representative
of each nilpotent orbit using matrices. Let us fix the corresponding
representatives of each of the nilpotent orbits to be as in Table
\ref{tab:nilpotent ele}. For each nilpotent element $e_{\mathfrak{so}}\in\mathfrak{so}(7)$,
we also give a semisimple element $h$ such that there is an $\mathfrak{sl}(2)$-triple
$\{e_{\mathfrak{so}},h_{\mathfrak{so}},f_{\mathfrak{so}}\}\subseteq\mathfrak{so}(7)$.
\end{singlespace}
\begin{singlespace}
\noindent \begin{center}
\begin{longtable}[c]{|>{\centering}m{1.5cm}||>{\centering}m{5.4cm}||>{\centering}m{6.6cm}|}
\caption{\label{tab:nilpotent ele}Nilpotent orbits in $\mathfrak{so}(7)$}
\tabularnewline
\endfirsthead
\hline 
Jordan types & nilpotent element $e_{\mathfrak{so}}$  & semisimple element $h_{\mathfrak{so}}$\tabularnewline
\hline 
\hline 
$(7)$ & $e_{(7)}=R_{e_{1},e_{-2}}+R_{e_{2},e_{-3}}+R_{e_{3},e_{0}}$ & $h_{(7)}=6R_{e_{1},e_{-1}}+4R_{e_{2},e_{-2}}+2R_{e_{3},e_{-3}}$\tabularnewline
\hline 
\hline 
$(5,1^{2})$ & $e_{(5,1^{2})}=R_{e_{1},e_{-2}+R_{e_{2},e_{0}}}$ & $h_{(5,1^{2})}=4R_{e_{1},e_{-1}}+2R_{e_{2},e_{-2}}$\tabularnewline
\hline 
\hline 
$(3^{2},1)$ & $e_{(3^{2},1)}=R_{e_{1},e_{-3}}+R_{e_{2},e_{3}}$ & $h_{(3^{2},1)}=2R_{e_{1},e_{-1}}+2R_{e_{2},e_{-2}}$\tabularnewline
\hline 
\hline 
$(3,2^{2})$ & $e_{(3,2^{2})}=R_{e_{1},e_{0}}+R_{e_{2},e_{3}}$ & $h_{(3,2^{2})}=2R_{e_{1},e_{-1}}+R_{e_{2},e_{-2}}+R_{e_{3},e_{-3}}$\tabularnewline
\hline 
\hline 
$(3,1^{4})$ & $e_{(3,1^{4})}=R_{e_{1},e_{0}}$ & $h_{(3,1^{4})}=2R_{e_{1},e_{-1}}$\tabularnewline
\hline 
\hline 
$(2^{2},1^{3})$ & $e_{(2^{2},1^{3})}=R_{e_{1},e_{2}}$ & $h_{(2^{2},1^{3})}=R_{e_{1},e_{-1}}+R_{e_{2},e_{-2}}$\tabularnewline
\hline 
\hline 
$(1^{7})$ & $e_{(1^{7})}=0$ & $h_{(1^{7})}=0$\tabularnewline
\hline 
\end{longtable}
\par\end{center}
\end{singlespace}

\begin{singlespace}
For each nilpotent orbit $e_{\mathfrak{so}}\in\mathfrak{so}(7)$ and
any element $x\in\mathfrak{so}(7)$, by calculating $\left[e_{\mathfrak{so}},x\right]=0$
we obtain centralizers $\mathfrak{so}(7)^{e_{\mathfrak{so}}}$ of
each $e_{\mathfrak{so}}\in\mathfrak{so}(7)$ in the following table:

\begin{longtable}[c]{|>{\centering}m{1.9cm}|>{\centering}m{10.3cm}||>{\centering}m{1.9cm}|}
\caption{\label{tab:Centralizers-so7}$\mathfrak{so}(7)^{e_{\mathfrak{so}}}$
of nilpotent orbits $e_{\mathfrak{so}}\in\mathfrak{so}(7)$}
\tabularnewline
\endfirsthead
\hline 
$e_{\mathfrak{so}}\in\mathfrak{so}(7)$ & $\mathfrak{so}(7)^{e_{\mathfrak{so}}}$  & $\dim\mathfrak{so}(7)^{e_{\mathfrak{so}}}$ \tabularnewline
\hline 
\hline 
$e_{(7)}$ & $\langle e_{(7)},R_{e_{1},e_{0}}-2R_{e_{2},e_{3}},R_{e_{1},e_{2}}\rangle$ & $3$\tabularnewline
\hline 
\hline 
$e_{(5,1^{2})}$ & $\langle e_{(5,1^{2})},R_{e_{1},e_{-3}},R_{e_{3},e_{-3}},R_{e_{1},e_{3}},R_{e_{1},e_{2}}\rangle$ & $5$\tabularnewline
\hline 
\hline 
$e_{(3^{2},1)}$ & $\langle e_{(3^{2},1)},R_{e_{1},e_{-1}}-R_{e_{2},e_{-2}}+R_{e_{3},e_{-3}},R_{e_{2},e_{-3}},R_{e_{2},e_{0}},R_{e_{1},e_{0}},R_{e_{1},e_{3}},R_{e_{1},e_{2}}\rangle$ & $7$\tabularnewline
\hline 
\hline 
$e_{(3,2^{2})}$ & \begin{singlespace}
\noindent $\langle R_{e_{1},e_{0}},R_{e_{2},e_{3}},R_{e_{2},e_{-2}}-R_{e_{3},e_{-3}},R_{e_{2},e_{-3}},R_{e_{3},e_{-2}},2R_{e_{1},e_{-3}}+R_{e_{2},e_{0}},-2R_{e_{1},e_{-2}}+R_{e_{3},e_{0}},R_{e_{1},e_{3}},R_{e_{1},e_{2}}\rangle$
\end{singlespace}
 & $9$\tabularnewline
\hline 
\hline 
$e_{(3,1^{4})}$ & $\langle e_{(3,1^{4})},R_{e_{2},e_{3}},R_{e_{2},e_{-3}},R_{e_{2},e_{-2}},R_{e_{3},e_{-3}},$
$R_{e_{3},e_{-2}},R_{e_{-3},e_{-2}},R_{e_{1},e_{2}},R_{e_{1},e_{3}},R_{e_{1},e_{-3}},R_{e_{1},e_{-2}}\rangle$ & $11$\tabularnewline
\hline 
\hline 
$e_{(2^{2},1^{3})}$ & $\langle e_{(2^{2},1^{3})},R_{e_{1},e_{-2}},R_{e_{1},e_{-1}}-R_{e_{2},e_{-2}},R_{e_{3},e_{-3}},R_{e_{2},e_{-1}},$
$R_{e_{-3},e_{0}},R_{e_{1},e_{3}},R_{e_{1},e_{0}},R_{e_{2},e_{-3}},R_{e_{3},e_{0}},R_{e_{2},e_{3}},R_{e_{1},e_{-3}},R_{e_{2},e_{0}}\rangle$ & $13$\tabularnewline
\hline 
\hline 
$e_{(1^{7})}$ & $\mathfrak{so}(7)$ & $21$\tabularnewline
\hline 
\end{longtable}
\end{singlespace}

\subsection{The spin representation of Lie algebra $\mathfrak{so}(7)$\label{subsec:A-spin-representation}}

\begin{singlespace}
\noindent According to \cite[Lemma 6.2.2]{Goodman2009}, there exists
an injective Lie algebra homomorphism $\varphi:\mathfrak{so}(7)\rightarrow C(V,\beta)$
such that $\varphi(R_{e_{i},e_{-j}})=e_{i}e_{-j}$ for $i\neq j$
and $\varphi(R_{e_{i},e_{-i}})=e_{i}e_{-i}-\frac{1}{2}$ for $i\neq0$.
Note that $\mathrm{C}(V,\beta)$ has a basis $\{e_{1}^{\delta_{1}}e_{2}^{\delta_{2}}e_{3}^{\delta_{3}}e_{0}^{\delta_{0}}e_{-3}^{\delta_{-3}}e_{-2}^{\delta_{-2}}e_{-1}^{\delta_{-1}}:\delta_{i}=0\text{ or }1\}$.
The multiplication in $\mathrm{C}(V,\beta)$ satisfies $e_{0}^{2}=1$,
$e_{i}^{2}=0$ for $i\neq0$, $e_{i}e_{j}=-e_{j}e_{i}$ for $i\neq-j$
and $e_{i}e_{-i}=-e_{-i}e_{i}+1$. Consider the subalgebra $D_{0,-}$
of $\mathrm{C}(V,\beta)$ defined by $D_{0,-}=\langle e_{0}^{\delta_{0}}e_{-3}^{\delta_{-3}}e_{-2}^{\delta_{-2}}e_{-1}^{\delta_{-1}}:\delta_{i}=0\text{ or }1\rangle.$
Let us define $S=C(V,\beta)\otimes_{D_{0,-}}\langle s\rangle$ where
$\langle s\rangle$ is the 1-dimensional $D_{0,-}$-module such that
$e_{-i}s=0\text{ for }i=1,2,3\text{ and }e_{0}s=s.$ In fact, $S$
is an $8$-dimensional representation for $C(V,\beta)$ with basis
$\{1\otimes s,e_{1}\otimes s,e_{2}\otimes s,e_{3}\otimes s,e_{1}e_{2}\otimes s,e_{1}e_{3}\otimes s,e_{2}e_{3}\otimes s,e_{1}e_{2}e_{3}\otimes s\}$.
We can restrict $S$ to be a representation of $\mathfrak{so}(7)=\mathfrak{so}(V,\beta)\subseteq C(V,\beta)$
and we call it the \textit{spin representation} for $\mathfrak{so}(7)$,
see \cite[Section 6.2.2]{Goodman2009}. In the remaining subsections,
we write the basis for $S$ as
\begin{equation}
\{s,e_{1}s,e_{2}s,e_{3}s,e_{1}e_{2}s,e_{1}e_{3}s,e_{2}e_{3}s,e_{1}e_{2}e_{3}s\}.\label{eq:basis V8}
\end{equation}
We sometimes denote basis elements $s,-e_{1}s,e_{2}s,-e_{3}s,e_{1}e_{2}s,e_{1}e_{3}s,e_{2}e_{3}s,e_{1}e_{2}e_{3}s$
of $V_{8}$ by $v_{---}$, $v_{+--}$, $v_{-+-}$, $v_{--+}$, $v_{++-}$,
$v_{+-+}$, $v_{-++}$, $v_{+++}$ respectively.
\end{singlespace}

\subsection{Construction of the Lie superalgebra $F(4)$\label{subsec:Construction-of-F(4)}}

\begin{singlespace}
\noindent In order to fully describe the construction of $F(4)$,
we first let $V_{2}=V$ where $V$ is defined in \S3. We define $\psi_{2}:V_{2}\times V_{2}\rightarrow\mathbb{C}$
to be a non-degenerate skew-symmetric bilinear form such that $\psi_{2}(v_{1},v_{-1})=1$.
We also define $p_{2}:V_{2}\times V_{2}\rightarrow\mathfrak{sl}(2)$
by $p_{2}(x,y)(z)=3(\psi_{2}(y,z)x-\psi_{2}(z,x)y)$ for $x,y,z\in V_{2}$.
We compute that $p_{2}(v_{1},v_{-1})=-3H,p_{2}(v_{1},v_{1})=6E$ and
$p_{2}(v_{-1},v_{-1})=-6F$. 
\end{singlespace}

\begin{singlespace}
Next let $V_{8}$ be the spin representation of $\mathfrak{so}(7)$
with a basis shown in (\ref{eq:basis V8}). Then let $\psi_{8}:V_{8}\times V_{8}\rightarrow\mathbb{C}$
be the non-degenerate symmetric bilinear form given by 
\[
\psi_{8}(v_{\sigma_{1},\sigma_{2},\sigma_{3}},v_{\sigma'_{1},\sigma'_{2},\sigma'_{3}})=\prod_{i=1}^{3}\delta_{\sigma_{i},-\sigma'_{i}}
\]
for $\sigma_{i},\sigma_{i}'\in\{+,-\}$, e.g. we have that $\psi_{8}(v_{+++},v_{+--})=0$
and $\psi_{8}(v_{+-+},v_{-+-})=1$. Define $p_{8}:V_{8}\times V_{8}\rightarrow\mathfrak{so}(7)$
to be the antisymmetric bilinear map gives explicitly on the basis
elements in Table \ref{tab:p8}. Note that values of $p_{8}(\cdotp,\cdotp)$
are calculated based on the assumption that 
\begin{equation}
(v_{+++},v_{++-})\longmapsto R_{e_{1},e_{2}}.\label{eq:v(+++)v(++-)=00003DRe1e2}
\end{equation}
 For example, by applying $R_{e_{3},e_{-2}}$ to both sides of (\ref{eq:v(+++)v(++-)=00003DRe1e2})
we get 
\[
(R_{e_{3},e_{-2}}v_{+++},v_{++-})+(v_{+++},R_{e_{3},e_{-2}}v_{++-})\longmapsto[R_{e_{3},e_{-2}},R_{e_{1},e_{2}}],
\]
 this implies that $(v_{+++},v_{+-+})\longmapsto R_{e_{1},e_{3}}$.
\end{singlespace}
\begin{singlespace}
\noindent \begin{center}
\begin{longtable}[c]{|>{\centering}m{1cm}|>{\centering}m{3cm}|>{\centering}m{3cm}|>{\centering}m{3cm}|>{\centering}m{3cm}|}
\caption{\label{tab:p8}$p_{8}:V_{8}\times V_{8}\rightarrow\mathfrak{so}(7)$}
\tabularnewline
\endfirsthead
\hline 
 & $v_{---}$ & $v_{+--}$ & $v_{-+-}$ & $v_{--+}$\tabularnewline
\hline 
\hline 
$v_{---}$ & \begin{onehalfspace}
\noindent \centering{}$0$
\end{onehalfspace}
 & \begin{onehalfspace}
\noindent \centering{}$-R_{e_{-3},e_{-2}}$
\end{onehalfspace}
 & \begin{onehalfspace}
\noindent \centering{}$-R_{e_{-3},e_{-1}}$
\end{onehalfspace}
 & \begin{onehalfspace}
\noindent \centering{}$-R_{e_{-2},e_{-1}}$
\end{onehalfspace}
\tabularnewline
\hline 
\hline 
$v_{+--}$ & \begin{onehalfspace}
\noindent \centering{}$R_{e_{-3},e_{-2}}$
\end{onehalfspace}
 & \begin{onehalfspace}
\noindent \centering{}$0$
\end{onehalfspace}
 & \begin{onehalfspace}
\noindent \centering{}$-\frac{1}{2}R_{e_{-3},e_{0}}$
\end{onehalfspace}
 & \begin{onehalfspace}
\noindent \centering{}$-\frac{1}{2}R_{e_{-2},e_{0}}$
\end{onehalfspace}
\tabularnewline
\hline 
\hline 
$v_{-+-}$ & \begin{onehalfspace}
\noindent \centering{}$R_{e_{-3},e_{-1}}$
\end{onehalfspace}
 & \begin{onehalfspace}
\noindent \centering{}$\frac{1}{2}R_{e_{-3},e_{0}}$
\end{onehalfspace}
 & \begin{onehalfspace}
\noindent \centering{}$0$
\end{onehalfspace}
 & \begin{onehalfspace}
\noindent \centering{}$-\frac{1}{2}R_{e_{-1},e_{0}}$
\end{onehalfspace}
\tabularnewline
\hline 
\hline 
$v_{--+}$ & \begin{onehalfspace}
\noindent \centering{}$R_{e_{-2},e_{-1}}$
\end{onehalfspace}
 & \begin{onehalfspace}
\noindent \centering{}$\frac{1}{2}R_{e_{-2},e_{0}}$
\end{onehalfspace}
 & \begin{onehalfspace}
\noindent \centering{}$\frac{1}{2}R_{e_{-1},e_{0}}$
\end{onehalfspace}
 & \begin{onehalfspace}
\noindent \centering{}$0$
\end{onehalfspace}
\tabularnewline
\hline 
\hline 
$v_{++-}$ & \begin{onehalfspace}
\noindent \centering{}$\frac{1}{2}R_{e_{-3},e_{0}}$
\end{onehalfspace}
 & \begin{onehalfspace}
\noindent \centering{}$R_{e_{1},e_{-3}}$
\end{onehalfspace}
 & \begin{onehalfspace}
\noindent \centering{}$-R_{e_{2},e_{-3}}$
\end{onehalfspace}
 & \begin{onehalfspace}
\noindent \centering{}$-\frac{1}{2}R_{e_{1},e_{-1}}-\frac{1}{2}R_{e_{2},e_{-2}}+\frac{1}{2}R_{e_{3},e_{-3}}$
\end{onehalfspace}
\tabularnewline
\hline 
\hline 
$v_{+-+}$ & \begin{onehalfspace}
\noindent \centering{}$-\frac{1}{2}R_{e_{-2},e_{0}}$
\end{onehalfspace}
 & \begin{onehalfspace}
\noindent \centering{}$-R_{e_{1},e_{-2}}$
\end{onehalfspace}
 & \begin{onehalfspace}
\noindent \centering{}$-\frac{1}{2}R_{e_{1},e_{-1}}+\frac{1}{2}R_{e_{2},e_{-2}}-\frac{1}{2}R_{e_{3},e_{-3}}$
\end{onehalfspace}
 & \begin{onehalfspace}
\noindent \centering{}$-R_{e_{3},e_{-2}}$
\end{onehalfspace}
\tabularnewline
\hline 
\hline 
$v_{-++}$ & \begin{onehalfspace}
\noindent \centering{}$\frac{1}{2}R_{e_{-1},e_{0}}$
\end{onehalfspace}
 & \begin{onehalfspace}
\noindent \centering{}$\frac{1}{2}R_{e_{1},e_{-1}}-\frac{1}{2}R_{e_{2},e_{-2}}-\frac{1}{2}R_{e_{3},e_{-3}}$
\end{onehalfspace}
 & \begin{onehalfspace}
\noindent \centering{}$-R_{e_{2},e_{-1}}$
\end{onehalfspace}
 & \begin{onehalfspace}
\noindent \centering{}$R_{e_{3},e_{-1}}$
\end{onehalfspace}
\tabularnewline
\hline 
\hline 
$v_{+++}$ & \begin{onehalfspace}
\noindent \centering{}$-\frac{1}{2}R_{e_{1},e_{-1}}-\frac{1}{2}R_{e_{2},e_{-2}}-\frac{1}{2}R_{e_{3},e_{-3}}$
\end{onehalfspace}
 & \begin{onehalfspace}
\noindent \centering{}$-\frac{1}{2}R_{e_{1},e_{0}}$
\end{onehalfspace}
 & \begin{onehalfspace}
\noindent \centering{}$\frac{1}{2}R_{e_{2},e_{0}}$
\end{onehalfspace}
 & \begin{onehalfspace}
\noindent \centering{}$-\frac{1}{2}R_{e_{3},e_{0}}$
\end{onehalfspace}
\tabularnewline
\hline 
 & \begin{onehalfspace}
\noindent \centering{}$v_{++-}$
\end{onehalfspace}
 & \begin{onehalfspace}
\noindent \centering{}$v_{+-+}$
\end{onehalfspace}
 & \begin{onehalfspace}
\noindent \centering{}$v_{-++}$
\end{onehalfspace}
 & \begin{onehalfspace}
\noindent \centering{}$v_{+++}$
\end{onehalfspace}
\tabularnewline
\hline 
$v_{---}$ & \begin{onehalfspace}
\noindent \centering{}$-\frac{1}{2}R_{e_{-3},e_{0}}$
\end{onehalfspace}
 & \begin{onehalfspace}
\noindent \centering{}$\frac{1}{2}R_{e_{-2},e_{0}}$
\end{onehalfspace}
 & \begin{onehalfspace}
\noindent \centering{}$-\frac{1}{2}R_{e_{-1},e_{0}}$
\end{onehalfspace}
 & \begin{onehalfspace}
\noindent \centering{}$\frac{1}{2}R_{e_{1},e_{-1}}+\frac{1}{2}R_{e_{2},e_{-2}}+\frac{1}{2}R_{e_{3},e_{-3}}$
\end{onehalfspace}
\tabularnewline
\hline 
$v_{+--}$ & \begin{onehalfspace}
\noindent \centering{}$-R_{e_{1},e_{-3}}$
\end{onehalfspace}
 & \begin{onehalfspace}
\noindent \centering{}$R_{e_{1},e_{-2}}$
\end{onehalfspace}
 & \begin{onehalfspace}
\noindent \centering{}$-\frac{1}{2}R_{e_{1},e_{-1}}+\frac{1}{2}R_{e_{2},e_{-2}}+\frac{1}{2}R_{e_{3},e_{-3}}$
\end{onehalfspace}
 & \begin{onehalfspace}
\noindent \centering{}$\frac{1}{2}R_{e_{1},e_{0}}$
\end{onehalfspace}
\tabularnewline
\hline 
$v_{-+-}$ & \begin{onehalfspace}
\noindent \centering{}$R_{e_{2},e_{-3}}$
\end{onehalfspace}
 & \begin{onehalfspace}
\noindent \centering{}$\frac{1}{2}R_{e_{1},e_{-1}}-\frac{1}{2}R_{e_{2},e_{-2}}+\frac{1}{2}R_{e_{3},e_{-3}}$ 
\end{onehalfspace}
 & \begin{onehalfspace}
\noindent \centering{}$R_{e_{2},e_{-1}}$
\end{onehalfspace}
 & \begin{onehalfspace}
\noindent \centering{}$-\frac{1}{2}R_{e_{2},e_{0}}$
\end{onehalfspace}
\tabularnewline
\hline 
$v_{--+}$ & \begin{onehalfspace}
\noindent \centering{}$\frac{1}{2}R_{e_{1},e_{-1}}+\frac{1}{2}R_{e_{2},e_{-2}}-\frac{1}{2}R_{e_{3},e_{-3}}$
\end{onehalfspace}
 & \begin{onehalfspace}
\noindent \centering{}$R_{e_{3},e_{-2}}$
\end{onehalfspace}
 & \begin{onehalfspace}
\noindent \centering{}$-R_{e_{3},e_{-1}}$
\end{onehalfspace}
 & \begin{onehalfspace}
\noindent \centering{}$\frac{1}{2}R_{e_{3},e_{0}}$
\end{onehalfspace}
\tabularnewline
\hline 
$v_{++-}$ & \begin{onehalfspace}
\noindent \centering{}$0$
\end{onehalfspace}
 & \begin{onehalfspace}
\noindent \centering{}$\frac{1}{2}R_{e_{1},e_{0}}$
\end{onehalfspace}
 & \begin{onehalfspace}
\noindent \centering{}$\frac{1}{2}R_{e_{2},e_{0}}$
\end{onehalfspace}
 & \begin{onehalfspace}
\noindent \centering{}$-R_{e_{1},e_{2}}$
\end{onehalfspace}
\tabularnewline
\hline 
$v_{+-+}$ & \begin{onehalfspace}
\noindent \centering{}$-\frac{1}{2}R_{e_{1},e_{0}}$
\end{onehalfspace}
 & \begin{onehalfspace}
\noindent \centering{}$0$
\end{onehalfspace}
 & \begin{onehalfspace}
\noindent \centering{}$\frac{1}{2}R_{e_{3},e_{0}}$
\end{onehalfspace}
 & \begin{onehalfspace}
\noindent \centering{}$-R_{e_{1},e_{3}}$
\end{onehalfspace}
\tabularnewline
\hline 
$v_{-++}$ & \begin{onehalfspace}
\noindent \centering{}$-\frac{1}{2}R_{e_{2},e_{0}}$
\end{onehalfspace}
 & \begin{onehalfspace}
\noindent \centering{}$-\frac{1}{2}R_{e_{3},e_{0}}$
\end{onehalfspace}
 & \begin{onehalfspace}
\noindent \centering{}$0$
\end{onehalfspace}
 & \begin{onehalfspace}
\noindent \centering{}$-R_{e_{2},e_{3}}$
\end{onehalfspace}
\tabularnewline
\hline 
$v_{+++}$ & \begin{onehalfspace}
\noindent \centering{}$R_{e_{1},e_{2}}$
\end{onehalfspace}
 & \begin{onehalfspace}
\noindent \centering{}$R_{e_{1},e_{3}}$
\end{onehalfspace}
 & \begin{onehalfspace}
\noindent \centering{}$R_{e_{2},e_{3}}$
\end{onehalfspace}
 & \begin{onehalfspace}
\noindent \centering{}$0$
\end{onehalfspace}
\tabularnewline
\hline 
\end{longtable}
\par\end{center}
\end{singlespace}

\begin{singlespace}
\noindent For example, we can read from this table that $p_{8}(v_{+++},v_{++-})=R_{e_{1},e_{2}}$.
\end{singlespace}

\begin{singlespace}
Together with the above definitions and notation, we are now able
to describe the structure of $F(4)$. Recall that the Lie superalgebra
of type $F(4)=\mathfrak{g}=\mathfrak{g}_{\bar{0}}\oplus\mathfrak{g}_{\bar{1}}$
where
\[
\mathfrak{g}_{\bar{0}}=\mathfrak{sl}(2)\oplus\mathfrak{so}(7)\text{ and }\mathfrak{g}_{\bar{1}}=V_{2}\otimes V_{8}.
\]

\end{singlespace}

\begin{singlespace}
\noindent We know that $\mathfrak{g}_{\bar{0}}$ is a Lie algebra
thus we have the bracket $\mathfrak{g}_{\bar{0}}\times\mathfrak{g}_{\bar{0}}\rightarrow\mathfrak{g}_{\bar{0}}$
and the bracket $\left[\cdotp,\cdotp\right]:\mathfrak{g}_{\bar{0}}\times\mathfrak{g}_{\bar{1}}\rightarrow\mathfrak{g}_{\bar{1}}$
is given by $[x+y,v_{2}\otimes v_{8}]=xv_{2}\otimes v_{8}+v_{2}\otimes yv_{8}$
for $x\in\mathfrak{sl}(2),y\in\mathfrak{so}(7)$, $v_{2}\in V_{2}$
and $v_{8}\in V_{8}$. Note that $Ev_{1}=0$, $Ev_{-1}=v_{1}$ and
$Hv_{i}=iv_{i}$ for $i\in\{\pm1\}$. The following table gives the
action of $\mathfrak{so}(7)$ on each basis element of $V_{8}$:

\noindent 
\begin{table}[H]
\begin{singlespace}
\noindent %
\begin{longtable}[c]{|c|c|c|c|c|c|c|c|c|}
\hline 
 & $s$ & $e_{1}s$ & $e_{2}s$ & $e_{3}s$ & $e_{1}e_{2}s$ & $e_{1}e_{3}s$ & $e_{2}e_{3}s$ & $e_{1}e_{2}e_{3}s$\tabularnewline
\hline 
\hline 
$R_{e_{1},e_{-1}}$ & $-\frac{1}{2}s$ & $\frac{1}{2}e_{1}s$ & $-\frac{1}{2}e_{2}s$ & $-\frac{1}{2}e_{3}s$ & $\frac{1}{2}e_{1}e_{2}s$ & $\frac{1}{2}e_{1}e_{3}s$ & $-\frac{1}{2}e_{2}e_{3}s$ & $\frac{1}{2}e_{1}e_{2}e_{3}s$\tabularnewline
\hline 
$R_{e_{1},e_{-2}}$ & $0$ & $0$ & $e_{1}s$ & $0$ & $0$ & $0$ & $e_{1}e_{3}s$ & $0$\tabularnewline
\hline 
$R_{e_{1},e_{-3}}$ & $0$ & $0$ & $0$ & $e_{1}s$ & $0$ & $0$ & $-e_{1}e_{2}s$ & $0$\tabularnewline
\hline 
$R_{e_{1},e_{0}}$ & $e_{1}s$ & $0$ & $-e_{1}e_{2}s$ & $-e_{1}e_{3}s$ & $0$ & $0$ & $e_{1}e_{2}e_{3}s$ & $0$\tabularnewline
\hline 
$R_{e_{1},e_{3}}$ & $e_{1}e_{3}s$ & $0$ & $-e_{1}e_{2}e_{3}s$ & $0$ & $0$ & $0$ & $0$ & $0$\tabularnewline
\hline 
$R_{e_{1},e_{2}}$ & $e_{1}e_{2}s$ & $0$ & $0$ & $e_{1}e_{2}e_{3}s$ & $0$ & $0$ & $0$ & $0$\tabularnewline
\hline 
$R_{e_{2},e_{-1}}$ & $0$ & $e_{2}s$ & $0$ & $0$ & $0$ & $e_{2}e_{3}s$ & $0$ & $0$\tabularnewline
\hline 
$R_{e_{2},e_{-2}}$ & $-\frac{1}{2}s$ & $-\frac{1}{2}e_{1}s$ & $\frac{1}{2}e_{2}s$ & $-\frac{1}{2}e_{3}s$ & $\frac{1}{2}e_{1}e_{2}s$ & $-\frac{1}{2}e_{1}e_{3}s$ & $\frac{1}{2}e_{2}e_{3}s$ & $\frac{1}{2}e_{1}e_{2}e_{3}s$\tabularnewline
\hline 
$R_{e_{2},e_{-3}}$ & $0$ & $0$ & $0$ & $e_{2}s$ & $0$ & $e_{1}e_{2}s$ & $0$ & $0$\tabularnewline
\hline 
$R_{e_{2},e_{0}}$ & $e_{2}s$ & $e_{1}e_{2}s$ & $0$ & $-e_{2}e_{3}s$ & $0$ & $-e_{1}e_{2}e_{3}s$ & $0$ & $0$\tabularnewline
\hline 
$R_{e_{2},e_{3}}$ & $e_{2}e_{3}s$ & $e_{1}e_{2}e_{3}s$ & $0$ & $0$ & $0$ & $0$ & $0$ & $0$\tabularnewline
\hline 
$R_{e_{3},e_{-1}}$ & $0$ & $e_{3}s$ & $0$ & $0$ & $-e_{2}e_{3}s$ & $0$ & $0$ & $0$\tabularnewline
\hline 
$R_{e_{3},e_{-2}}$ & $0$ & $0$ & $e_{3}s$ & $0$ & $e_{1}e_{3}s$ & $0$ & $0$ & $0$\tabularnewline
\hline 
$R_{e_{3},e_{-3}}$ & $-\frac{1}{2}s$ & $-\frac{1}{2}e_{1}s$ & $-\frac{1}{2}e_{2}s$ & $\frac{1}{2}e_{3}s$ & $-\frac{1}{2}e_{1}e_{2}s$ & $\frac{1}{2}e_{1}e_{3}s$ & $\frac{1}{2}e_{2}e_{3}s$ & $\frac{1}{2}e_{1}e_{2}e_{3}s$\tabularnewline
\hline 
$R_{e_{3},e_{0}}$ & $e_{3}s$ & $e_{1}e_{3}s$ & $e_{2}e_{3}s$ & $0$ & $e_{1}e_{2}e_{3}s$ & $0$ & $0$ & $0$\tabularnewline
\hline 
$R_{e_{-1},e_{0}}$ & $0$ & $-s$ & $0$ & $0$ & $e_{2}s$ & $e_{3}s$ & $0$ & $-e_{2}e_{3}s$\tabularnewline
\hline 
$R_{e_{-2},e_{0}}$ & $0$ & $0$ & $-s$ & $0$ & $-e_{1}s$ & $0$ & $e_{3}s$ & $e_{1}e_{3}s$\tabularnewline
\hline 
$R_{e_{-3},e_{0}}$ & $0$ & $0$ & $0$ & $-s$ & $0$ & $-e_{1}s$ & $-e_{2}s$ & $-e_{1}e_{2}s$\tabularnewline
\hline 
$R_{e_{-3},e_{-1}}$ & $0$ & $0$ & $0$ & $0$ & $0$ & $s$ & $0$ & $-e_{2}s$\tabularnewline
\hline 
$R_{e_{-3},e_{-2}}$ & $0$ & $0$ & $0$ & $0$ & $0$ & $0$ & $s$ & $e_{1}s$\tabularnewline
\hline 
$R_{e_{-2},e_{-1}}$ & $0$ & $0$ & $0$ & $0$ & $s$ & $0$ & $0$ & $e_{3}s$\tabularnewline
\hline 
\end{longtable}
\end{singlespace}

\caption{\label{tab:action-of-so7}The action of $\mathfrak{so}(7)$ on $V_{8}$}

\end{table}

\end{singlespace}

\begin{singlespace}
For $x_{2},y_{2}\in V_{2},x_{8},y_{8}\in V_{8}$, we define
\[
\left[x_{2}\otimes x_{8},y_{2}\otimes y_{8}\right]=\psi_{2}(x_{2},y_{2})p_{8}(x_{8},y_{8})+\psi_{8}(x_{8},y_{8})p_{2}(x_{2},y_{2}).
\]

\end{singlespace}

\subsection{Root system and Dynkin diagrams of $F(4)$\label{subsec:Root-system-F(4)}}

\begin{singlespace}
\noindent In this subsection, we use the structure of the root system
of $F(4)$ given in \cite[Appendix A]{Iohara2001}. Note that roots
of $F(4)$ are given by $\Phi=\Phi_{\bar{0}}\cup\Phi_{\bar{1}}$ where
\[
\Phi_{\bar{0}}=\{\pm\delta,\pm\varepsilon_{i}\pm\varepsilon_{j},\pm\varepsilon_{i}:i\neq j,i,j=1,2,3\}\text{ and }\Phi_{\bar{1}}=\{\frac{1}{2}(\pm\delta\pm\varepsilon_{1}\pm\varepsilon_{2}\pm\varepsilon_{3})\}
\]
where $\{\delta,\varepsilon_{1},\varepsilon_{2},\varepsilon_{3}\}$
is an orthogonal basis such that $(\delta,\delta)=-6$, $(\varepsilon_{i},\varepsilon_{j})=2$
if $i=j$ and $(\varepsilon_{i},\varepsilon_{j})=0$ otherwise.
\end{singlespace}

\begin{singlespace}
We list all roots and the corresponding root vectors in the table
below:
\end{singlespace}
\begin{singlespace}
\noindent \begin{center}
\begin{tabular}{|>{\centering}m{1.5cm}||>{\centering}m{0.4cm}||>{\centering}m{0.5cm}||>{\centering}m{2cm}||c||>{\centering}m{1.2cm}||>{\centering}m{1cm}||>{\centering}m{3cm}|}
\hline 
Roots & $\delta$ & $-\delta$ & $\varepsilon_{i}+\varepsilon_{j}(i<j)$ & $-\varepsilon_{i}-\varepsilon_{j}(i<j)$ & $\varepsilon_{i}-\varepsilon_{j}$ & $\pm\varepsilon_{i}$ & $\frac{1}{2}(i\delta+\sigma_{1}\varepsilon_{1}+\sigma_{2}\varepsilon_{2}+\sigma_{3}\varepsilon_{3})$\tabularnewline
\hline 
\hline 
Root vectors & $E$ & $F$ & $R_{e_{i},e_{j}}$ & $R_{e_{-j},e_{-i}}$ & $R_{e_{i},e_{-j}}$ & $R_{e_{\pm i},e_{0}}$ & $v_{i}\otimes v_{\sigma_{1},\sigma_{2},\sigma_{3}},\sigma_{i}\in\{\pm\}$\tabularnewline
\hline 
\end{tabular}
\par\end{center}
\end{singlespace}

\begin{singlespace}
The following table covers all possible Dynkin diagrams with respect
to different systems of simple roots based on \cite[Section 2.18]{Frappat1996}.
Moreover, we have that all odd roots in the root system of $F(4)$
are isotropic.
\end{singlespace}

\begin{longtable}[c]{|>{\centering}m{7cm}||>{\centering}m{7cm}|}
\caption{Dynkin diagrams for $F(4)$}
\tabularnewline
\endfirsthead
\hline 
Simple systems $\varPi=\{\alpha_{1},\alpha_{2},\alpha_{3},\alpha_{4}\}$ & Dynkin diagrams\tabularnewline
\hline 
\hline 
$\{\frac{1}{2}(\delta-\varepsilon_{1}-\varepsilon_{2}-\varepsilon_{3}),\varepsilon_{3},\varepsilon_{2}-\varepsilon_{3},\varepsilon_{1}-\varepsilon_{2}\}$ & Figure 6.1
\noindent \centering{}\includegraphics{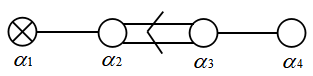}\tabularnewline
\hline 
\hline 
$\{\frac{1}{2}(-\delta+\varepsilon_{1}+\varepsilon_{2}+\varepsilon_{3}),\frac{1}{2}(\delta-\varepsilon_{1}-\varepsilon_{2}+\varepsilon_{3}),\varepsilon_{2}-\varepsilon_{3},\varepsilon_{1}-\varepsilon_{2}\}$  & Figure 6.2
\noindent \centering{}\includegraphics{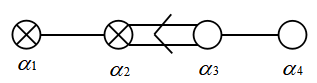}\tabularnewline
\hline 
\hline 
$\{\varepsilon_{1}-\varepsilon_{2},\frac{1}{2}(\delta-\varepsilon_{1}+\varepsilon_{2}-\varepsilon_{3}),\frac{1}{2}(-\delta+\varepsilon_{1}+\varepsilon_{2}-\varepsilon_{3}),\varepsilon_{3}\}$  & Figure 6.3
\noindent \centering{}\includegraphics[bb=0bp 0bp 130bp 86bp,scale=0.8]{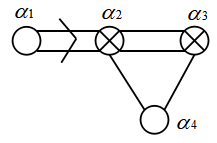}\tabularnewline
\hline 
\hline 
$\{\frac{1}{2}(\delta+\varepsilon_{1}-\varepsilon_{2}-\varepsilon_{3}),\frac{1}{2}(\delta-\varepsilon_{1}+\varepsilon_{2}+\varepsilon_{3}),\frac{1}{2}(-\delta+\varepsilon_{1}-\varepsilon_{2}+\varepsilon_{3}),\varepsilon_{2}-\varepsilon_{3}\}$ & Figure 6.4
\noindent \centering{}\includegraphics[bb=0bp 0bp 122bp 72.7899bp,scale=0.9]{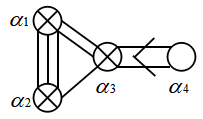}\tabularnewline
\hline 
\hline 
$\{\delta,\frac{1}{2}(-\delta+\varepsilon_{1}-\varepsilon_{2}-\varepsilon_{3}),\varepsilon_{3},\varepsilon_{2}-\varepsilon_{3}\}$ & Figure 6.5
\noindent \centering{}\includegraphics{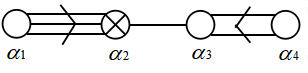}\tabularnewline
\hline 
\hline 
$\{\delta,\frac{1}{2}(-\delta-\varepsilon_{1}+\varepsilon_{2}+\varepsilon_{3}),\varepsilon_{1}-\varepsilon_{2},\varepsilon_{2}-\varepsilon_{3}\}$ & Figure 6.6
\noindent \centering{}\includegraphics{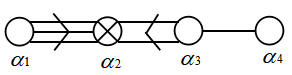}\tabularnewline
\hline 
\end{longtable}

\subsection{Centres of centralizers of nilpotent elements $e\in\mathfrak{g}_{\bar{0}}$
and corresponding labelled Dynkin diagrams\label{subsec:centre-of-Centralizer-F(4)}}

\begin{singlespace}
\noindent Let $e=e_{\mathfrak{sl}}+e_{\mathfrak{so}}\in\mathfrak{g}_{\bar{0}}$
be nilpotent where $e_{\mathfrak{sl}}\in\mathfrak{sl}(2)$ and $e_{\mathfrak{so}}\in\mathfrak{so}(7)$.
We know that there are two representatives of nilpotent orbits in
$\mathfrak{sl}(2)$, i.e. $\{0,E\}$. Based on Table \ref{tab:nilpotent ele},
there are $7$ representatives of nilpotent orbits in $\mathfrak{so}(7)$.
Hence, there are in total $14$ possibilities for $e$. We give basis
elements for $\mathfrak{g}_{\bar{1}}^{e}$ and $\mathfrak{z}(\mathfrak{g}^{e})$
and list the labelled Dynkin diagrams $\varDelta$ with respect to
$e$ in Table \ref{tab:F(4)}. Note that we have already calculated
$\mathfrak{so}(7)^{e_{\mathfrak{so}}}$ in Table \ref{tab:Centralizers-so7}
and $\mathfrak{sl}(2)^{e_{\mathfrak{sl}}}=\langle E\rangle$ for $e_{\mathfrak{sl}}=E$,
$\mathfrak{sl}(2)^{e_{\mathfrak{sl}}}=\mathfrak{sl}(2)$ for $e_{\mathfrak{sl}}=0$.
The numbers in the column labelled ``$\varDelta$'' represent labels
$a_{i}$ corresponding to $\alpha_{i}$ for $i=1,2,3,4$ in $\varDelta$.
\end{singlespace}
\begin{singlespace}
\noindent \begin{center}
\begin{longtable}[c]{|>{\centering}m{1.5cm}|>{\centering}m{6cm}|>{\centering}m{3cm}|>{\centering}m{3.5cm}|}
\caption{\label{tab:F(4)}$\mathfrak{g}^{e}$, $\mathfrak{z}(\mathfrak{g}^{e})$
and $\varDelta$ for $\mathfrak{g}=F(4)$}
\tabularnewline
\endfirsthead
\hline 
$e$ & \begin{singlespace}
\noindent $\mathfrak{g}_{\bar{1}}^{e}$ 
\end{singlespace}
 & $\mathfrak{z}(\mathfrak{g}^{e})$ & $\varDelta$\tabularnewline
\hline 
\hline 
$E+e_{(7)}$ & \begin{singlespace}
\noindent $\langle v_{1}\otimes e_{1}e_{2}e_{3}s,v_{1}\otimes e_{1}e_{2}s-v_{-1}\otimes e_{1}e_{2}e_{3}s,v_{1}\otimes e_{1}s-v_{1}\otimes e_{2}e_{3}s\rangle$
\end{singlespace}
 & \begin{singlespace}
\noindent $\langle e,v_{1}\otimes e_{1}e_{2}e_{3},R_{e_{1},e_{2}}\rangle$
\end{singlespace}
 & Figure 6.4: $1,1,1,2$\tabularnewline
\hline 
\hline 
$E+e_{(5,1^{2})}$ & \begin{singlespace}
\noindent $\langle v_{1}\otimes e_{1}e_{2}e_{3}s,v_{1}\otimes e_{1}e_{2}s,v_{1}\otimes e_{1}s-v_{-1}\otimes e_{1}e_{2}s,v_{1}\otimes e_{1}e_{3}s+v_{-1}\otimes e_{1}e_{2}e_{3}s\rangle$
\end{singlespace}
 & $\langle e,R_{e_{1},e_{2}}\rangle$ & Figure 6.3: $2,0,2,0$

Figure 6.4: $2,0,0,2$

Figure 6.5: $2,0,0,2$\tabularnewline
\hline 
\hline 
$E+e_{(3^{2},1)}$ & $\langle v_{1}\otimes e_{1}e_{2}s,v_{1}\otimes e_{1}e_{2}e_{3}s,v_{1}\otimes e_{1}s-v_{-1}\otimes e_{1}e_{2}e_{3}s,v_{1}\otimes e_{2}s,v_{-1}\otimes e_{1}e_{2}s+v_{1}\otimes e_{2}e_{3}s,v_{1}\otimes e_{1}e_{3}s\rangle$ & $\langle e,R_{e_{1},e_{2}}\rangle$ & Figure 6.3: $0,1,1,0$\tabularnewline
\hline 
\hline 
$E+e_{(3,2^{2})}$ & \begin{singlespace}
\noindent $\langle v_{1}\otimes e_{1}e_{2}e_{3}s,v_{1}\otimes e_{1}e_{2}s,v_{1}\otimes e_{1}e_{3}s,v_{1}\otimes e_{1}s-v_{-1}\otimes e_{1}e_{2}e_{3}s,v_{1}\otimes e_{2}e_{3}s-v_{-1}\otimes e_{1}e_{2}e_{3}s,v_{1}\otimes e_{3}s+v_{-1}\otimes e_{1}e_{3}s,v_{1}\otimes e_{2}s+v_{-1}\otimes e_{1}e_{2}s\rangle$
\end{singlespace}
 & $\langle e\rangle$ & Figure 6.2: $1,0,0,1$

Figure 6.3: $1,0,0,1$

Figure 6.4: $1,1,0,0$\tabularnewline
\hline 
\hline 
$E+e_{(3,1^{4})}$ & \begin{singlespace}
\noindent $\langle v_{1}\otimes s-v_{-1}\otimes e_{1}s,v_{1}\otimes e_{2}s+v_{-1}\otimes e_{1}e_{2}s,$
$v_{1}\otimes e_{3}s+v_{-1}\otimes e_{1}e_{3}s,v_{1}\otimes e_{1}s,v_{1}\otimes e_{1}e_{3}s,v_{1}\otimes e_{1}e_{2}e_{3}s$
$v_{1}\otimes e_{2}e_{3}s-v_{-1}\otimes e_{1}e_{2}e_{3}s,v_{1}\otimes e_{1}e_{2}s\rangle$
\end{singlespace}
 & $\langle e\rangle$ & Figure 6.1: $0,0,0,2$

Figure 6.2: $0,0,0,2$

Figure 6.3: $2,0,0,0$

Figure 6.4: $2,0,0,0$\tabularnewline
\hline 
\hline 
$E+e_{(2^{2},1^{3})}$ & \begin{singlespace}
\noindent $\langle v_{1}\otimes e_{1}e_{2}e_{3}s,v_{1}\otimes e_{1}e_{2}s,$
$v_{1}\otimes s-v_{-1}\otimes e_{1}e_{2}s,v_{1}\otimes e_{3}s-v_{-1}\otimes e_{1}e_{2}e_{3}s,$
$v_{1}\otimes e_{1}s,v_{1}\otimes e_{2}s,v_{1}\otimes e_{1}e_{3}s,v_{1}\otimes e_{2}e_{3}s\rangle$
\end{singlespace}
 & $\langle e\rangle$ & Figure 6.1: $0,0,1,0$

Figure 6.2: $0,0,1,0$

Figure 6.3: $0,1,0,0$\tabularnewline
\hline 
\hline 
$E$ & \begin{singlespace}
\noindent $\langle v_{1}\otimes s,v_{1}\otimes e_{1}s,v_{1}\otimes e_{2}s,v_{1}\otimes e_{3}s,v_{1}\otimes e_{1}e_{2}s,v_{1}\otimes e_{1}e_{3}s,v_{1}\otimes e_{2}e_{3}s,v_{1}\otimes e_{1}e_{2}e_{3}s\rangle$
\end{singlespace}
 & $\langle e\rangle$ & Figure 6.1: $1,0,0,0$\tabularnewline
\hline 
\hline 
$e_{(7)}$ & \begin{singlespace}
\noindent $\langle v_{1}\otimes e_{1}s-v_{1}\otimes e_{2}e_{3}s,v_{-1}\otimes e_{1}s-v_{-1}\otimes e_{2}e_{3}s,v_{1}\otimes e_{1}e_{2}e_{3}s,v_{-1}\otimes e_{1}e_{2}e_{3}s\rangle$
\end{singlespace}
 & $\langle e,R_{e_{1},e_{2}}\rangle$ & Figure 6.4: $0,0,2,2$

Figure 6.5: $0,0,2,2$

Figure 6.6: $0,0,2,2$\tabularnewline
\hline 
\hline 
$e_{(5,1^{2})}$ & \begin{singlespace}
\noindent \centering{}$\langle v_{1}\otimes e_{1}e_{2}s,v_{-1}\otimes e_{1}e_{2}s,v_{1}\otimes e_{1}e_{2}e_{3}s,v_{-1}\otimes e_{1}e_{2}e_{3}s\rangle$
\end{singlespace}
 & $\langle e,R_{e_{1},e_{2}}\rangle$ & Figure 6.5: $0,1,0,2$\tabularnewline
\hline 
\hline 
$e_{(3^{2},1)}$ & \begin{singlespace}
\noindent \centering{}$\langle v_{1}\otimes e_{2}s,v_{-1}\otimes e_{2}s,v_{1}\otimes e_{1}e_{3}s,v_{-1}\otimes e_{1}e_{3}s,v_{1}\otimes e_{1}e_{2}s,v_{-1}\otimes e_{1}e_{2}s,v_{1}\otimes e_{1}e_{2}e_{3}s,v_{-1}\otimes e_{1}e_{2}e_{3}s\rangle$
\end{singlespace}
 & $\langle e,R_{e_{1},e_{2}}\rangle$ & Figure 6.3: $0,0,2,0$

Figure 6.4: $0,0,0,2$

Figure 6.5: $0,0,0,2$

Figure 6.6: $0,0,0,2$\tabularnewline
\hline 
\hline 
$e_{(3,2^{2})}$ & \begin{singlespace}
\noindent \centering{}$\langle v_{1}\otimes e_{1}s-v_{1}\otimes e_{2}e_{3}s,v_{-1}\otimes e_{1}s-v_{-1}\otimes e_{2}e_{3}s,v_{1}\otimes e_{1}e_{2}s,v_{-1}\otimes e_{1}e_{2}s,v_{1}\otimes e_{1}e_{3}s,v_{-1}\otimes e_{1}e_{3}s,v_{1}\otimes e_{1}e_{2}e_{3}s,v_{-1}\otimes e_{1}e_{2}e_{3}s\rangle$
\end{singlespace}
 & $\langle e\rangle$ & Figure 6.4: $0,0,1,0$

Figure 6.5: $0,0,1,0$

Figure 6.6: $0,0,1,0$\tabularnewline
\hline 
\hline 
$e_{(3,1^{4})}$ & $\langle v_{1}\otimes e_{1}e_{2}e_{3}s,$ $v_{-1}\otimes e_{1}e_{2}e_{3}s,v_{1}\otimes e_{1}s,v_{-1}\otimes e_{1}s,v_{1}\otimes e_{1}e_{2}s,$
$v_{-1}\otimes e_{1}e_{2}s,v_{1}\otimes e_{1}e_{3}s,v_{-1}\otimes e_{1}e_{3}s\rangle$ & $\langle e\rangle$ & Figure 6.5: $0,1,0,0$\tabularnewline
\hline 
\hline 
$e_{(2^{2},1^{3})}$ & $\langle v_{1}\otimes e_{1}s,v_{-1}\otimes e_{1}s,v_{-1}\otimes e_{1}e_{2}e_{3}s,$
$v_{1}\otimes e_{1}e_{2}e_{3}s,v_{1}\otimes e_{1}e_{2}s,v_{-1}\otimes e_{1}e_{2}s,v_{1}\otimes e_{2}s,$
$v_{-1}\otimes e_{2}s,v_{1}\otimes e_{1}e_{3}s,v_{-1}\otimes e_{1}e_{3}s,v_{1}\otimes e_{2}e_{3}s,v_{-1}\otimes e_{2}e_{3}s\rangle$ & $\langle e\rangle$ & Figure 6.3: $0,0,1,0$

Figure 6.4: $0,0,0,1$

Figure 6.5: $0,0,0,1$

Figure 6.6: $0,0,0,1$\tabularnewline
\hline 
\hline 
$0$ & \begin{singlespace}
\noindent $\mathfrak{g}_{\bar{1}}$
\end{singlespace}
 & $\{0\}$ & Figures 6.1, 6.2, 6.3, 6.4, 6.5, 6.6: All labels are zeros\tabularnewline
\hline 
\end{longtable}
\par\end{center}
\end{singlespace}

\begin{singlespace}
We also calculate the $\mathfrak{g}^{e}(0)$-module structure on each
$\mathfrak{g}^{e}(j)$ for $j>0$. Recall that we denote by $V^{\mathfrak{sl}}(j)$
the $\mathfrak{sl}(2)$-module with highest weight $j$. We also let
$V^{\mathfrak{osp}}(j)$ be the $\mathfrak{osp}(1|2)$-module with
highest weight $j$. Let $\mathfrak{t}$ be a 1-dimensional Lie algebra
and $\mathfrak{t}_{j}$ be the $\mathfrak{t}$-module such that $t\cdot a=ja$
for $t\in\mathfrak{t}$, $a\in\mathfrak{t}_{j}$. For $e=e_{(3^{2},1)}$
and $e_{(2^{2},1^{3})}$, the $\mathfrak{g}^{e}(0)$-module structure
on $\mathfrak{g}^{e}(j)$ is not included as it requires representations
of $\mathfrak{sl}(2|1)$ and $D(2,1;\alpha)$.
\end{singlespace}
\begin{singlespace}
\noindent \begin{center}
\begin{longtable}[c]{|>{\centering}m{2cm}|>{\centering}m{3.5cm}|>{\centering}m{8cm}|}
\caption{\label{tab:g^e(0)-F(4)}The $\mathfrak{g}^{e}(0)$-module structure
on $\mathfrak{g}^{e}(j)$ for $j>0$ }
\tabularnewline
\endfirsthead
\hline 
$e$ & $\mathfrak{g}^{e}(0)$ & $\mathfrak{g}^{e}(j)$ for $j>0$\tabularnewline
\hline 
\hline 
$E+e_{(7)}$ & $0$ & $\dim\mathfrak{g}^{e}(10)=\dim\mathfrak{g}^{e}(7)=\dim\mathfrak{g}^{e}(6)=\dim\mathfrak{g}^{e}(5)=\dim\mathfrak{g}^{e}(1)=1,\dim\mathfrak{g}^{e}(2)=2.$\tabularnewline
\hline 
\hline 
$E+e_{(5,1^{2})}$ & $\mathfrak{t}$ & $\mathfrak{g}^{e}(1)=0,\mathfrak{g}^{e}(2)=\mathfrak{t}_{0}\oplus\mathfrak{t}_{-1}\oplus\mathfrak{t}_{1},$$\mathfrak{g}^{e}(4)=\mathfrak{t}_{-2}\oplus\mathfrak{t}_{-1}\oplus\mathfrak{t}_{1}\oplus\mathfrak{t}_{2},\mathfrak{g}^{e}(6)=\mathfrak{t}_{0}.$\tabularnewline
\hline 
\hline 
$E+e_{(3^{2},1)}$ & $\mathfrak{t}$ & $\mathfrak{g}^{e}(1)=\mathfrak{t}_{-3}\oplus\mathfrak{t}_{-1}\oplus\mathfrak{t}_{1}\oplus\mathfrak{t}_{3}$;$\mathfrak{g}^{e}(2)=\mathfrak{t}_{-4}\oplus\mathfrak{t}_{-2}\oplus\mathfrak{t}_{0}\oplus\mathfrak{t}_{2}\oplus\mathfrak{t}_{-4};$$\mathfrak{g}^{e}(3)=\mathfrak{t}_{-1}\oplus\mathfrak{t}_{1},\mathfrak{g}^{e}(4)=0.$\tabularnewline
\hline 
\hline 
$E+e_{(3,2^{2})}$ & $\mathfrak{osp}(1|2)$ & $\mathfrak{g}^{e}(1)=V^{\mathfrak{osp}}(1)\oplus V^{\mathfrak{osp}}(0),$
$\mathfrak{g}^{e}(2)=V^{\mathfrak{osp}}(1)\oplus V^{\mathfrak{osp}}(0)\oplus V^{\mathfrak{osp}}(0),$
$\mathfrak{g}^{e}(3)=V^{\mathfrak{osp}}(1).$\tabularnewline
\hline 
\hline 
$E+e_{(3,1^{4})}$ & $\mathfrak{osp}(1|2)\oplus\mathfrak{osp}(1|2)$ & $\mathfrak{g}^{e}(2)=\left(V^{\mathfrak{osp}}(1)\otimes V^{\mathfrak{osp}}(1)\right)\oplus\left(V^{\mathfrak{osp}}(0)\oplus V^{\mathfrak{osp}}(0)\right)$\tabularnewline
\hline 
\hline 
$E+e_{(2^{2},1^{3})}$ & $\mathfrak{sl}(2)\oplus\mathfrak{osp}(1|2)$ & $\mathfrak{g}^{e}(1)=V^{\mathfrak{sl}}(1)\otimes V^{\mathfrak{osp}}(2),$$\mathfrak{g}^{e}(2)=\left(V^{\mathfrak{sl}}(0)\otimes V^{\mathfrak{osp}}(0)\right)\oplus\left(V^{\mathfrak{sl}}(0)\otimes V^{\mathfrak{osp}}(1)\right).$\tabularnewline
\hline 
\hline 
$E$ & $\mathfrak{so}(7)$ & $\mathfrak{g}^{e}(1)=V_{8}$, $\mathfrak{g}^{e}(2)=\left\langle e\right\rangle .$\tabularnewline
\hline 
\hline 
$e_{(7)}$ & $\mathfrak{osp}(1|2)$ & $\mathfrak{g}^{e}(2)=\mathfrak{g}^{e}(10)=V^{\mathfrak{osp}}(0),\mathfrak{g}^{e}(6)=V^{\mathfrak{osp}}(1)$.\tabularnewline
\hline 
\hline 
$e_{(5,1^{2})}$ & $\mathfrak{sl}(2)\oplus\mathfrak{t}$ & $\mathfrak{g}^{e}(2)=V^{\mathfrak{sl}}(0)\otimes\mathfrak{t}_{0},$
$\mathfrak{g}^{e}(3)=\left(V^{\mathfrak{sl}}(1)\otimes\mathfrak{t}_{-1}\right)\oplus\left(V^{\mathfrak{sl}}(1)\otimes\mathfrak{t}_{1}\right),$$\mathfrak{g}^{e}(4)=\left(V^{\mathfrak{sl}}(0)\otimes\mathfrak{t}_{-2}\right)\oplus\left(V^{\mathfrak{sl}}(0)\otimes\mathfrak{t}_{2}\right),$
$\mathfrak{g}^{e}(6)=V^{\mathfrak{sl}}(0)\otimes\mathfrak{t}_{0}.$\tabularnewline
\hline 
\hline 
$e_{(3^{2},1)}$ & $\mathfrak{sl}(2|1)$ & Omitted\tabularnewline
\hline 
\hline 
$e_{(3,2^{2})}$ & $\mathfrak{sl}(2)\oplus\mathfrak{osp}(1\vert2)$ & $\mathfrak{g}^{e}(1)=V^{\mathfrak{sl}}(1)\otimes V^{\mathfrak{osp}}(1),\mathfrak{g}^{e}(3)=V^{\mathfrak{sl}}(1)\otimes V^{\mathfrak{osp}}(0),$
$\mathfrak{g}^{e}(2)=\left(V^{\mathfrak{sl}}(0)\otimes V^{\mathfrak{osp}}(0)\right)\oplus\left(V^{\mathfrak{sl}}(0)\otimes V^{\mathfrak{osp}}(1)\right).$\tabularnewline
\hline 
\hline 
$e_{(3,1^{4})}$ & $\mathfrak{sl}(2)\oplus\mathfrak{sl}(2)\oplus\mathfrak{sl}(2)$ & $\mathfrak{g}^{e}(1)=\left(V^{\mathfrak{sl}}(1)\otimes V^{\mathfrak{sl}}(0)\otimes V^{\mathfrak{sl}}(1)\right)\oplus\left(V^{\mathfrak{sl}}(1)\otimes V^{\mathfrak{sl}}(1)\otimes V^{\mathfrak{sl}}(0)\right),$$\mathfrak{g}^{e}(2)=\left(V^{\mathfrak{sl}}(0)\otimes V^{\mathfrak{sl}}(1)\otimes V^{\mathfrak{sl}}(1)\right)\oplus\left(V^{\mathfrak{sl}}(0)\otimes V^{\mathfrak{sl}}(0)\otimes V^{\mathfrak{sl}}(0)\right).$\tabularnewline
\hline 
\hline 
$e_{(2^{2},1^{3})}$ & $D(2,1;2)$ & Omitted\tabularnewline
\hline 
\hline 
$0$ & $\mathfrak{g}$ & $0$\tabularnewline
\hline 
\end{longtable}
\par\end{center}
\end{singlespace}

\begin{singlespace}
In the remaining part of this subsection, we give explicit calculations
for finding $\mathfrak{g}^{e}$ and $\mathfrak{z}(\mathfrak{g}^{e})$
and obtain the corresponding labelled Dynkin diagrams for $e=e_{(7)}$.
The results for all other cases are obtained using the same approach. 

For $e=e_{(7)}$, a basis for $\mathfrak{so}(7)^{e_{(7)}}$ is given
in Table \ref{tab:Centralizers-so7}. Hence, $\mathfrak{g}_{\bar{0}}^{e}=\mathfrak{sl}(2)\oplus\mathfrak{so}(7)^{e_{(7)}}$
and it has dimension $3+3=6$. 

Next we determine $\mathfrak{g}_{\bar{1}}^{e}$. By letting $[e_{(7)},a_{1}v_{1}\otimes s+a_{2}v_{1}\otimes e_{1}s+a_{3}v_{1}\otimes e_{2}s+a_{4}v_{1}\otimes e_{3}s+a_{5}v_{1}\otimes e_{1}e_{2}s+a_{6}v_{1}\otimes e_{1}e_{3}s+a_{7}v_{1}\otimes e_{2}e_{3}s+a_{8}v_{1}\otimes e_{1}e_{2}e_{3}s+b_{1}v_{-1}\otimes s+b_{2}v_{-1}\otimes e_{1}s+b_{3}v_{-1}\otimes e_{2}s+b_{4}v_{-1}\otimes e_{3}s+b_{5}v_{-1}\otimes e_{1}e_{2}s+b_{6}v_{-1}\otimes e_{1}e_{3}s+b_{7}v_{-1}\otimes e_{2}e_{3}s+b_{8}v_{1-}\otimes e_{1}e_{2}e_{3}s]=0$,
we have that $a_{2}+a_{7}=0$, $b_{2}+b_{7}=0$ and $a_{i}=b_{i}=0$
for $i=1,3,4,5,6$. Therefore, we obtain that $\mathfrak{g}_{\bar{1}}^{e}$
has a basis $\{v_{1}\otimes e_{1}s-v_{1}\otimes e_{2}e_{3}s,v_{-1}\otimes e_{1}s-v_{-1}\otimes e_{2}e_{3}s,v_{1}\otimes e_{1}e_{2}e_{3}s,v_{-1}\otimes e_{1}e_{2}e_{3}s\}$.
\end{singlespace}

\begin{singlespace}
\noindent According to Table \ref{tab:nilpotent ele}, there is an
$\mathfrak{sl}(2)$-triple $\{e,h,f\}$ in $\mathfrak{g}_{\bar{0}}$
such that $h=h_{(7)}=\mathrm{diag}(6,4,2,0,-2,-4,-6)=6R_{e_{1},e_{-1}}+4R_{e_{2},e_{-2}}+2R_{e_{3},e_{-3}}$.
Then the $\mathrm{ad}h$-eigenvalues of basis elements in $\mathfrak{g}^{e}$
are shown in the following table:
\end{singlespace}
\begin{doublespace}
\noindent \begin{center}
\begin{tabular}{|c|>{\centering}m{10cm}|}
\hline 
$\mathrm{ad}h$-eigenvalues & Basis elements in $\mathfrak{g}^{e}$ \tabularnewline
\hline 
\hline 
$0$ & $E$, $H$, $F$, $v_{1}\otimes e_{1}s-v_{1}\otimes e_{2}e_{3}s$,
$v_{-1}\otimes e_{1}s-v_{-1}\otimes e_{2}e_{3}s$\tabularnewline
\hline 
$2$ & $e_{(7)}$\tabularnewline
\hline 
$6$ & $R_{e_{1},e_{0}}-2R_{e_{2},e_{3}}$, $v_{1}\otimes e_{1}e_{2}e_{3}s$,
$v_{-1}\otimes e_{1}e_{2}e_{3}s$\tabularnewline
\hline 
$10$ & $R_{e_{1},e_{2}}$\tabularnewline
\hline 
\end{tabular}
\par\end{center}
\end{doublespace}

\begin{singlespace}
By computing commutators in $\mathfrak{g}^{e}(0)$, we deduce that
$\mathfrak{g}^{e}(0)=\mathfrak{osp}(1|2)$ by Lemma \ref{lem:osp(1,2)}
where $F,v_{-1}\otimes e_{1}s-v_{-1}\otimes e_{2}e_{3}s,H,v_{1}\otimes e_{1}s-v_{1}\otimes e_{2}e_{3}s,E$
correspond to $u_{-2},u_{-1},u_{0},u_{1},u_{2}$ in Lemma \ref{lem:osp(1,2)}
respectively. Moreover, both $\mathfrak{g}^{e}(2)$ and $\mathfrak{g}^{e}(10)$
are isomorphic to $V^{\mathfrak{osp}}(0)$ and $\mathfrak{g}^{e}(6)=V^{\mathfrak{osp}}(1)$.
Hence, we deduce that 
\begin{align*}
\mathfrak{z} & =\mathfrak{z}(0)\oplus\mathfrak{z}(2)\oplus\mathfrak{z}(6)\oplus\mathfrak{z}(10)\\
 & \subseteq\left(\mathfrak{g}^{e}(0)\right)^{\mathfrak{g}^{e}(0)}\oplus\left(\mathfrak{g}^{e}(2)\right)^{\mathfrak{g}^{e}(0)}\oplus\left(\mathfrak{g}^{e}(6)\right)^{\mathfrak{g}^{e}(0)}\oplus\left(\mathfrak{g}^{e}(10)\right)^{\mathfrak{g}^{e}(0)}=\langle e_{(7)},R_{e_{1},e_{2}}\rangle.
\end{align*}

\end{singlespace}

\noindent We check that $e_{(7)},R_{e_{1},e_{2}}\in\mathfrak{z}$,
therefore $\mathfrak{z}=\langle e_{(7)},R_{e_{1},e_{2}}\rangle$.

\begin{singlespace}
Next we look at the labelled Dynkin diagrams with respect to $e$.
We obtain that roots in $\mathfrak{g}(>0)$ are $\{\varepsilon_{1}+\varepsilon_{2},\varepsilon_{1}+\varepsilon_{3},\varepsilon_{2}+\varepsilon_{3},\varepsilon_{1}-\varepsilon_{3},\varepsilon_{1}-\varepsilon_{2},\varepsilon_{2}-\varepsilon_{3},\varepsilon_{1},\varepsilon_{2},\varepsilon_{3},\frac{1}{2}(\pm\delta+\varepsilon_{1}+\varepsilon_{2}-\varepsilon_{3}),\frac{1}{2}(\pm\delta+\varepsilon_{1}-\varepsilon_{2}+\varepsilon_{3}),\frac{1}{2}(\pm\delta+\varepsilon_{1}+\varepsilon_{2}+\varepsilon_{3})\}$
and roots in $\mathfrak{g}(0)$ are $\Phi(0)=\{\pm\delta,\frac{1}{2}(\pm\delta+\varepsilon_{1}-\varepsilon_{2}-\varepsilon_{3}),\frac{1}{2}(\pm\delta-\varepsilon_{1}+\varepsilon_{2}+\varepsilon_{3})\}$.
Hence, there are three systems of simple roots of $\mathfrak{g}(0)$
up to conjugacy: $\varPi_{1}(0)=\{\frac{1}{2}(\delta+\varepsilon_{1}-\varepsilon_{2}-\varepsilon_{3}),\frac{1}{2}(\delta-\varepsilon_{1}+\varepsilon_{2}+\varepsilon_{3})\}$,
$\varPi_{2}(0)=\{\delta,\frac{1}{2}(-\delta+\varepsilon_{1}-\varepsilon_{2}-\varepsilon_{3})\}$
and $\varPi_{3}(0)=\{\delta,\frac{1}{2}(-\delta-\varepsilon_{1}+\varepsilon_{2}+\varepsilon_{3})\}$.
By extending $\varPi_{i}(0)$ to simple root systems of $\mathfrak{g}$,
we get three systems of positive roots $\Phi^{+}$ and simple roots
$\varPi$ and thus there exist three conjugacy classes of Borel subalgebras
that satisfy $\mathfrak{b}=\mathfrak{h}\oplus\bigoplus_{\alpha\in\Phi^{+}}\mathfrak{g}_{\alpha}\subseteq\bigoplus_{j\geq0}\mathfrak{g}(j)$.
Hence, the systems of simple roots are:

$\varPi_{1}=\{\alpha_{1}=\frac{1}{2}(\delta+\varepsilon_{1}-\varepsilon_{2}-\varepsilon_{3}),\alpha_{2}=\frac{1}{2}(\delta-\varepsilon_{1}+\varepsilon_{2}+\varepsilon_{3}),\alpha_{3}=\frac{1}{2}(-\delta+\varepsilon_{1}-\varepsilon_{2}+\varepsilon_{3}),\alpha_{4}=\varepsilon_{2}-\varepsilon_{3}\}$.
We compute $\mu_{12}=3$, $\mu_{13}=2$, $\mu_{23}=1$ and $\mu_{34}=2$
using Formula (\ref{eq:lines-=0003BC}). Therefore, the labelled Dynkin
diagram with repect to $\varPi_{1}$ is the Dynkin diagram in Figure
6.4 with labels $0,0,2,2$.

$\varPi_{2}=\{\alpha_{1}=\delta,\alpha_{2}=\frac{1}{2}(-\delta+\varepsilon_{1}-\varepsilon_{2}-\varepsilon_{3}),\alpha_{3}=\varepsilon_{3},\alpha_{4}=\varepsilon_{2}-\varepsilon_{3}\}$.
We compute $\mu_{12}=3$, $\mu_{23}=1$ and $\mu_{34}=2$ using Formula
(\ref{eq:lines-=0003BC}). Therefore, the labelled Dynkin diagram
with repect to $\varPi_{2}$ is the Dynkin diagram in Figure 6.5 with
labels $0,0,2,2$.

$\varPi_{3}=\{\alpha_{1}=\delta,\alpha_{2}=\frac{1}{2}(-\delta-\varepsilon_{1}+\varepsilon_{2}+\varepsilon_{3}),\alpha_{3}=\varepsilon_{1}-\varepsilon_{2},\alpha_{4}=\varepsilon_{2}-\varepsilon_{3}\}$.
We compute $\mu_{12}=3$, $\mu_{23}=2$ and $\mu_{34}=1$ using Formula
(\ref{eq:lines-=0003BC}). Therefore, the labelled Dynkin diagram
with repect to $\varPi_{3}$ is the Dynkin diagram in Figure 6.6 with
labels $0,0,2,2$.
\end{singlespace}

\subsection{Analysis of results}

\begin{singlespace}
\noindent For each nilpotent element $e\in\mathfrak{g}_{\bar{0}}$,
a semisimple element $h\in\mathfrak{g}_{\bar{0}}$ is given in Table
\ref{tab:nilpotent ele}. Denote a root system for $\mathfrak{g}^{h}$
by $\Phi_{h}$, i.e. $\Phi_{h}=\{\alpha\in\Phi:\alpha(h)=0\}$ and
a simple root system for $\mathfrak{g}^{h}$ by $\varPi_{h}$. We
also denote the labelled Dynkin diagram for $\varPi_{h}$ by $\varDelta_{h}$
.
\end{singlespace}

\begin{singlespace}
In order to prove Theorem 1 for $\mathfrak{g}$, we need to determine
$\mathfrak{g}^{h}$ for each case such that $\varDelta$ has no label
equal to $1$. Note that $\mathfrak{g}^{h}$ is of the form $\mathfrak{s}\oplus\bigoplus_{\alpha\in\Phi_{h}}\mathfrak{g}_{\alpha}$
where $\mathfrak{s}$ is a subalgebra of the Cartan subalgebra $\mathfrak{h}$
of $\mathfrak{g}$. Then $\mathfrak{z}(\mathfrak{g}^{h})=\{t\in\mathfrak{h}:\alpha(t)=0\text{ for all }\alpha\in\Phi_{h}\}$,
thus $\mathfrak{z}(\mathfrak{g}^{h})$ is a subalgebra of $\mathfrak{h}$
with dimension $\mathrm{rank}\Phi-\mathrm{rank}\Phi_{h}$. Note that
$\varDelta$ has no label equal to $1$ when the nilpotent elements
are $E+e_{(5,1^{2})}$, $E+e_{(3,1^{4})}$, $e_{(7)},e_{(3^{2},1)}$.
We take $E+e_{(5,1^{2})}$ and $e_{(3^{2},1)}$ as examples to show
explicit calculation on $\dim\mathfrak{z}(\mathfrak{g}^{h})$. The
results for all other cases are obtained by the same method.

When $e=E+e_{(5,1^{2})}$, we have that $\Phi_{h}=\{\pm\varepsilon_{3},\pm\frac{1}{2}(\delta-\varepsilon_{1}+\varepsilon_{2}-\varepsilon_{3}),\pm\frac{1}{2}(\delta-\varepsilon_{1}+\varepsilon_{2}+\varepsilon_{3})\}$.
Then up to conjugacy the simple root systems are $\varPi_{h}^{1}=\{\varepsilon_{3},\frac{1}{2}(\delta-\varepsilon_{1}+\varepsilon_{2}-\varepsilon_{3})\}$
and $\varPi_{h}^{2}=\{\frac{1}{2}(\delta-\varepsilon_{1}+\varepsilon_{2}+\varepsilon_{3}),-\frac{1}{2}(\delta+\varepsilon_{1}-\varepsilon_{2}+\varepsilon_{3})\}$.
Hence $\Phi_{h}$ is of type $\mathfrak{sl}(2|1)$  and $\mathfrak{g}^{h}=\mathfrak{s}\oplus\mathfrak{sl}(2|1)$
where $\mathfrak{s}$ is a complement of $\mathfrak{h}\cap\mathfrak{sl}(2|1)$
in $\mathfrak{h}$. Note that $\mathfrak{sl}(2|1)$ has no centre,
thus $\dim\mathfrak{z}(\mathfrak{g}^{h})=4-2=2=n_{2}(\varDelta)=\dim\mathfrak{z}(\mathfrak{g}^{e})$.

When $e=e_{(3^{2},1)}$, we have that $\Phi_{h}=\{\pm\delta,\pm(\varepsilon_{1}-\varepsilon_{2}),\pm\varepsilon_{3},\pm\frac{1}{2}(\delta+\varepsilon_{1}-\varepsilon_{2}-\varepsilon_{3}),\pm\frac{1}{2}(\delta-\varepsilon_{1}+\varepsilon_{2}-\varepsilon_{3}),\pm\frac{1}{2}(\delta+\varepsilon_{1}-\varepsilon_{2}+\varepsilon_{3}),\pm\frac{1}{2}(\delta-\varepsilon_{1}+\varepsilon_{2}+\varepsilon_{3})\}$.
Then up to conjugacy the simple root systems are $\varPi_{h}^{1}=\{\varepsilon_{1}-\varepsilon_{2},\frac{1}{2}(\delta-\varepsilon_{1}+\varepsilon_{2}-\varepsilon_{3}),\varepsilon_{3}\}$,
$\varPi_{h}^{2}=\{\frac{1}{2}(\delta+\varepsilon_{1}-\varepsilon_{2}-\varepsilon_{3}),\frac{1}{2}(\delta-\varepsilon_{1}+\varepsilon_{2}+\varepsilon_{3}),\frac{1}{2}(-\delta+\varepsilon_{1}-\varepsilon_{2}+\varepsilon_{3})\}$,
$\varPi_{h}^{3}=\{\delta,\frac{1}{2}(-\delta+\varepsilon_{1}-\varepsilon_{2}-\varepsilon_{3}),\varepsilon_{3}\}$
and $\varPi_{h}^{4}=\{\delta,\frac{1}{2}(-\delta-\varepsilon_{1}+\varepsilon_{2}+\varepsilon_{3}),\varepsilon_{1}-\varepsilon_{2}\}$.
Hence $\Phi_{h}$ is of type $D(2,1;2)$ according to \S4.3 and $\mathfrak{g}^{h}=\mathfrak{s}\oplus D(2,1;2)$
where $\mathfrak{s}$ is a complement of $\mathfrak{h}\cap D(2,1;2)$
in $\mathfrak{h}$. Note that $D(2,1;2)$ has no centre, thus $\dim\mathfrak{z}(\mathfrak{g}^{h})=4-3=1=n_{2}(\varDelta)$
but $\dim\mathfrak{z}(\mathfrak{g}^{h})\neq\dim\mathfrak{z}(\mathfrak{g}^{e})$.
We will further discuss this case in \S6.7.

In order to prove Theorem 2 for $\mathfrak{g}$, we only need to look
at cases such that $\varDelta$ has labels equal to $2$ as for the
remaining cases $\mathfrak{g}_{0}=\mathfrak{g}$, $e_{0}=e$ and $n_{2}(\varDelta)=0$.
These cases are $E+e_{(7)},E+e_{(5,1^{2})},E+e_{(3,1^{4})},e_{(7)},e_{(5,1^{2})},e_{(3^{2},1)}.$
We take $E+e_{(7)}$ and $e_{(3^{2},1)}$ as examples to show explicit
analysis. The results for all other cases are obtained by the same
method.

When $e=E+e_{(7)}$, we have that $\mathfrak{g}_{0}=D(2,1;2)$ and
$e_{0}=(E,E,E)$ according to \S4.3. Therefore, we obtain that $\dim\mathfrak{g}^{e}-\dim\mathfrak{g}_{0}^{e_{0}}=7-6=1=n_{2}(\varDelta)$
but $\dim\mathfrak{z}(\mathfrak{g}^{e})-\dim\mathfrak{z}(\mathfrak{g}_{0}^{e_{0}})=3-1\neq n_{2}(\varDelta)$.
We will further discuss this case in \S6.7.

When $e=e_{(3^{2},1)}$, we have that $\mathfrak{g}_{0}=D(2,1;2)$
and $e_{0}=0$ by looking at $\varDelta_{0}$. Note that $\dim D(2,1;2)=17$.
Therefore, we obtain that $\dim\mathfrak{g}^{e}-\dim\mathfrak{g}_{0}^{e_{0}}=18-17=1=n_{2}(\varDelta)$
but $\dim\mathfrak{z}(\mathfrak{g}^{e})-\dim\mathfrak{z}(\mathfrak{g}_{0}^{e_{0}})=2\neq n_{2}(\varDelta)$.
We will further discuss this case in \S6.7.
\end{singlespace}

\subsection{Adjoint action on $F(4)$\label{subsec:Adjoint-action-of-F(4)}}

\begin{singlespace}
\noindent Let $G$ be a linear algebraic group $G=\mathrm{SL}_{2}(\mathbb{C})\times\mathrm{Spin}_{7}(\mathbb{C})$.
In this subsection, we determine $\left(\mathfrak{z}(\mathfrak{g}^{e})\right)^{G^{e}}$
in order to complete the proof of the theorems. Note that we only
need to look at cases $e=E+e_{(7)}$, $E+e_{(5,1^{2})}$, $E+e_{(3^{2},1)}$,
$e_{(7)}$, $e_{(5,1^{2})}$, $e_{(3^{2},1)}$ as for all other cases
we have $\left(\mathfrak{z}(\mathfrak{g}^{e})\right)^{G^{e}}=\mathfrak{z}(\mathfrak{g}^{e})=\langle e\rangle$.
\end{singlespace}

\begin{singlespace}
For $e=e_{(7)},e_{(5,1^{2})},e_{(3^{2},1)}$, the results can be obtained
via \cite[Proposition 4.2]{Lawther2008} since all calculation take
place in $\mathfrak{so}(7)$. In the remaining part of this subsection,
we include details for the case $e=e_{(3^{2},1)}$ as an example.
When $e=e_{(3^{2},1)}$, recall that $\mathfrak{z}(\mathfrak{g}^{e})=\langle e_{(3^{2},1)},R_{e_{1},e_{2}}\rangle$.
According to \cite[Theorem 6.3.5]{Goodman2009}, there exists a homomorphism
$\mathrm{id}\times\pi:G\rightarrow\mathrm{SL}_{2}(\mathbb{C})\times\mathrm{SO}_{7}(\mathbb{C})$
with $\ker\pi=\{1\}\times\{\pm1\}$. Now let us denote $\mathrm{SL}_{2}(\mathbb{C})\times\mathrm{SO}_{7}(\mathbb{C})$
by $K$. We know that $\ker\pi$ acts trivially on $\mathfrak{g}_{\bar{0}}$,
thus we obtain an induced action of $K$ on $\mathfrak{g}_{\bar{0}}$.
Let $\mathfrak{z}(\mathfrak{g}^{e})_{\bar{0}}=\mathfrak{z}(\mathfrak{g}^{e})\cap\mathfrak{g}_{\bar{0}}$
and $\mathfrak{z}(\mathfrak{g}^{e})_{\bar{1}}=\mathfrak{z}(\mathfrak{g}^{e})\cap\mathfrak{g}_{\bar{1}}$,
note that 
\[
\left(\mathfrak{z}(\mathfrak{g}^{e})\right)^{G^{e}}=\left(\mathfrak{z}(\mathfrak{g}^{e})_{\bar{0}}\oplus\mathfrak{z}(\mathfrak{g}^{e})_{\bar{1}}\right)^{G^{e}}=\left(\mathfrak{z}(\mathfrak{g}_{\bar{0}}^{e})\oplus\mathfrak{z}(\mathfrak{g}_{\bar{1}}^{e})\right)^{G^{e}}=\left(\mathfrak{z}(\mathfrak{g}_{\bar{0}}^{e})\right)^{G^{e}}\oplus\left(\mathfrak{z}(\mathfrak{g}_{\bar{1}}^{e})\right)^{G^{e}}.
\]
Furthermore, we have that $(\mathfrak{z}(\mathfrak{g}_{\bar{0}}^{e}))^{G^{e}}=(\mathfrak{z}(\mathfrak{g}_{\bar{0}}^{e}))^{K^{e}}$.
Thus when $\mathfrak{z}(\mathfrak{g}_{\bar{1}}^{e})=0$, it suffices
to look at $(\mathfrak{z}(\mathfrak{g}_{\bar{0}}^{e}))^{K^{e}}$.
It is obvious that $e\subseteq(\mathfrak{z}(\mathfrak{g}_{\bar{0}}^{e}))^{K^{e}}$.
Since $\mathrm{SL}_{2}(\mathbb{C})$ is connected, we have that $K^{e}$
is the semidirect product of the subgroup $C^{e}$ and the normal
subgroup $R^{e}$. Furthermore, we recall that $K^{e}/(K^{e})^{\circ}\cong C^{e}/(C^{e})^{\circ}$.
Now we have that $C^{e}\cong\left(\mathrm{O}_{1}(\mathbb{C})\times\mathrm{O}_{2}(\mathbb{C})\right)\cap\mathrm{SO}_{7}(\mathbb{C})$
where $\mathrm{O}_{1}(\mathbb{C})$ (resp. $\mathrm{O}_{2}(\mathbb{C})$)
has the connected component $\mathrm{SO}_{1}(\mathbb{C})$ (resp.
$\mathrm{SO}_{2}(\mathbb{C})$). Let us consider the element $g\in K^{e}/(K^{e})^{\circ}$
given by 
\[
g=\begin{pmatrix}0 & 1 & 0 & 0 & 0 & 0 & 0\\
1 & 0 & 0 & 0 & 0 & 0 & 0\\
0 & 0 & 0 & 0 & 1 & 0 & 0\\
0 & 0 & 0 & -1 & 0 & 0 & 0\\
0 & 0 & 1 & 0 & 0 & 0 & 0\\
0 & 0 & 0 & 0 & 0 & 0 & 1\\
0 & 0 & 0 & 0 & 0 & 1 & 0
\end{pmatrix}.
\]
We calculate that $g\cdot e_{(3^{2},1)}=ge_{(3^{2},1)}g^{-1}=e$ and
$g\cdot R_{e_{1},e_{2}}=gR_{e_{1},e_{2}}g^{-1}=-R_{e_{1},e_{2}}$.
Hence, we have that $\left(\mathfrak{z}(\mathfrak{g}^{e})\right)^{K^{e}}\subseteq\left(\mathfrak{z}(\mathfrak{g}^{e})\right)^{g}=\langle e\rangle$
and we deduce that $\left(\mathfrak{z}(\mathfrak{g}^{e})\right)^{K^{e}}=\langle e\rangle$
. Therefore, $\dim\left(\mathfrak{z}(\mathfrak{g}^{e})\right)^{G^{e}}=n_{2}(\varDelta)=1$
and $\dim\left(\mathfrak{z}(\mathfrak{g}^{e})\right)^{G^{e}}-\dim\left(\mathfrak{z}(\mathfrak{g}_{0}^{e_{0}})\right)^{G_{0}^{e_{0}}}=n_{2}(\varDelta)=1$.
When $e=e_{(7)},e_{(5,1^{2})}$, using similar arguments we obtain
that $\left(\mathfrak{z}(\mathfrak{g}^{e})\right)^{G^{e}}=\langle e,R_{e_{1},e_{2}}\rangle$.

When $e=E+e_{(7)}$, recall that $\mathfrak{z}(\mathfrak{g}^{e})=\langle E+e_{(7)},v_{1}\otimes e_{1}e_{2}e_{3},R_{e_{1},e_{2}}\rangle.$
We know that $G^{e}=\left(\{\pm1\}\ltimes R^{E}\right)\times\mathrm{Spin}_{7}(\mathbb{C})^{e_{(7)}}$
where $R^{E}$ is a connected normal subgroup of $G^{e}$. Now we
take $g=-1\in\mathrm{SL}_{2}(\mathbb{C})$ such that $g\in G^{e}$
and $g\notin(G^{e})^{\circ}$. We know that $g$ acts trivially on
$\mathfrak{z}(\mathfrak{g}^{e})_{\bar{0}}$, thus $e,R_{e_{1},e_{2}}\in\left(\mathfrak{z}(\mathfrak{g}^{e})\right)^{g}$.
However, $v_{1}\otimes e_{1}e_{2}e_{3}\notin\mathfrak{z}(\mathfrak{g}^{e})^{g}$
since the action of $g$ on $v_{1}\otimes e_{1}e_{2}e_{3}$ sends
it to $-v_{1}\otimes e_{1}e_{2}e_{3}$. Hence, we have that $\left(\mathfrak{z}(\mathfrak{g}^{e})\right)^{G^{e}}\subseteq\left(\mathfrak{z}(\mathfrak{g}^{e})\right)^{g}=\langle e,R_{e_{1},e_{2}}\rangle$.
Therefore, we deduce that $\left(\mathfrak{z}(\mathfrak{g}^{e})\right)^{G^{e}}=\langle e,R_{e_{1},e_{2}}\rangle$
and $\dim\left(\mathfrak{z}(\mathfrak{g}^{e})\right)^{G^{e}}-\dim\left(\mathfrak{z}(\mathfrak{g}_{0}^{e_{0}})\right)^{G_{0}^{e_{0}}}=n_{2}(\varDelta)=1$.
Using similar arguments we obtain that $\left(\mathfrak{z}(\mathfrak{g}^{e})\right)^{G^{e}}=\langle e,R_{e_{1},e_{2}}\rangle$
for $e=E+e_{(5,1^{2})}$ and $\left(\mathfrak{z}(\mathfrak{g}^{e})\right)^{G^{e}}=\langle e\rangle$
for $e=E+e_{(3^{2},1)}$. The above argument completes the proof of
Theorems 1 and 2 for $F(4)$.

By combining results in Table \ref{tab:F(4)}, we have that $\dim\left(\mathfrak{z}(\mathfrak{g}^{e})\right)^{G^{e}}=\left\lceil \frac{1}{2}\sum_{i=1}^{4}a_{i}\right\rceil +\varepsilon$
where $\varepsilon=-1$ for $e=E+e_{(7)}$ and $\varepsilon=0$ for
all other cases. This proves the statement of Theorem 3 for $F(4)$.
\end{singlespace}

\noindent \address{School of Mathematics, University of Birmingham, Birmingham, B152TT, UK}
\end{document}